\documentclass{amsart}

\usepackage{amsmath}
\usepackage{amssymb}
\usepackage{amsthm}
\usepackage[all]{xy}
\usepackage{graphicx}

\theoremstyle{plain}

\newtheorem{thm}{Theorem}[section]
\newtheorem{lem}[thm]{Lemma}
\newtheorem{prop}[thm]{Proposition}

\theoremstyle{definition}

\newtheorem{defn}[thm]{Definition}
\newtheorem{rem}[thm]{Remark}
\newtheorem{ex}[thm]{Example}

\newcommand{\R}{\mathbb{R}}
\newcommand{\Z}{\mathbb{Z}}

\newcommand{\D}{\mathcal{D}}

\newcommand{\abs}[1]{\lvert {#1} \rvert}

\newcommand{\calH}{\mathcal{H}}

\newcommand{\emb}[2]{\mathcal{K}_{#1 ,#2}}

\newcommand{\gr}{\mathrm{gr}}

\newcommand{\I}{\mathcal{I}}
\newcommand{\id}{\mathrm{id}}

\newcommand{\Int}{\mathrm{Int}\,}

\newcommand{\ord}{\mathrm{ord}\,}

\newcommand{\pair}[2]{\langle {#1},\, {#2}\rangle}

\numberwithin{equation}{section}
\numberwithin{figure}{section}
\numberwithin{table}{section}

\title[Configuration space integrals and the Haefliger invariant]
{Configuration space integrals for embedding spaces and the Haefliger invariant}
\author{Keiichi Sakai}
\date{\today}
\address{Graduate School of Mathematical Sciences, University of Tokyo}
\email{ksakai@ms.u-tokyo.ac.jp}
\urladdr{http://www.ms.u-tokyo.ac.jp/~ksakai/index.html}

\subjclass[2000]{Primary~58D10; Secondary~57Q45, 81Q30}
\keywords{The space of long $j$-knots; Graph complex; Configuration space integral; Haefliger invariant}

\begin{document}

\begin{abstract}
Let $\emb{n}{j}$ be the space of long $j$-knots in $\R^n$.
In this paper we introduce a graph complex $\D^*$ and a linear map $I:\D^* \to \Omega^*_{DR}(\emb{n}{j})$ via
configuration space integral, and prove that
(1) when both $n>j\ge 3$ are odd, $I$ is a cochain map if restricted to graphs with at most one loop component,
(2) when $n-j \ge 2$ is even, $I$ is a cochain map if restricted to tree graphs, and
(3) when $n-j \ge 3$ is odd, $I$ added a correction term produces a $(2n-3j-3)$-cocycle of $\emb{n}{j}$ which gives a new
formulation of the Haefliger invariant when $n=6k$, $j=4k-1$ for some $k$.
\end{abstract}

\maketitle

\section{Introduction}\label{sec_intro}

The aim of this research is to study the topology of the space of long knots.

\begin{defn}\label{def_knotspace}
A {\em long $j$-knot} in $\R^n$ is an embedding $f: \R^j \hookrightarrow \R^n$ which is standard at infinity:
\[
 f(x) = (x,0) \in \R^j \times \{ 0 \}^{n-j} , \quad x \not\in [-1,1]^j .
\]
Denote by $\emb{n}{j}$ the space of long $j$-knots in $\R^n$, equipped with $C^{\infty}$-topology.
\end{defn}

In this and the forthcoming papers \cite{KWatanabe08} we will develop the methods to study $\emb{n}{j}$ originated in
perturbative Chern-Simons theory.
This method was used by various authors \cite{AltschulerFreidel97, BottTaubes94, Kohno94} to give some integral
expressions of the finite type invariants for knots in $\R^3$.
Cattaneo, Cotta-Ramusino and Longoni \cite{CCL02} generalized their constructions to define a cochain map from a
certain graph complex to the de Rham complex of the space of (long) $1$-knots in $\R^n$, $n>3$, via fiber-integrations
over configuration spaces.
Not only the trivalent graphs correspond to finite type invariants, but non-trivalent graphs can also work in their
framework.
Indeed, some non-trivalent graph cocycles produce non-trivial cohomology classes \cite{Longoni04, K07}.

Another generalization was done by Rossi \cite{Rossi_thesis} and Cattaneo-Rossi \cite{CattaneoRossi05}, who proved the
invariance of (order two) Bott invariant and order three invariant for long $m$-knots in $\R^{m+2}$ for $m \ge 2$.
Following their work, Watanabe \cite{Watanabe07} proved that there is one finite type invariant for long ``ribbon''
$m$-knots \cite{HabiroKanenobuShima99} in $\R^{m+2}$, $m \ge 3$ odd, at each even order.
These invariants also come from trivalent graphs via the perturbative method.

In this paper we introduce graph complexes in the same manner as \cite{CCL02} generated by more general graphs
than those in \cite{CattaneoRossi05, Rossi_thesis, Watanabe07}.
We prove that some graph cocycles produce cohomology classes of $\emb{n}{j}$ via configuration space integrals.

\begin{thm}\label{thm_main2}
There exist graph complexes\footnote
{
The combinatorial meaning of the number $k$ will be specified in Definition \ref{definition_space_graphs}.
}
 $\D^{k,*}_g$, $k \ge 1$ and $g\ge 0$, spanned by graphs of first Betti number $g$ (after its
`small loops' are removed) and linear maps $I : \D^{k,l}_g \to \Omega^{k(n-j-2)+(g-1)(j-1)+l}_{DR}(\emb{n}{j})$ given by
configuration space integral.
The map $I$ is a cochain map when $n-j \ge 2$ is even and $g=0$, or both $n>j\ge 3$ are odd and $g=1$.
\end{thm}

Unfortunately it is not known how `big' the graph cohomology $H^* (\D^*_0 )$ is.
But in \cite{KWatanabe08} we will prove that $H^* (\D^*_1 )$ is not trivial, and the map $I$ produces non-trivial
cohomology classes of $\emb{n}{j}$ and $\overline{\mathcal{K}}_{n,j}$ for various $n$ and $j$, which generalize the
`finite type invariants' for `long ribbon knots' \cite{HabiroKanenobuShima99, Watanabe07}.
Here $\overline{\mathcal{K}}_{n,j}$ is the space of long $j$-knots `modulo immersions' (see for example
\cite{KWatanabe08, Sinha04}).
When $n-j$ is odd, we have no idea to prove that $I$ is a cochain map, since we cannot ignore the contributions of some
`hidden faces' of the boundary of configuration spaces.
So we will need some correction terms (different from that given in Theorem \ref{thm_main1}).
See \cite{KWatanabe08}.

The graph complex introduced in \cite{CCL02} is conjectured to give all the real cohomology classes of $\emb{n}{1}$,
while it is not clear whether the graphs appeared in \cite{CattaneoRossi05, Rossi_thesis, Watanabe07} are enough for
$H^*_{DR}(\emb{n}{j})$.
But our graph complex $\D^* =\bigoplus_g \D^*_g$ contains more general graphs, so we might be able to conjecture that
our graph complex would describe whole $H^*_{DR}(\emb{n}{j})$.
The next Theorem \ref{thm_main1} confirms the conjecture; via perturbative method we can give a new formulation of the
{\em Haefliger invariant} \cite{Haefliger62,Haefliger66}.

\begin{thm}\label{thm_main1}
Suppose $n-j \ge 3$ is odd and $2n-3j-3 \ge 0$.
Then a graph cocycle $H \in \D^{2,0}_0$ produces a non-trivial cohomology class
$\calH :=[I(H)+c]\in H^{2n-3j-3}_{DR}(\emb{n}{j})$, where $c$ is some correction term.
Moreover when $n=6k$ and $j=4k-1$, $\calH \in H^0_{DR}(\emb{6k}{4k-1})$ is nothing but the Haefliger invariant
(up to sign).
\end{thm}

The correction term $c$ is added to kill some contribution of the `anomalous face' of a compactified configuration space.
This contribution is an obstruction for $I$ to be a cochain map.
See \S \ref{sec_Haefliger} and \S \ref{sec_vanish} for details.

The Haefliger invariant, an isotopy invariant for (long) $(4k-1)$-knots in $6k$-space, was originally defined by
using  a $4k$-manifold bounded by the knot in $(6k+1)$-space.
Instead of such additional data, we use the generators of cohomology of configuration space.

This paper is organized as follows.
We define the graph complexes in \S \ref{section_graph} and describe the integration map $I$ in detail in \S
\ref{sec_conf}.
Several vanishing results are proved in \S \ref{sec_vanish} to prove Theorem \ref{thm_main2}.
The class $\calH \in H^{2n-3j-3}_{DR}(\emb{n}{j})$ is studied in \S \ref{sec_Haefliger}.
To prove $d\calH =0$ we need some of the results in \S \ref{sec_vanish}, which hold even if $n-j$ is odd.

\subsection*{Acknowledgment}
The author expresses his great appreciation to Professor Toshitake Kohno for his encouragement, to Tadayuki Watanabe
for many useful suggestions and comments, and to Masamichi Takase for his idea to improve Theorem \ref{thm_main1}.
The author is partially supported by the Grant-in-Aid for Young Scientists (B), MEXT, Japan, by The Sumitomo Foundation,
and by The Iwanami Fujukai Foundation.

\section{Graph complexes}\label{section_graph}

This and the next sections provide preliminaries for not only the subsequent sections but the forthcoming paper
\cite{KWatanabe08}.
In this section we introduce more general graphs than those in \cite{CattaneoRossi05, Rossi_thesis, Watanabe07},
including `degenerate graphs.'

\subsection{Graphs}\label{subsection_graph}

The {\em graphs} appearing here have two types of vertices.
One is the {\it external} vertex (or shortly {\it e-vertex}), which is depicted by $\circ$, while the other is the
{\it internal} vertex (shortly {\it i-vertex}), depicted by $\bullet$.

Similarly there are two types of edges.
One is $\theta${\em -edge}, which is depicted by a dotted line, and another is $\eta${\it -edge}, depicted by a
solid line.
We suppose that all the endpoints of $\eta$-edges are i-vertices.

A $\theta$-edge $e$ can form a {\it small loop} at an i-vertex $p$, that is, $e$ may have exactly one i-vertex $p$ as
its endpoint.
A single $\eta$-edge is not allowed to be a small loop, but an $\eta$-edge together with a $\theta$-edge can form
a loop at an i-vertex, called a {\em double loop}.
When we count the edges of a graph, we count a double loop twice, regarding it as consisting of an $\eta$-loop and a
$\theta$-loop.
But a double loop raises the first Betti number of the graph by one.

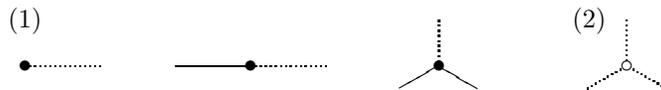
\begin{figure}[htb]
\[
 \begin{xy}
  (0,6)*{(1)},(75,6)*{(2)},
  (0,0)*{\bullet}, (30,0)*{\bullet},(55,0)*{\bullet},(80,0)*{\circ}="A",
  {\ar@{.}(0,0);(10,0)},
  {\ar@{-}(20,0);(30,0)},{\ar@{.}(30,0);(40,0)},
  {\ar@{.}(55,6);(55,0)},{\ar@{-}(55,0);(49.8038,-3)},{\ar@{-}(55,0);(60.1961,-3)},
  {\ar@{.}(80,6);"A"},{\ar@{.}"A";(74.8038,-3)},{\ar@{.}"A";(85.1961,-3)}
 \end{xy}
\]
\caption{Admissible vertices}\label{fig_admissible}
\end{figure}

\begin{defn}\label{def_admissible}
A vertex $v$ of a graph is {\em admissible} if it is
\begin{enumerate}
\item an i-vertex of valence $\ge 1$ with at least one $\theta$-edge (possibly a loop) emanating from $v$, or
\item an e-vertex of valence $\ge 3$, with only $\theta$-edges (which are not loops) emanating from $v$
\end{enumerate}
(see Figure \ref{fig_admissible}).
A graph is called {\em admissible} if all its vertices are admissible.
\end{defn}

\begin{figure}[htb]
\[
 \begin{xy}
                        (15,20)*{\bullet}="A", (30,20)*{\bullet}="B",
  (0,14)*{\bullet}="C",                                               (45,14)*{\circ}="D", (55,14)*{\bullet}="E",
                        (15, 8)*{\bullet}="F", (30, 8)*{\bullet}="G",
                        (15, 0)*{\bullet}="H", (30, 0)*{\bullet}="I",
  {\ar@{.>}"A";"C"^<{1}_>{4}}, {\ar "A";"B"_>{5}},{\ar@{.}@(ru,lu)"B";"B"}, {\ar@{.>}"B";"D"^>{9}},
  {\ar@{.>}"D";"E"^>{7}},
  {\ar@{.>}"C";"F"^>{8}},{\ar "G";"F"_<{3}},{\ar@{.>}"G";"D"},
  {\ar@{.>}"F";"H"^>{6}},{\ar@{.>}"G";"I"^>{2}}
 \end{xy}
\]
\caption{An example of an admissible graph for odd $n$ and $j$}\label{example_1}
\end{figure}
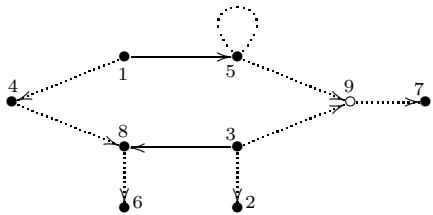

Figure \ref{example_1} shows an example of an admissible graph.
There might be an i-vertex which is adjacent to more than one $\theta$-edges, an e-vertex of valency $\ge 4$, and so on.
Such vertices did not appear in \cite{CattaneoRossi05,Rossi_thesis,Watanabe07}.

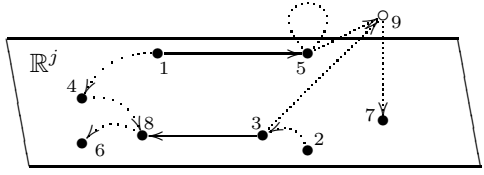
\begin{figure}[htb]
\[
 \begin{xy}
  {\ar@{-}(0,20);(3,3)},{\ar@{-}(3,3);(63,3)},{\ar@{-}(0,20);(60,20)},{\ar@{-}(60,20);(63,3)},
  (5,17)*{\R^j},
  (50,23)*{\circ}="A",
  (20,18)*{\bullet}="B", (40,18)*{\bullet}="C",
  (10,12)*{\bullet}="D", (18,7)*{\bullet}="E", (34,7)*{\bullet}="F",(50,9)*{\bullet}="G",
  (10,6)*{\bullet}="H",  (40,5)*{\bullet}="I",
  {\ar@{.>}@/_/"B";"D"_>{4}},{\ar "B";"C"_<{1}_>{5}},{\ar@{.}@(ru,lu)"C";"C"}, {\ar@{.>}"C";"A"},
  {\ar@{.>}@/^/"D";"E"},{\ar "F";"E"_<{3}_>{8}}, {\ar@{.>}"F";"A"},{\ar@{.>}@/_/"I";"F"_<{2}},
  {\ar@{.>}"A";"G"^<{9}_>{7}},{\ar@{.>}@/_/"E";"H"^>{6}}
 \end{xy}
\]
\caption{The graph from Figure \ref{example_1} in $\R^n$}\label{example_2}
\end{figure}

In \S \ref{sec_conf} the graphs will be regarded as in $\R^n$ (Figure \ref{example_2}); the set of vertices will become
a configuration of points in $\R^n$, where the i-vertices are on $\R^j$ (or a long $j$-knot) embedded in $\R^n$.
Each edge will correspond to a `direction map' determined by its endpoints, or to the volume form of sphere pulled back
by the direction map.
The $\eta$-edges correspond to directions in $\R^j$, while $\theta$-edges to those in $\R^n$.

Afterward we will need to fix the orientations of configuration spaces and the signs of the volume forms.
For these purposes we `decorate' the graphs as follows (see Figure \ref{fig_decoration}).
\begin{enumerate}
\item
 The i-vertices are labeled by $1,\dots ,s$ and the e-vertices are labeled by $s+1,\dots ,s+t$ for some suitable $s$
 and $t$.
\item
 When both $n$ and $j$ are odd, all the edges are oriented, all the loops are ordered and each loop is given the sign
 $\pm 1$.
\item
 When $n$ is odd and $j$ is even, all the $\theta$-edges are oriented, while all the $\eta$-edges are labeled.
\item
 When $n$ is even and $j$ is odd, all the $\theta$-edges (including those in double loops) are labeled, while all the
 $\eta$-edges are oriented.
 The double loops are given other labels $1,2,\dots$ than those for all the $\theta$-edges.
 A double loop is given a sign $\pm 1$.
\item
 When both $n$ and $j$ are even, all the edges (including those in double loops) are labeled.
 The small / double loops are given other labels.
\end{enumerate}
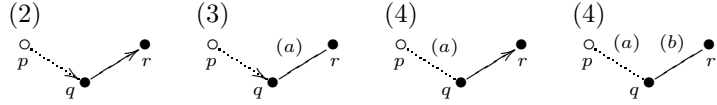
\begin{figure}[htb]
\[
 \begin{xy}
  (0,9)*{(2)},(0,5)*{\circ}="A", (8,0)*{\bullet}="B", (16,5)*{\bullet}="C",
  {\ar@{.>}"A";"B"_<{p}_>{q}},{\ar "B";"C"_>{r}},
  (25,9)*{(3)},(25,5)*{\circ}="D", (33,0)*{\bullet}="E", (41,5)*{\bullet}="F",
  {\ar@{.>}"D";"E"_<{p}_>{q}},{\ar@{-} "E";"F"^(.5){(a)}_>{r}},
  (50,9)*{(4)},(50,5)*{\circ}="G", (58,0)*{\bullet}="H", (66,5)*{\bullet}="I",
  {\ar@{.}"G";"H"_<{p}^(.5){(a)}_>{q}},{\ar "H";"I"_>{r}},
  (75,9)*{(4)},(75,5)*{\circ}="J", (83,0)*{\bullet}="K", (91,5)*{\bullet}="L",
  {\ar@{.}"J";"K"_<{p}^(.4){(a)}_>{q}},{\ar@{-} "K";"L"^(.6){(b)}_>{r}},
 \end{xy}
\]
\caption{Decorations of graphs}\label{fig_decoration}
\end{figure}

\subsection{Space of graphs}\label{subsection_space_of_graph}

Below we assume that all the graphs are admissible and decorated unless otherwise stated.

\begin{defn}\label{definition_space_graphs}
For a graph $\Gamma$, define
\begin{align*}
 \ord \Gamma &:= \sharp \{ \theta \text{-edges of } \Gamma \} - \sharp \{ \text{e-vertices of } \Gamma \} , \\
 \deg \Gamma &:= 2 \sharp \{ \theta \text{-edges of } \Gamma \} - 3 \sharp \{ \text{e-vertices of } \Gamma \}
 - \sharp \{ \text{i-vertices of } \Gamma \}
\end{align*}
(for example, $\ord \Gamma =7$ and $\deg \Gamma =5$ for the graph $\Gamma$ in Figure \ref{example_1}; see below for
more explanations).
We denote by $\D^{k,l}$ the vector space spanned by admissible decorated graphs with $\ord =k$, $\deg =l$ modulo the
subspace spanned by
\[
 \Gamma' - (-1)^{j\text{sign}\, \sigma +n\text{sign}\, \tau +a+b+\text{sign}\, \rho}\Gamma
 \quad \text{and} \quad \Gamma'' ,
\]
where $\Gamma'$ is obtained from $\Gamma$ by
\begin{itemize}
\item
 permuting the labels of i- and e-vertices of $\Gamma$ by $\sigma \in \mathfrak{S}_s$ and $\tau \in \mathfrak{S}_t$
 respectively ($s$ and $t$ are the numbers of i- and e-vertices of $\Gamma$ respectively),
\item
 reversing $a$ oriented edges of $\Gamma$, and
\item
 switching $b$ signs and permuting the labels of small / double loops by $\rho$,
\end{itemize}
and $\Gamma''$ is a graph with `multiple $\eta$- (or $\theta$-) edges,' i.e., there are two vertices $p,q$ which are
joined by two or more $\eta$- (or $\theta$-) edges (we allow $p,q$ joined by two edges, one $\eta$-edge and one
$\theta$-edge).
Moreover we introduce one more relation; $\Gamma \sim 0$ in $\D^*$ if $n-j$ is odd and $\Gamma$ is a graph with
at least one small loop, or if $n$ is odd and $\Gamma$ has a double loop.
\end{defn}

\begin{rem}\label{rem_graph}
The sign $j\text{sign}\, \sigma +n\text{sign}\, \tau$ will correspond to the orientation sign of the configuration
space.
In \S \ref{sec_conf} we will associate the volume forms of spheres of even dimensions with the oriented edges.
Reversing an oriented edge corresponds to pull-back via the antipodal map, hence yields a sign $-1$.
Unoriented edges correspond to differential forms of odd degrees, so permuting the labels yields a sign.
\end{rem}

The meaning of $\ord \Gamma$ is as follows.
In Figure \ref{example_2}, if we contract $\R^j$ together with $\eta$-edges regarded as in $\R^j$, then we obtain a
one dimensional CW complex whose edges are $\theta$-edges. Its first Betti number is equal to $\ord \Gamma$.

To explain the meaning of $\deg \Gamma$, we need some terminologies.

\begin{defn}\label{non_degenerte_graphs}
An admissible vertex is said to be {\it non-degenerate} if it is
\begin{itemize}
\item an i-vertex with exactly one $\theta$-edge (and possibly many $\eta$-edges) emanating from it, or
\item a tri-valent e-vertex.
\end{itemize}
All other vertices are said to be {\it degenerate}.
\end{defn}

For example, all the vertices in Figure \ref{fig_admissible}, and the vertices $1$, $2$, $6$, $7$ and $9$ of the graph
in Figure \ref{example_1} are non-degenerate.

\begin{rem}
It can be easily understood that ``non-degenerate vertex'' is the same notion as ``trivalent vertex'' in \cite{CCL02}.
All the vertices of the graphs appeared in \cite{CattaneoRossi05, Rossi_thesis, Watanabe07} are non-degenerate.
\end{rem}

\begin{lem}
We have $\deg \Gamma \ge 0$ for any admissible graph $\Gamma$, and $\deg \Gamma =0$ if and only if all the vertices
of $\Gamma$ are non-degenerate.
\end{lem}

\begin{proof}
This Lemma is obvious by the definition of admissible vertices;
at least one $\theta$-edge emanates from any i-vertex of a graph and at least three $\theta$-edges emanate from
any i-vertex.
This implies
\[
 2 \sharp \{ \theta \text{-edges}\} \ge 3 \sharp \{ \text{e-vertices}\} + \sharp \{ \text{i-vertices}\}
\]
and the equality holds if and only if exactly one (resp.\ three) $\theta$-edge emanates from any i-vertices
(resp.\ e-vertices), that is, all the vertices are non-degenerate.
\end{proof}

\begin{rem}\label{number_of_vertices}
Let $\Gamma$ be an admissible graph of $\ord \Gamma =k$, $\deg \Gamma =l$.
Then $\Gamma$ has $2k-l$ vertices; it is a direct consequence of the definition.
In particular, if $\Gamma$ is non-degenerate ($l=0$), then the number of all vertices is $2k$.
In \cite{Watanabe07} the half of the number of the vertices of a (non-degenerate) graph is called its `degree.'
Thus our terminology `$\ord$' is a generalization of the `degree' in \cite{Watanabe07}.
\end{rem}

\subsection{Coboundary operation}\label{subsection_coboundary}

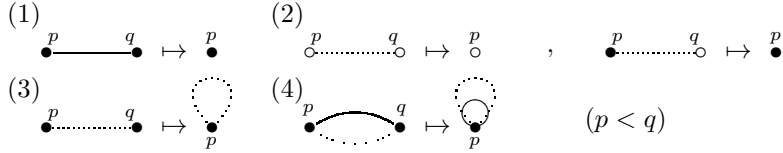
\begin{figure}[htb]
\[
 \begin{xy}
  (0,20)*{(1)}, (3,15)*{\bullet}="A", (15,15)*{\bullet}="B", {\ar@{-}"A";"B"^<{p}^>{q}},(20,15)*{\mapsto},
  (25,15)*{\bullet},(25,17)*{{\sb p}},
  (35,20)*{(2)}, (38,15)*{\circ}="C", (50,15)*{\circ}="D", {\ar@{.}"C";"D"^<{p}^>{q}},(55,15)*{\mapsto},
  (60,15)*{\circ},(60,17)*{{\sb p}}, (70,15)*{,},
  (78,15)*{\bullet}="E", (90,15)*{\circ}="F", {\ar@{.}"E";"F"^<{p}^>{q}},(95,15)*{\mapsto},
  (100,15)*{\bullet},(100,17)*{{\sb p}},
  (0,10)*{(3)}, (3,5)*{\bullet}="G", (15,5)*{\bullet}="H", {\ar@{.}"G";"H"^<{p}^>{q}},(20,5)*{\mapsto},
  (25,5)*{\bullet},(25,3)*{{\sb p}},{\ar@{.}@(ru,lu)(25,5);(25,5)},
  (35,10)*{(4)}, (38,5)*{\bullet}="I", (50,5)*{\bullet}="J", {\ar@{.}@/_/"I";"J"}, {\ar@{-}@/^/"I";"J"^<{p}^>{q}},
  (55,5)*{\mapsto}, (60,5)*{\bullet},(60,3)*{{\sb p}},{\ar@{.}@(ru,lu)(60,5);(60,5)},
  (60,7)*+[Fo]{\phantom{a}},
  (80,6)*{(p<q)}
  \end{xy}
\]
\caption{Contractions of edges}\label{edge_contraction}
\end{figure}

\begin{defn}
Let $\Gamma$ be a graph and $e = \overrightarrow{pq}$ its (possibly oriented) edge (but not a loop).
Define a new graph $\Gamma / e$ as follows (see Figure \ref{edge_contraction}).

{\bf (1)} When $e$ is an $\eta$-edge (then endpoints $p,q$ are both internal), define $\Gamma / e$ by contracting $e$,
that is, identifying the endpoints $p,q$ of $e$ and removing the edge $e$.
The decoration of $\Gamma /e$ is derived from that of $\Gamma$; the vertex of $\Gamma / e$ where the contraction
occurred is re-labeled by $\min \{ p,q \}$, and all the labels of the vertices of $\Gamma$ bigger than $\max \{ p,q\}$
are decreased by one.
The labels of the other vertices remain unchanged.
When $j$ is even and $e$ is the $i$-th edge, then the labels of other edges bigger than $i$ is decreased by one.

{\bf (2)} When $e$ is a $\theta$-edge and at least one of $p,q$ is an e-vertex, then $\Gamma / e$ is defined in the same
way as above.
If both $p,q$ are external, then the vertex where the contraction occurred is also external.
If one of $p,q$ is internal, then the resulting vertex is internal.

{\bf (3)} When $e$ is a $\theta$-edge with both $p,q$ being i-vertices, then $\Gamma / e$ is obtained from $\Gamma$ by
identifying the vertices $p$ and $q$, but not removing the edge $e$.
The edge $e$ becomes a small loop at the i-vertex $\min \{ p,q\}$.
The labeling of $\Gamma / e$ is determined similarly as above.
When $n$ and $j$ are odd, its sign is $+1$ (resp.\ $-1$) if $p<q$ (resp.\ $p>q$).
This small loop is labeled by $a$ if $\Gamma$ has $(a-1)$ small loops.

{\bf (4)} When $e$ is the $\eta$-edge of the multiple edges joining two i-vertices, then $\Gamma / e$ is obtained from
$\Gamma$ by identifying the vertices $p$ and $q$, and attaching a double loop at $p$.
This double loop is labeled by $a$ if $\Gamma$ has $(a-1)$ loops.
When $j$ is odd, the sign $\pm 1$ is given similarly as in (3).
We do not define $\Gamma /e$ for a $\theta$-edge of the multiple edges.
\end{defn}

We should notice that a graph $\Gamma /e$ of type (3) in Figure \ref{edge_contraction} is ruled out in $\D^*$
when $n-j$ is odd, and similarly a graph $\Gamma /e$ of type (4) is ruled out when $n$ is odd (see Definition
\ref{definition_space_graphs}).

We would like to define the operator $\delta$ by
\[
 \delta \Gamma = \sum_{e \in E(\Gamma ) \setminus \{ \text{loops} \}} (-1)^{\tau (e)} \Gamma / e,
\]
by giving some suitable signs $\tau (e)$, where $E(\Gamma )$ is the set of edges of $\Gamma$.

\begin{prop}\label{def_signs}
If we define the signs $\tau (e)$ as in (1)-(5) below, then the operator $\delta$ is well-defined and determines
a coboundary operation
\[
 \delta : \D^{k,l} \longrightarrow \D^{k,l+1},
\]
that is, $\delta \circ \delta =0$.
Thus $\{ \D^{k,*}, \delta \}$ is a cochain complex for any $k$.

{\bf (1)} Let both $n$ and $j$ be odd.
For any oriented edge $e=\overrightarrow{pq}$, define $\tau (e)$ by
\begin{equation}\label{eq_sign_odd_odd}
 \tau (e) :=
 \begin{cases}
  q   & p<q, \\
  p+1 & p>q.
 \end{cases}
\end{equation}

{\bf (2)} Let both $n$ and $j$ be even.
For the $i$-th edge $e$, define $\tau (e)$ by
\[
 \tau (e) :=
 \begin{cases}
  i   & e \text{ is of type (1), (2) in Figure \ref{edge_contraction}}, \\
  u+1 & e \text{ is of type (3) in Figure \ref{edge_contraction}},
 \end{cases}
\]
where $u$ is the number of small / double loops of $\Gamma$.

{\bf (3)} Let $n$ be even and $j$ be odd.
For any oriented $\eta$-edge $e=\overrightarrow{pq}$, define $\tau (e)$ by \eqref{eq_sign_odd_odd}.
For the $i$-th $\theta$-edge $e=pq$, $p<q$ with $q$ being an e-vertex, define $\tau (e):=i+s+1$, where $s$ is the number
of i-vertices of $\Gamma$.

{\bf (4)} Let $n$ be odd and $j$ be even.
For the $i$-th $\eta$-edge $e$, define $\tau (e):=i+t+1$ where $t$ is the number of e-vertices of $\Gamma$.
For any oriented $\theta$-edge $e=\overrightarrow{pq}$, define $\tau (e)$ by \eqref{eq_sign_odd_odd}.

{\bf (5)} Consider the case that $n$ is even, and $p$ and $q$ are i-vertices joined by `multiple edges,' one $\eta$-edge
and one $\theta$-edge.
If $j$ is odd and the $\eta$-edge is oriented from $p$ to $q$, then $\tau$ is given by \eqref{eq_sign_odd_odd}.
If $j$ is even, then $\tau =u+1$.
\end{prop}

\begin{rem}
The signs in Proposition \ref{def_signs} correspond to those of induced orientations of the boundary strata of
configuration spaces; see \S \ref{subsec_orientation}.
\end{rem}

The proof is completely similar to \cite[Theorem 4.2]{CCL02}; choose two edges $e_1$ and $e_2$, and contract them
in two different orders, then we obtain the same graph with opposite signs.
Notice that if $n-j$ is odd (resp.\ $n$ is odd) the case (3) (resp.\ (4)) in Figure \ref{edge_contraction} does not occur.

\section{Configuration space integral}\label{sec_conf}

\subsection{Configuration spaces}\label{subsec_conf}

For any space $M$, denote the space of configurations in $M$ by
\[
 C^o_k (M) := \{ (x_1 ,\dots ,x_k )\in M^k \, |\, x_p \ne x_q \text{ if } p \ne q \} .
\]
Let $\Gamma$ be a (possibly non-admissible) decorated graph with $s$ i-vertices, $t$ e-vertices and $u$ loops.
Define the configuration space associated with (vertices of) $\Gamma$ by
\[
 C^o_{\Gamma} := \left\{
 \genfrac{}{}{0pt}{}
  {(f;x_1 ,\dots ,x_s ; y_{s+1} ,\dots ,y_{s+t})}{\in \emb{n}{j} \times C^o_s (\R^j ) \times C^o_t (\R^n )}
 \right. \left|
 \genfrac{}{}{0pt}{}{f(x_p ) \ne y_q ,}{\forall p,q}
 \right\}  \times (S^{j-1})^u .
\]
We think of i-vertex $p$ (resp.\ e-vertex $q$) as corresponding to $x_p \in \R^j$ (resp.~$y_q \in \R^n$)
for all $1 \le p \le s$ (resp.\ $s+1 \le q \le s+t$).
The $S^{j-1}$-factors will be used to define a differential form $\omega_e$ for a loop $e$ (see below).
There is a projection
\[
 \pi_{\Gamma} : C^o_{\Gamma} \longrightarrow \emb{n}{j}.
\]
We will denote its fiber over $f$ by $C^o_{\Gamma}(f)$.

\subsection{Differential forms associated to graphs}\label{subsection_form}
Let $e=\overrightarrow{pq}$ be an (oriented) edge or a loop of $\Gamma$.
To $e$ we will assign a differential form $\omega_e \in \Omega^*_{DR}(C^o_{\Gamma})$ as follows.

First consider the case that $e$ is not a loop (thus $p \ne q$).
When $e$ is an $\eta$-edge (then $p,q \le s$), define the `direction map' $\varphi^{\eta}_e : C^o_{\Gamma} \to S^{j-1}$
by
\[
 \varphi^{\eta}_e (f,x,y,v) := \frac{x_q -x_p}{\abs{x_q -x_p}} \in S^{j-1}.
\]
When $e$ is a $\theta$-edge, define $\varphi^{\theta}_e : C^o_{\Gamma} \to S^{n-1}$ by
\[
 \varphi^{\theta}_e (f,x,y,v) := \frac{z_q -z_p}{\abs{z_q -z_p}} \in S^{n-1},
\]
where
\[
 z_p =
 \begin{cases}
  f(x_p ) & \text{if the vertex $p$ is internal (thus $p \le s$)}, \\
  y_p     & \text{if the vertex $p$ is external (thus $p>s$)}.
 \end{cases}
\]

\noindent
{\bf Notation.}
We denote by $vol_{S^{N-1}}$ the {\em (anti-)symmetric} volume form of $S^{N-1}$ ($N=j$ or $n$).
Namely, $vol_{S^{N-1}}$ is an $(N-1)$-form of $S^{N-1}$ with total integral one and
$\iota^* vol_{S^{N-1}} = (-1)^N vol_{S^{N-1}}$ for the antipodal map $\iota : S^{N-1} \to S^{N-1}$.

\smallskip

Define the differential forms $\omega_e \in \Omega^*_{DR}(C^o_{\Gamma})$ as
$\eta_e := (\varphi^{\eta}_e )^* vol_{S^{j-1}}$ or $\theta_e := (\varphi^{\theta}_e )^* vol_{S^{n-1}}$,
according to whether $e$ is an $\eta$-edge or a $\theta$-edge.
Notice that, if $e$ is not oriented, then the map $\varphi_e$ has ambiguity of signs, but in such a case the
corresponding volume form is of odd degree and is invariant under the antipodal map of spheres.
Hence the form $\omega_e$ is well defined.

When $n-j$ is even, we also assign differential forms to small loops.
For the $a$-th small loop $e$ with sign $\varepsilon$ (which is always $+1$ when $n$ is even) at the i-vertex $p$,
define $D_a : C^o_{\Gamma} \to S^{n-1}$ by
\[
 D_a (f,x,y,v) := \epsilon \cdot \frac{df_{x_p}(v_a )}{\abs{df_{x_p}(v_a )}},
\]
here $df_{x_p} : T_{x_p}\R^j \to T_{f(x_p )}\R^n$ is the derivation map.
The differential form associated with $e$ is $\omega_e :=D^*_a vol_{S^{n-1}} \in \Omega^{n-1}_{DR}(C^o_{\Gamma})$.
Such a form is not needed when $n-j$ is odd (see Definition \ref{definition_space_graphs}).

When $n$ is even, to the $a$-th double loop at the i-vertex $p$ with sign $\epsilon$ (which is $+1$ when $j$ is even),
we assign a map $\tilde{D}_a : C^o_{\Gamma} \to S^{j-1} \times S^{n-1}$ defined by
\[
 \tilde{D}_a (f,x,y,v):= \left( \epsilon v_a , \frac{df_{x_p}(v_a )}{\abs{df_{x_p}(v_a )}}\right)
\]
and a differential form
$\omega_e :=\tilde{D}_a^* (vol_{S^{j-1}} \times vol_{S^{n-1}})\in \Omega^{n+j-2}_{DR}(C^o_{\Gamma})$.

Define the differential form $\omega_{\Gamma}$ by
\[
 \omega_{\Gamma} := \bigwedge_{e \in E(\Gamma )}\omega_e \in \Omega^*_{DR} (C^o_{\Gamma}),
\]
here $\omega_e$'s for labeled edges must be ordered according to the labels, since by definition they are odd forms.
Since non-labeled edges correspond to even forms, we need not to mention the order of corresponding forms.

\subsection{Fiber integration and compactified configuration spaces}\label{subsection_integral}
We would like to define $I(\Gamma ) \in \Omega^*_{DR} (\emb{n}{j})$ by integrating $\omega_{\Gamma}$ along the fiber
of $\pi_{\Gamma} : C^o_{\Gamma} \to \emb{n}{j}$.
It is not clear whether such an integral converges, because the fiber of $\pi_{\Gamma}$ is not compact.
This difficulty has been resolved in \cite{AxelrodSinger94,BottTaubes94}; for any manifold $M$, we can construct a
compact manifold $C_k (M)$ with corners so that its interior is $C^o_k (M)$, by `blowing up' all the diagonals of $M^k$.
We can smoothly extend the direction maps like $\varphi$'s, the projections $C^o_k (M) \to C^o_{k-l}(M)$ ($l>0$)
and the evaluation map
\[
 ev : \emb{n}{j} \times C^o_k (\R^j ) \longrightarrow C^o_k (\R^n ), \quad (f;x_1 ,\dots ,x_k ) \longmapsto
 (f(x_1 ), \dots ,f(x_k ))
\]
onto $C_k (M)$.

Let $s$ and $t$ be the numbers of i- and e-vertices of $\Gamma$ respectively.
Define $C'_{\Gamma}$ by the pull-back square
\[
 \xymatrix{
  C'_{\Gamma} \ar[r] \ar[d] & C_{s+t}(\R^n ) \ar[d]^-{pr_s} \\
  \emb{n}{j} \times C_s (\R^j ) \ar[r]^-{ev} & C_s (\R^n )
 }
\]
where $pr_s$ is the first $s$ projection, and define $C_{\Gamma} := C'_{\Gamma} \times (S^{j-1})^u$ ($u$ is the number of
loops of $\Gamma$).
The differential form $\omega_{\Gamma}$ is defined on $C_{\Gamma}$ since the direction maps $\varphi$'s are well defined
on $C_{\Gamma}$.
The natural projection $\pi_{\Gamma} : C_{\Gamma} \to \emb{n}{j}$ is defined and is a fibration with compact fibers.
Thus the integration
\[
 I(\Gamma ) :=(\pi_{\Gamma})_* \omega_{\Gamma} \in \Omega^*_{DR}(\emb{n}{j})
\]
along the fiber is well defined.
The degree of $I(\Gamma )$ is given by
\begin{align*}
 &\deg \omega_{\Gamma} -\dim \text{fib}(\pi_{\Gamma}) \\
 &\ =(n-1)\sharp \{ \theta \} + (j-1)\sharp \{ \eta \} -j\sharp \{ \bullet \} - n\sharp \{ \circ \} -(j-1)u \\
 &\ =(n-3)(\sharp \{ \theta \} - \sharp \{ \circ \}) + (j-1)(\sharp \{ \eta \} - \sharp \{ \bullet \}-u)
  + 2 \sharp \{ \theta \} - 3\sharp \{ \circ \} -\sharp \{ \bullet \} \\
 &\ =(n-3)\ord \Gamma + (j-1)(\sharp \{ \eta \} - \sharp \{ \bullet \}-u) + \deg \Gamma .
\end{align*}
Suppose that a one-dimensional CW complex $\Gamma \setminus \{ \text{small loops}\}$ has $g$ loop components (that is,
the first Betti number of it is $g$).
Then
\[
 \sharp \{ \eta \} - \sharp \{ \bullet \} -u = -\sharp \{ \theta \} + \sharp \{ \circ \}+(g-1) = -\ord \Gamma +(g-1),
\]
and hence
\[
 \deg I(\Gamma ) = (n-j-2)\ord \Gamma +(g-1)(j-1) + \deg \Gamma .
\]

\begin{prop}[\cite{CCL02}]
The form $I(\Gamma )$ depends only on the equivalence class of $\Gamma \in \D^*$.
\end{prop}

\begin{proof}
This is because the space $\D^*$ is arranged so that the map $I$ is compatible with the permutations of coordinates
and the antipodal map of spheres (see \cite{CCL02} for details).
But there are two point we should stress here; one is that we have to temporarily gave the labels to loops when $j$ is
odd.
But the choices of the labeling do not change the form $I(\Gamma )$ since the permutation on the $(S^{j-1})^u$-factor
of the fiber $C_{\Gamma}(f)$ is an orientation preserving diffeomorphism and preserves the form $\omega_{\Gamma}$.

Another is that we ruled out the graphs with small loops when $n-j \ge 3$ is odd, and those with double loops when $n$ is
odd (see Definition \ref{definition_space_graphs}).
This is because the differential forms arising from such graphs are zero (see Lemmas \ref{lem_small_loop_vanish},
\ref{lem_double_loop_vanish}).
\end{proof}

\begin{lem}\label{lem_small_loop_vanish}
When $n-j$ is odd, then the differential form $I(\Gamma )$ is zero for any graph $\Gamma$ with at least one small loop.
\end{lem}

\begin{proof}
Let $e$ be a small loop of $\Gamma$.
Define a fiberwise involution $F:C_{\Gamma} \to C_{\Gamma}$ as the antipodal map on $S^{j-1}$ corresponding to $e$,
and the identity on the other factors.
The orientation sign of the involution on the fiber of $\pi_{\Gamma}$ is $(-1)^j$, while
$F^* \omega_{\Gamma}=(-1)^n \omega_{\Gamma}$ since
\begin{align*}
 F^* \omega_e &= \left( \frac{df}{\abs{df}} \circ \iota_{S^{j-1}} \right)^* vol_{S^{n-1}}
  = \left( \iota_{S^{n-1}} \circ \frac{df}{\abs{df}} \right)^* vol_{S^{n-1}} \\
 &= (-1)^n \left( \frac{df}{\abs{df}} \right)^* vol_{S^{n-1}} =(-1)^n \omega_e
\end{align*}
and other $\omega_e$'s do not change.
Thus
\[
 (\pi_{\Gamma})_* \omega_{\Gamma} = (\pi_{\Gamma} \circ F)_* F^* \omega_{\Gamma}
 = (-1)^{j+n}(\pi_{\Gamma})_* \omega_{\Gamma} = -(\pi_{\Gamma})_* \omega_{\Gamma}
\]
since $n+j$ is odd, and hence $(\pi_{\Gamma})_* \omega_{\Gamma}=0$.
\end{proof}

\begin{lem}\label{lem_double_loop_vanish}
When $n$ is odd, then the differential form $I(\Gamma )$ is zero for any graph $\Gamma$ with at least one double loop.
\end{lem}

\begin{proof}
The fiberwise involution $F$ in the proof of Lemma \ref{lem_small_loop_vanish} has an orientation sign $(-1)^j$ on each
fiber, but in the case here $F^* \omega_{\Gamma}=(-1)^{n+j}\omega_{\Gamma}$ since $F^*$ affects $vol_{S^{j-1}}$ and
$vol_{S^{n-1}}$ as the antipodal map.
Thus
\[
 (\pi_{\Gamma})_* \omega_{\Gamma} = (-1)^{n+2j}(\pi_{\Gamma})_* \omega_{\Gamma}
 = -(\pi_{\Gamma})_* \omega_{\Gamma}.
\]
\end{proof}

Thus we obtain linear maps
\[
 I : \D^{k,l}_g \longrightarrow \Omega^{(n-j-2)k+(g-1)(j-1)+l}_{DR}(\emb{n}{j}),
\]
where $\D^{k,l}_g$ is the subspace of $\D^{k,l}$ spanned by the graphs $\Gamma$ whose first Betti numbers are $g$
after its small loops (not double loops) are removed.
It can be easily seen that $\D^{k,*}_g$ forms a subcomplex of $\D^{k,*}$ for any $g$ (we regard the contraction (4)
in Figure \ref{edge_contraction} as preserving the first Betti number).

Here we restate the last half of Theorem \ref{thm_main2}.

\begin{thm}\label{cochain}
Suppose $n-j\ge 2$ is even and $j \ge 2$.
The integration map $I : \D^{k,*}_g \to \Omega^{(n-j-2)k-(j-1)+*}_{DR}(\emb{n}{j})$ is a cochain map if
(1) both $n>j \ge 2$ are even and $g=0$, or (2) both $n > j \ge 3$ are odd and $g=0,1$.
\end{thm}

Theorem \ref{cochain} is a direct consequence of Theorems \ref{3.1_first_half}, \ref{3.1_second_half}, which are
proved in similar ways as the results in \cite{CCL02, Rossi_thesis, Watanabe07}.
The key step is the generalized Stokes theorem (see \cite{CCL02});
\[
 dI(\Gamma ) = (\pi_{\Gamma})_* (d\omega_{\Gamma})
 +(-1)^{\deg \omega_{\Gamma }+1} (\pi^{\partial}_{\Gamma})_* \omega_{\Gamma} 
 = (-1)^{\deg \omega_{\Gamma }+1} (\pi^{\partial}_{\Gamma})_* \omega_{\Gamma}
\]
where $\pi^{\partial}_{\Gamma}$ is the restriction of $\pi_{\Gamma}$ onto the boundaries of fibers.
The second equality holds since $\omega_{\Gamma}$ is a product of closed forms.
So we need to study the boundaries of fibers to prove Theorem \ref{cochain}.
The proof will be given in \S \ref{sec_vanish}.

Here we state one more result, which concerns the choices of volume forms.

\begin{prop}[\cite{CCL02}]\label{prop_indep_vol}
Suppose $g$, $n$ and $j$ satisfy (1) or (2) in Theorem \ref{cochain} and $n-j>2$.
Let $v_0$ and $v_1$ be two (anti-)symmetric volume forms of $S^{n-1}$ with total integral one, and $I_0$, $I_1$ the
corresponding integration maps.
Then $I_1 (\Gamma ) - I_0 (\Gamma )$ is an exact form for any  graph cocycle $\Gamma = \sum a_i \Gamma_i$.
\end{prop}

\begin{proof}
Choose any $w' \in \Omega^{n-2}_{DR}(S^{n-1})$ such that $v_1 -v_0 =dw'$, and put $w=(w' +(-1)^n \iota^*_{S^{n-1}}w')/2$.
Then $dw=v_1 -v_0$ and $\iota^*_{S^{n-1}}w = (-1)^n w$.
Define
\[
 \tilde{v} := v_0 +d(tw) \in \Omega^{n-1}_{DR}(S^{n-1}\times [0,1]),
\]
where $t$ is the coordinate of $[0,1]$.
Then $\tilde{v}$ is a (anti-)symmetric closed form on $S^{n-1}\times [0,1]$.
Moreover the restriction of $\tilde{v}$ onto $S^{n-1}\times \{ \varepsilon \}$, $\varepsilon =0,1$, is $v_{\varepsilon}$.

Assigning to $\theta$-edges $e$ of $\Gamma_i$ the differential forms
\[
 \tilde{\theta}_e := (\varphi^{\theta}_e \times \id )^* \tilde{v} \in \Omega^{n-1}_{DR}(C_{\Gamma_i} \times [0,1])
\]
instead of $\theta_e$, and integrating their product along the fiber of $\pi_{\Gamma_i}\times \id_{[0,1]}$,
we obtain a differential form of $\emb{n}{j} \times [0,1]$, which we denote by
$\tilde{I}(\Gamma_i )$.
Its restriction onto $\emb{n}{j} \times \{ \varepsilon \}$, $\varepsilon =0,1$, is $I_{\varepsilon}(\Gamma_i )$.

In Lemma \ref{tilde_I_closed} we will prove that $\tilde{I} (\Gamma ) \in \Omega^*_{DR}(\emb{n}{j} \times [0,1])$
is a closed form if $n-j >2$.
This completes the proof; if we denote the first projection by $p: \emb{n}{j} \times [0,1] \to \emb{n}{j}$, then
by the generalized Stokes theorem
\[
 d_{\emb{n}{j}}p_* \tilde{I}(\Gamma )
 = p_* d_{\emb{n}{j} \times [0,1]} \tilde{I}(\Gamma ) \pm p^{\partial}_* \tilde{I}(\Gamma ) 
 = \pm (\tilde{I}(\Gamma )|_{\emb{n}{j}\times \{1 \}} - \tilde{I}(\Gamma )|_{\emb{n}{j}\times \{0 \}})
\]
and thus $I_1 (\Gamma ) - I_0 (\Gamma ) = \pm d p_* \tilde{I}(\Gamma )$.
\end{proof}

In completely similar way we can prove the following.

\begin{prop}\label{prop_indep_vol2}
Suppose $g$, $n$ and $j$ are as in Proposition \ref{prop_indep_vol}.
Then the cohomology class $[I(\Gamma )]$ for a graph cocycle $\Gamma$ does not depend on the choice of $vol_{S^{j-1}}$.
\end{prop}

\section{The Haefliger invariant}\label{sec_Haefliger}

Here we assume $n-j \ge 3$ is odd, $j \ge 2$ and $2n-3j-3 \ge 0$.
Then at present we cannot prove that $I$ is a cochain map, but in this section we describe a $(2n-3j-3)$-cocycle
$\calH :=I(H)+c$ of $\emb{n}{j}$ using a graph cocycle $H \in \D^{2,0}_0$ and some `correction term' $c$.

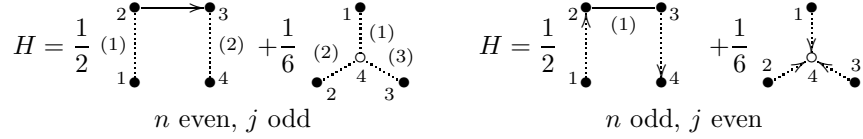
\begin{figure}[htb]
\[
 \begin{xy}
  (-6,5)*{H=\dfrac{1}{2}},(24,5)*{+\dfrac{1}{6}},
  (5,0)*{\bullet}="A", (5,10)*{\bullet}="B", (15,10)*{\bullet}="C", (15,0)*{\bullet}="D",
  {\ar@{.} "A";"B"^<{1}^(.5){(1)}^>{2}},{\ar "B";"C"},{\ar@{.} "C";"D"^<{3}^(.5){(2)}^>{4}},
  (35,3.33)*{\circ}="E",(35,10)*{\bullet}="F",(29.227,0)*{\bullet}="G",(40.773,0)*{\bullet}="H",
  {\ar@{.} "E";"F"_(.5){(1)}^>{1}},{\ar@{.} "E";"G"^<{4}_(.5){(2)}^>{2}},{\ar@{.} "E";"H"^(.65){(3)}_>{3}},
  (18,-5)*{n \text{ even, } j \text{ odd}},
  (56,5)*{H=\dfrac{1}{2}},(84,5)*{+\dfrac{1}{6}},
  (65,0)*{\bullet}="a", (65,10)*{\bullet}="b", (75,10)*{\bullet}="c", (75,0)*{\bullet}="d",
  {\ar@{.>} "a";"b"^<{1}^>{2}},{\ar@{-} "b";"c"_(.5){(1)}},{\ar@{.>} "c";"d"^<{3}^>{4}},
  (95,3.33)*{\circ}="e",(95,10)*{\bullet}="f",(89.227,0)*{\bullet}="g",(100.773,0)*{\bullet}="h",
  {\ar@{.>} "f";"e"_<{1}},{\ar@{.>} "g";"e"^<{2}_>{4}},{\ar@{.>} "h";"e"_<{3}},
  (78,-5)*{n \text{ odd, } j \text{ even}}
 \end{xy}
\]
\caption{Graph cocycle $H$}\label{fig_H}
\end{figure}

The graph cocycle $H$ is shown in Figure \ref{fig_H}.
We call the first graph (with four i-vertices) $H_1$ and the second $H_2$; $H =H_1 /2 +H_2 /6$.
By using the rules in Proposition \ref{def_signs}, we can see that $H$ is indeed a cocycle.

\begin{thm}\label{thm_closed}
Suppose $n-j \ge 3$ is odd, $j \ge 2$ and $2n-3j-3 \ge 0$.
There exists a differential form $c\in \Omega^{2n-3j-3}_{DR}(\emb{n}{j})$ such that $\calH := I(H)+c$ is closed.
\end{thm}

Though rigorous definition of $c$ and the proof of Theorem \ref{thm_closed} can be found in \S \ref{subsec_proof_closed},
we give a rough explanation here.
Let $\I_j (\R^n )$ be the space of linear injections $\R^j \to \R^n$.
Then $c:=-p_* D^* \mu$, where $D:\emb{n}{j} \times \R^j \to \I_j (\R^n )$ is the derivation map,
$p:\emb{n}{j}\times \R^j \to \emb{n}{j}$ is the first projection and the form
$\mu \in \Omega^{2n-2j-3}_{DR}(\I_j (\R^n ))$ is given so that $d\mu$ describes the contribution of a boundary stratum
of $C_{H_2}$ corresponding to the collision of all the four points.
See Definition \ref{def_c} for details.

Below we prove Theorem \ref{thm_main1}, which states that $\calH$ gives a non-zero cohomology class and is the
Haefliger invariant when $2n-3j-3 = 0$.

\subsection{Additivity}\label{subsec_additive}

Here we assume $2n-3j-3=0$, which implies $n=6k$ and $j=4k-1$ for some $k \ge 1$.
Then $\calH =I(H)+c$ is an isotopy invariant for long $(4k-1)$-knots in $\R^{6k}$.
What we will show in this subsection is the `additivity' of the invariant $\calH$.

\begin{prop}\label{prop_additive}
The invariant $\calH$ is additive under the connect-sum;
for any $f_+ ,f_- \in \emb{6k}{4k-1}$, the Kronecker pairing satisfies
\[
 \pair{\calH}{f_+ \sharp f_-} = \pair{\calH}{f_+} + \pair{\calH}{f_-}.
\]
\end{prop}

\noindent
{\bf Notation.}
We will show in Lemma \ref{lem_indep_H_vol} that the invariant $\calH$ is independent of the choice of (anti-)symmetric
volume forms.
So we may choose (anti-)symmetric volume forms $vol_{S^{N-1}}$, $N=4k-1$ or $6k$ such that
\begin{itemize}
\item
 their supports are contained in the sufficiently small neighborhoods of the poles
 $p^{N-1}_{\pm} := (0,\dots ,0,\pm 1) \in S^{N-1}$, and
\item
 they are ``$O(N-1)$-invariant,'' that is, $A^* vol_{S^{N-1}}=(\det A)vol_{S^{N-1}}$ for any $A\in O(N-1)$ regarded as
 in $O(N-1)\oplus 1 \subset O(N)$ (in other words $A\in O(N)$ fixes $x_N$-axis).
\end{itemize}

\smallskip

For any $f \in \emb{n}{j}$, the {\em support} of $f$, denoted by $\text{supp}(f) \subset [-1,1]^j$, is defined by
\[
 \text{supp}(f) := \overline{\{ x \in \R^j \, | \, f(x) \ne (x,0) \in \R^j \times \{ 0 \}^{n-j} \}} .
\]
Since we are considering an isotopy invariant, we may suppose
$\text{supp}(f_{\pm}) \subset B^{4k-1}_{\pm}(\varepsilon )$ and
$f_{\pm}(\text{supp}(f_{\pm})) \subset B^{6k}_{\pm}(\varepsilon /2) \cup (\R^{4k-1}\times \{ 0\}^{2k+1})$, where
\[
 B^m_{\pm} (\varepsilon ) := \{ (x_1 ,\dots ,x_m ) \in \R^m \, | \, ( x_1 \pm 1/2 )^2 +x^2_2 +\dots +x^2_m
 < \varepsilon^2 \} ,
\]
and $\varepsilon >0$ is a sufficiently small number.

First we compute $\pair{\calH}{f_+} = \pair{I(H_1 )}{f_+}/2 +\pair{I(H_2 )}{f_+}/6 + \pair{c}{f_+}$.
The first term $\pair{I(H_1 )}{f_+}$ is equal to
\[
 \int_{C_{H_1}(f_+ )}\omega_{H_1} ,
\]
where $C_{H_1}(f) \approx C^o_4 (\R^{4k-1})$ is the fiber of $\pi_{H_1}$ over $f \in \emb{6k}{4k-1}$.
Recall that the integrand is $\theta_{12}\eta_{23}\theta_{34}$ restricted to the fiber.
Since we choose $vol_{S^{N-1}}$ with support `localized' around $p^{N-1}_{\pm}$, we must know for which configurations
$x=(x_1 ,\dots ,x_4 ) \in C_{H_1}(f_+ )= C_4 (\R^{4k-1})$ all the three vectors $\varphi^{\theta}_{12}(x)$,
$\varphi^{\theta}_{34}(x)$ and $\varphi^{\eta}_{23}(x)$ are near $p^{6k-1}_{\pm}$ or $p^{4k-2}_{\pm}$;
$\omega_{H_1}$ does not vanish only on the subspace of such configurations.

\begin{lem}\label{lem_localize1}
Define $X_+ := C^o_4 (B^{4k-1}_+ (\varepsilon )) \subset C_{H_1}(f_+ )$.
Then
\[
 \pair{I(H_1 )}{f_+} = \int_{X_+}\omega_{H_1} .
\]
\end{lem}

\begin{proof}
For example, consider $(x_1 ,\dots ,x_4 ) \in C_{H_1}(f_+ ) \setminus X_+ $ with $x_1 \not\in B^{4k-1}_+ (\varepsilon )$.
Then by our choice of $f_+$, the vector $\varphi^{\theta}_{12} (x_1 ,\dots ,x_4 )$ cannot be in the support of
$vol_{S^{6k-1}}$ (a neighborhood of the pole $p^{6k-1}$).
So the form $\theta_{12}$ vanishes on such configurations.
In this way we can see that the integrand does not vanish only on $X_+$, since outside of $X_+$ the images of
$\varphi^{\theta}_{12}$ and $\varphi^{\theta}_{34}$ cannot intersect with the supports of $vol_{S^{6k-1}}$
simultaneously.
\end{proof}

Next we compute $\pair{I(H_2 )}{f_+}$ in similar fashion.
This equals
\[
 \int_{C_{H_2}(f_+ )}\omega_{H_2},
\]
where $\omega_{H_2} = \theta_{14}\theta_{24}\theta_{34}$.
Recall $C_{H_2}(f_+ ) \subset C^o_3 (\R^{4k-1}) \times \R^{6k}$.

\begin{lem}\label{lem_localize2}
Define a subspace $Y_+ \subset C_{H_2}(f_+ )$ by
\[
 Y_+ := \{ (x_1 ,x_2 ,x_3 ;y) \in C_{H_2}(f_+ ) \, | \, \text{at least two }x_p \in B^{4k-1}_+ (\varepsilon )\} .
\]
Then
\[
 \pair{I(H_2 )}{f_+} = \int_{Y_+} \omega_{H_2}.
\]
\end{lem}

\begin{proof}
This is because the image of
\[
 \varphi_{H_2}:= \varphi^{\theta}_{14} \times \varphi^{\theta}_{24} \times \varphi^{\theta}_{34}
 : C_{H_2}(f_+ ) \setminus Y_+ \longrightarrow (S^{6k-1})^3
\]
is of positive codimension, and hence the pullback of the top form of $(S^{6k-1})^3$ must vanish on
$C_{H_2}(f_+ ) \setminus Y_+$.
Indeed, outside $Y_+$ at least two points $f_+ (x_p )$, $f_+ (x_q )$ are on the standardly embedded $\R^{4k-1}$, hence
the image of $C_{H_2}\setminus Y_+$ is the set of $(u_1 ,u_2 ,u_3 ) \in (S^{6k-1})^3$ such that one vector
$u_p$ has to be in the linear subspace $(\R^{4k-1}\times \{ 0 \}^{2k+1})+\R u_q$ of $\R^{6k}$, hence is of
codimension $\ge 2k$.
\end{proof}

Now consider the third term $\pair{c}{f_+}$.
This is an integration of $D(f_+ )^* \mu$ over $\R^{4k-1}$, where $D(f_+ ):\R^{4k-1} \to \I_{4k-1}(\R^{6k})$ is given by
$x\mapsto (df_+ )_x$ (for the definition of $\mu$, see Definition \ref{def_c} and the explanation after
Theorem \ref{thm_closed}).
Since $D(f_+ )$ is constant outside $\text{supp}(f_+ )$, we obtain the following.

\begin{lem}\label{lem_localize3}
$\displaystyle \pair{c}{f_+} = \int_{B^{4k-1}_+ (\varepsilon)}D(f_+ )^* \mu$.
\end{lem}

Finally we can compute $\pair{\calH}{f_+ \sharp f_-}$.
The long knot $f_+ \sharp f_-$ is isotopic to
\[
 f_+ \sharp f_- (x) :=
 \begin{cases}
  f_{\pm}(x)  & \text{if } x \in B^{4k-1}_{\pm}(\varepsilon ) , \\
  (x,0) & \text{if } x \not\in B^{4k-1}_+ (\varepsilon ) \sqcup B^{4k-1}_- (\varepsilon ).
 \end{cases}
\]
First, by completely similar arguments to Lemma \ref{lem_localize1}, only the configurations $x=(x_1 ,\dots ,x_4 )$
such that $x_1$ and $x_2$ are in the same $B^{4k-1}_+ (\varepsilon )$ or $B^{4k-1}_- (\varepsilon )$, and so is
$(x_3 ,x_4 )$, can contribute to $\pair{I(H_1)}{f_+ \sharp f_-}$.
Moreover, if $x_2$ and $x_3$ are in the different $B^{4k-1}(\varepsilon )$'s, then the image of the direction map
$\varphi^{\eta}_{23}(x)$ cannot be in $\text{supp}(vol_{S^{4k-2}})$.
Thus only the configurations with all four $x_p$'s being in the same $B^{4k-1}_+ (\varepsilon )$ or
$B^{4k-1}_- (\varepsilon )$ can contribute to the pairing, and hence
\begin{equation}\label{eq_sum_1}
 \pair{I(H_1 )}{f_+ \sharp f_-} = \int_{X_+ \sqcup X_-}\omega_{H_1} =\int_{X_+}\omega_{H_1} + \int_{X_-}\omega_{H_1},
\end{equation}
where $X_- :=C^o_4 (B^{4k-1}_- (\varepsilon ))$.

Next, similarly to Lemma \ref{lem_localize3},
\begin{equation}\label{eq_sum_2}
 \pair{c}{f_+ \sharp f_-} = \int_{B^{4k-1}_+ (\varepsilon ) \sqcup B^{4k-1}_- (\varepsilon )}D(f_+ \sharp f_- )^* \mu
 =\sum_{i=\pm}\int_{B^{4k-1}_i (\varepsilon )}D(f_i )^* \mu .
\end{equation}

\begin{figure}[htb]
\[
 \begin{xy}
 (0,15)*{\R^n},(10,3)*{\R^j},
 {\ar@{-} (0,0);(50,0)},{\ar@{-}(0,0);(10,7)},{\ar@{-}(10,7);(60,7)},{\ar@{-}(50,0);(60,7)},
 (29,6)*{\bullet}="A",(30,2)*{\bullet}="B",(33,4)*{\bullet}="C", (31,20)*{\circ}="D",
 {\ar@{.>}"A";"D"^(.15){x_1}^>{y}},{\ar@{.>}"B";"D"^<{x_2}},{\ar@{.>}"C";"D"_<{x_3}}
 \end{xy}
\]
\caption{Configurations which contribute to $\pair{I(H_2 )}{f_+ \sharp f_-}$}
\label{fig_contribution}
\end{figure}
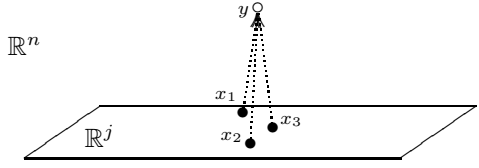

As for $I(H_2 )$, similar argument to Lemma \ref{lem_localize2} shows
\begin{equation}\label{eq_sum_3}
 \pair{I(H_2 )}{f_+ \sharp f_-} = \int_{Y_+ (f_+ \sharp f_- )}\omega_{H_2} + \int_{Y_- (f_+ \sharp f_- )}\omega_{H_2}
 + \int_{Y_0 (f_+ \sharp f_- )}\omega_{H_2},
\end{equation}
where $Y_- (f_+ \sharp f_- ) \subset C_{H_2}(f_+ \sharp f_- )$ is defined in similar way as $Y_+$ with $B_+$'s replaced
by $B_-$'s, and $Y_0 (f_+ \sharp f_- ) \subset C_{H_2}(f_+ \sharp f_- )$ consists of $(x_1 ,x_2 ,x_3 ;y)$ with exactly
one $x_p$ in each $B^{4k-1}_{\pm}(\varepsilon )$.
By our choice of volume forms, the integrand $\omega_{H_2}$ does not vanish only on the subspace consisting of
$(x_1 ,x_2 ,x_3 ,y)$ with
\begin{enumerate}
\item
 $y$ is very near `$(0,\dots ,0,\pm \infty )$,' or
\item
 the first $6k-1$ factors of $f_+ (x_1 )$, $f_+ (x_2 )$, $f_+ (x_3 )$ and $y$ are close to each other (see Figure
 \ref{fig_contribution}).
\end{enumerate}
On $Y_0 (f_+ \sharp f_- )$ the condition (2) cannot hold, so we can restrict the integral over $Y_0 (f_+ \sharp f_- )$
to that over the subspace of $Y_0 (f_+ \sharp f_- )$ with $\abs{y}$ large.
Then since the image of $f_+ \sharp f_-$ is very close to that of the trivial knot $f_0$
(if $\varepsilon$ is very small relative to $\abs{y}$),
\[
 \int_{Y_0 (f_+ \sharp f_- )}\omega_{H_2} - \int_{C_{H_2}(f_0 )}\omega_{H_2} = O(\varepsilon ).
\]
But since $\text{supp}(f_0 )=\emptyset$, Lemma \ref{lem_localize2} shows that $\int_{C_{H_2}(f_0 )}\omega_{H_2}=0$.
Hence
\begin{equation}\label{eq_sum_4}
 \int_{Y_0 (f_+ \sharp f_- )}\omega_{H_2} = O(\varepsilon ).
\end{equation}
$Y_{\pm} (f_+ \sharp f_- )$ is decomposed into two subspaces; in one space, two $x_p$, $x_q$ are in
$B^{4k-1}_{\pm}(\varepsilon )$ and remaining $x_r$ is in $B^{4k-1}_{\mp}(\varepsilon )$, and in the other space $x_r$
is not in $B^{4k-1}_{\mp} (\varepsilon )$.
In the former space, the condition (2) above cannot hold, so $\abs{y}$ has to be large.
So replacing $f_+ \sharp f_-$ with $f_{\pm}$ changes the integral by $O(\varepsilon )$;
\begin{equation}\label{eq_sum_5}
 \int_{Y_{\pm}(f_+ \sharp f_- )}\omega_{H_2} = \int_{Y_{\pm}(f_{\pm})}\omega_{H_2} + O(\varepsilon ).
\end{equation}
Summing up \eqref{eq_sum_1}, \eqref{eq_sum_2} and \eqref{eq_sum_3} and substituting \eqref{eq_sum_4} and
\eqref{eq_sum_5}, we have
\begin{align*}
 \pair{I(H)}{f_+ \sharp f_-}
  &= \sum_{i=\pm} \left(\frac{1}{2}\int_{X_i}\omega_{H_1} + \frac{1}{6}\int_{Y_i (f_i )}\omega_{H_2}
   + \int_{B^{4k-1}_i (\varepsilon)}D(f_i )^* \mu  \right) + O(\varepsilon ) \\
 &= \pair{I(H)}{f_+} + \pair{I(H)}{f_-} + O(\varepsilon ).
\end{align*}
This completes the proof of Proposition \ref{prop_additive}, since $\pair{\calH}{f_+ \sharp f_-}$ is independent of
$\varepsilon$.

\begin{rem}
By similar argument, the order two invariant \cite{Rossi_thesis, Watanabe07} for $\emb{m+2}{m}$ ($m \ge 3$ is odd) is
proved to be additive under the connect-sum.
\end{rem}

\subsection{Evaluation}\label{subsec_evaluation}
Here we prove the last part of Theorem \ref{thm_main1}.

We know by \cite{Haefliger66} that $\pi_0 (\emb{6k}{4k-1})$ forms a group under the connect-sum and is isomorphic to
$\Z$.
Our goal here is the following.

\begin{thm}\label{thm_Haefliger}
When $n=6k$ and $j=4k-1$ for some $k \ge 1$, then $\pair{\calH}{S}$ is equal to $\pm 1$ for the generator $S$ of
$\pi_0 (\emb{6k}{4k-1})$.
\end{thm}

Since both the Haefliger invariant and our $\calH$ are additive under the connect-sum (see \cite{Haefliger62} and
Proposition \ref{prop_additive}), Theorem \ref{thm_Haefliger} says that $\calH$ is the Haefliger invariant.

We will use the generator of $\pi_0 (\emb{6k}{4k-1})$ given by Budney \cite{Budney08} and Roseman-Takase
\cite{RosemanTakase07}.

\subsubsection{Deform-spinning}\label{subsubsec_deform_spinning}
Given an $N$-fold based loop $\gamma \in \Omega^N \emb{n}{j}$ (for any $n,j$) represented by a smooth map
$\gamma : \R^N \to \emb{n}{j}$ with $\gamma (x)\equiv \iota$ (the trivial long $j$-knot) for any $x \not\in [-1,1]^j$,
we have $S_{\gamma}\in \emb{N+n}{N+j}$ defined by
\[
 S_{\gamma}(x,t) :=(x,\gamma (x)(t)) \in \R^N \times \R^n , \quad x\in \R^N ,\ t\in \R^j.
\]
This is a special case of `deform-spinning' construction \cite{Roseman89}.
Putting $N=1$, we obtain the {\em graphing map} $\gr : \Omega \emb{n}{j} \to \emb{n+1}{j+1}$ given in \cite{Budney08}.
It has been shown in \cite{Budney08} that all the following maps
\[
 \pi_{4k-2}(\emb{2k+2}{1}) \xrightarrow{\gr} \pi_{4k-3}(\emb{2k+3}{2}) \xrightarrow{\gr} \dots \xrightarrow{\gr}
 \pi_0 (\emb{6k}{4k-1}) \cong \Z
\]
are isomorphisms, so $\pi_0 (\emb{6k}{4k-1})\cong \Z$ is generated by $S_{\psi}:=\gr^{4k-2}(\psi )$;
\begin{equation}\label{eq_deform_spinning}
 S_{\psi}(z) := (x,\psi(x)(t)) \in \R^{4k-2}\times\R^{2k+2}, \quad z=(x,t) \in \R^{4k-2} \times \R^1 ,
\end{equation}
where $\psi : \R^{4k-2} \to \emb{2k+2}{1}$ is a map which gives the generator of $\pi_{4k-2}(\emb{2k+2}{1})$.
In fact $\pi_{4k-2}(\emb{2k+2}{1})$ is the first non-vanishing homotopy group of $\emb{2k+2}{1}$ (for example see
\cite{Budney08, Tourtchine04_2}), and the image of its generator $\psi$ via Hurewicz isomorphism is given as follows
(probably the dual cocycle to $\psi$ first appeared in \cite{Vassiliev01}; see also
\cite{Budney08, CCL02, Longoni04, K07}).

Consider a `long immersion' $f:\R^1 \to \R^{2k+2}$ which has only two transversal doublepoints
$z_i = f(\xi_i )=f(\xi_{i+2} )$, $i=1,2$, $\xi_1 < \xi_2 < \xi_3 < \xi_4$ (see Figure \ref{fig_f}).
Consider the unit sphere $S^{2k-1}_i$ in the complementary $2k$-space to $\R f' (\xi_i ) \oplus \R f' (\xi_{i+2})$ in
$T_{z_i}\R^{2k+2}$.
At each doublepoint $z_i$, we can `lift' the segments $f(\xi_{i+2} -\varepsilon , \xi_{i+2} +\varepsilon )$ to the
directions of complementary $S^{2k-1}_i$ to obtain non-singular embeddings (see Figure \ref{fig_blowup}).
More precisely, for $(w_1 ,w_2 ) \in (S^{2k-1})^2$, define $\psi (w_1 ,w_2 )\in \emb{2k+2}{1}$ by
\[
 \psi (w_1 ,w_2 )(t) :=
 \begin{cases}
  f(t)+\left( \delta \exp \left[ \frac{1}{(t-\xi_{i+2} )^2 -\varepsilon^2}\right] \right) w_i
   & \abs{t-\xi_{i+2}} < \varepsilon ,\ i=1,2 \\
  f(t) & \text{otherwise}
 \end{cases}
\]
where $\delta$ and $\varepsilon$ are positive small numbers.
Thus a map $\psi =\psi_{\delta}: (S^{2k-1})^2 \to \emb{2k+2}{1}$ is defined, and it is known that $\psi$ induces an
isomorphism $\psi_* : H_{4k-2}((S^{2k-1})^2) \xrightarrow{\cong} H_{4k-2}(\emb{2k+2}{1})$.

\begin{figure}[htb]
\[
 \begin{xy}
 {\ar@{.>}(0,0);(85,0)^>{1}},{\ar@{.>}(40,-10);(40,12)^>{3,\dots ,2k+2}},{\ar@{.>}(20,-10);(60,10)^>{2}},
 {\ar@{-}(5,0);(50,0)},{\ar@{-}(30,5);(39,5)},{\ar@{-}(41,5);(60,5)},
 {\ar@{-}@(r,r)(50,0);(60,5)},{\ar@{-}@(l,l)(30,5);(15,-5)},
 {\ar@{-}(15,-5);(40,-5)},
 {\ar@{-}@(r,l)(40,-5);(63,3)},
 {\ar@(r,l)(67,3);(80,0)},
 (12,2)*{z_1}, (55,-2)*{z_2}, (68,5)*{c}
 \end{xy}
\]
\caption{Immersion $f$ in $\R^{2k+2}$}\label{fig_f}
\end{figure}
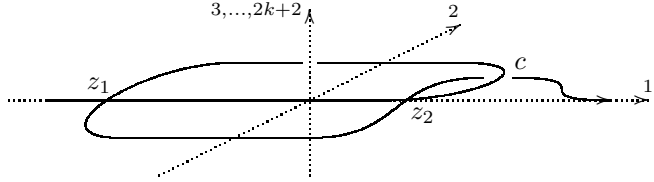
\begin{figure}[htb]
\[
 \begin{xy}
 {\ar (0,0);(40,0)}, {\ar (4,-8);(36,8)}, {\ar@{.>}(20,-6);(20,12)},
 (0,3)*{f}, (22,-2)*{z_i}, (20,9)*{\bullet}, (17,9)*{w_i},
 (50,0)*{\leadsto},
 {(60,0);(100,0) **\crv{(65,0)&(70,3)&(80,12)&(90,3)&(95,0)}?>*\dir{>}},
 {(60,0);(100,0) **\crv{~*{.}(65,0)&(70,-3)&(80,-12)&(90,-3)&(95,0)}?>*\dir{>}},
 (60,3)*{\psi (w_i )},
 {\ar@{-}(64,-8);(88,4)}, {\ar (90,5);(96,8)}, {\ar@{.>}(80,-6);(80,12)},
 (70,-12)*{I_i}="I", {\ar@{.}"I";(78,-1)},
 (95,-10)*{W_{i+2}}="W", {\ar@{.}"W";(85,7.2)},{\ar@{.}"W";(85,-7.2)}
 \end{xy}
\]
\caption{Definition of $\psi$}\label{fig_blowup}
\end{figure}
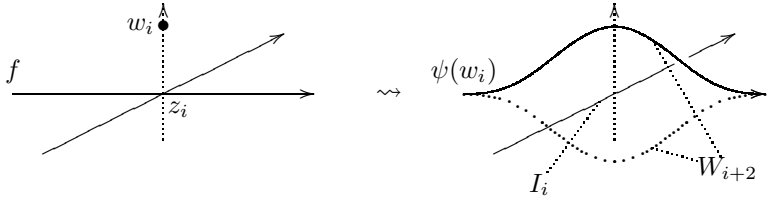

The following is proved by inspection (see Figure \ref{fig_blowup}).

\begin{lem}\label{lem_linking_number}
The union
\[
 W_3 := \bigcup_{w_1 \in S^{2k-1}}\psi (w_1 ,w_2 ) (\xi_3 -\varepsilon ,\xi_3 + \varepsilon ) \subset \R^{2k+2}
\]
is homeomorphic to $S^{2k}$.
Consider the map
\[
 \varphi : W_3 \times I_1 \longrightarrow S^{2k+1}, \quad (x,y) \longmapsto \frac{y-x}{\abs{y-x}},
\]
where $I_1 := f(\xi_1 -\varepsilon ,\xi_1 +\varepsilon )$.
The linking number of $W_3$ with $I_1$ is one; more precisely, the limit $\delta \to 0$ of the integration of
$\varphi^* vol_{S^{2k+1}}$ over $W_3 \times I_1$ is $\pm 1$.
Similar statement holds for
\[
 W_4 := \bigcup_{w_2 \in S^{2k-1}}\psi (w_1 ,w_2 ) (\xi_4 -\varepsilon ,\xi_4 + \varepsilon ) \quad \text{and}\quad
 I_2 := f(\xi_2 -\varepsilon ,\xi_2 +\varepsilon ).
\]
\end{lem}

Fixing the natural projection $\pi : [-1/2,1/2]^{2k-1} \to S^{2k-1}$ such that
\begin{itemize}
\item
 $q^{2k-1}_{\pm}:=\pi^{-1}(p^{2k-1}_{\pm}) \in (-1/2,1/2)^{2k-1}$ (recall
 $p^{2k-1}_{\pm}:=(0,\dots ,0,\pm 1) \in S^{2k-1}$),
\item
 $\pi (\partial [-1/2,1/2]^{2k-1})$ consists of a single point $u_0 := (1,0,\dots ,0)$,
\end{itemize}
we can think of $\psi$ as $\psi :[-1/2,1/2]^{4k-2} \to \emb{2k+2}{1}$.
We can extend $\psi$ to whole $\R^{4k-2}$ so that $\psi (x)\equiv \iota$ if $x \not\in [-1,1]^{4k-2}$
(see also \S \ref{sssec_choice_psi}), since $\emb{2k+2}{1}$ is $(4k-3)$-connected as mentioned above.
Thus we have $S_{\psi}\in \emb{6k}{4k-1}$ defined as (\ref{eq_deform_spinning}).

\begin{thm}[\cite{Budney08, RosemanTakase07}]
If $\psi$ is as above, then $S_{\psi}$ generates $\pi_0 (\emb{6k}{4k-1})$.
\end{thm}

\begin{rem}
In \cite{RosemanTakase07} a generator of $\pi_0 (\text{Emb}\, (S^{4k-1},S^{6k}))$ was defined by the
deform-spinning $\psi :(S^{2k-1})^2 \to \emb{2k+2}{1}$ along the torus.
But it is not difficult to see that such a spinning also gives an element of $\emb{6k}{4k-1}$ which is isotopic to our
$S_{\psi}$ given by the graphing map \cite{Budney08}.
\end{rem}

The above construction can be done for the generator of $\pi_{2d-2}(\emb{d+2}{1})$ for any $d>1$ to obtain
$S_{\psi}\in \emb{3d}{2d-1}$.
Since we assumed that $n-j=d+1$ is odd, we put $d=2k$ here (otherwise there exist orientation reversing automorphisms
of graphs $H_i$, and hence they vanish).
When $d$ is odd, any generator of $\pi_0 (\emb{3d}{2d-1})\cong \Z /2$ would not be detected by any de Rham cohomology
classes.

\subsubsection{A suitable choice of $\psi : \R^{4k-2}\to \emb{2k+2}{1}$}\label{sssec_choice_psi}
Below we will compute
\begin{equation}\label{eq_evaluation}
 \pair{\calH}{S_{\psi}} = \frac{1}{2}\int_{C_{H_1}(S_{\psi})}\omega_{H_1}
 + \frac{1}{6}\int_{C_{H_2}(S_{\psi})}\omega_{H_2} + \int_{\R^{4k-1}} D(S_{\psi})^* \mu ,
\end{equation}
where $D(S_{\psi}):\R^{4k-1} \to \I_{4k-1}(\R^{6k})$ is defined by $D(S_{\psi})(x) := (dS_{\psi})_x$, and $\mu$ is
defined in the sentence after the proof of Lemma~\ref{lem_hatomega_closed}.
To simplify the computation of \eqref{eq_evaluation}, we need to choose a favorable extension
$\psi : \R^{4k-2} \to \emb{2k+2}{1}$.

First, we take the immersion $f$ as in Figure \ref{fig_f};
almost all the image of $\psi (x)$ is contained in the 2-plane $\R^2 \times \{ 0\}^{2k}\subset \R^{2k+2}$,
except for the neighborhoods of $z_1$, $z_2$ and another crossing $c$ which corresponds to $f(\zeta_1 )$ and
$f(\zeta_2 )$, $\xi_2 <\zeta_1 <\xi_3$, $\xi_4 <\zeta_2 <1$.
We suppose that the arc $\psi (x)(\zeta_1 -\varepsilon ,\zeta_1 +\varepsilon)$ is in the 2-plane,
$\psi (x)(\zeta_2 -\varepsilon ,\zeta_2 +\varepsilon)$ is in the 3-plane $\R^3 \times \{ 0\}^{2k-1}\subset \R^{2k+2}$,
and that the resolutions of segments $\psi (x)(\xi_{i+2} -\varepsilon ,\, \xi_{i+2} +\varepsilon )$ occur in
$\{ 0\}^2 \times \R^{2k}$.

Next we suppose that the image of $\psi (x)$ is `almost' in $\R^2$;
\begin{equation}\label{eq_thin}
 \psi (x)(\R^1 ) \subset \R^2 \times (-\delta ,\delta )^{2k}\quad \text{for any }x\in \R^{4k-2}.
\end{equation}
Lastly we suppose $\psi$ is `symmetric' in the following sense; define
\begin{align*}
 X&:= \{ (s,t,u) \in \R^{2k-1}\times \R^{2k-1}\times \R^1 \, | \,
  s \not\in [-1/2 ,1/2]^{2k-1}, t \in [-1/2 ,1/2]^{2k-1}\} ,\\
 Y&:= \{ (s,t,u) \in \R^{2k-1}\times \R^{2k-1}\times \R^1 \, | \, 
  s \in [-1/2 ,1/2]^{2k-1}, t \not\in [-1/2 ,1/2]^{2k-1}\} ,\\
 Z&:= \{ (s,t,u) \in \R^{2k-1}\times \R^{2k-1}\times \R^1 \, | \,  s,t \not\in [-1/2 ,1/2]^{2k-1} \}
\end{align*}
and $X',Y',Z'$ as the images of $X,Y,Z$ respectively via the first $4k-2$ projection $pr_{4k-2}:\R^{4k-1}\to \R^{4k-2}$
(see Figure \ref{fig_outside_zero}).
Define $\psi : X' \to \emb{2k+2}{1}$ by
\[
 \psi (s,t)=h_{\abs{s}/\abs{s'}}(\pi (t))
\]
where $h_{\alpha} : S^{2k-1}\to \emb{2k+2}{1}$ ($\alpha \ge 1$) is a one parameter family with $h_1 =\psi (u_0 ,-)$ and
$h_{\alpha} \equiv \iota$ (the trivial long knot) for $\alpha \ge 2$, and $s' \in \partial [-1/2,1/2]^{2k-1}$ is the
unique element with $s=\alpha s'$ for some $\alpha \ge 1$.
Using similar one parameter family connecting $\psi (-,u_0 )$ (resp.\ $\psi (u_0 ,u_0 )$) and $\iota$, we can define
$\psi$ on $Y'$ (resp.\ on $Z'$) (see Example \ref{ex_outside_zero} and Figure \ref{fig_outside_zero}).
We can modify $\psi$ so that it is smooth.
The following properties of $\psi$ will be important below;
\begin{align}
 \psi (s,t) &= \psi (-s,t) \quad \text{if }(s,t)\in X'\cup Z', \label{eq_psi_symmetry_1} \\
 \psi (s,t) &= \psi (s,-t) \quad \text{if }(s,t)\in Y'\cup Z'. \label{eq_psi_symmetry_2}
\end{align}

\begin{ex}\label{ex_outside_zero}
Consider the case $k=1$ (see Figure \ref{fig_outside_zero}).
\begin{figure}[htb]
\[
 \begin{xy}
 {\ar@{.>}(-25,0);(25,0)^(.38){-1/2}^>{s}_(.53){O}_(.65){1/2}},{\ar@{.>}(0,-20);(0,20)^(.3){-1/2}_(.7){1/2}_>{t}},
 {\ar@{-}(-20,10);(-10,10)},{\ar@{-}(-10,10);(-10,15)},
 {\ar@{.}(-10,12);(-12,10)},{\ar@{.}(-10,14);(-14,10)},{\ar@{.}(-12,14);(-16,10)},{\ar@{.}(-14,14);(-18,10)},
  {\ar@{.}(-16,14);(-18,12)},
 {\ar@{-}(20,10);(10,10)},{\ar@{-}(10,10);(10,15)},
 {\ar@{.}(10,12);(12,10)},{\ar@{.}(10,14);(14,10)},{\ar@{.}(12,14);(16,10)},{\ar@{.}(14,14);(18,10)},
  {\ar@{.}(16,14);(18,12)},
 {\ar@{-}(-20,-10);(-10,-10)},{\ar@{-}(-10,-10);(-10,-15)},
 {\ar@{.}(-10,-12);(-12,-10)},{\ar@{.}(-10,-14);(-14,-10)},{\ar@{.}(-12,-14);(-16,-10)},{\ar@{.}(-14,-14);(-18,-10)},
  {\ar@{.}(-16,-14);(-18,-12)},
 {\ar@{-}(20,-10);(10,-10)},{\ar@{-}(10,-10);(10,-15)},
 {\ar@{.}(10,-12);(12,-10)},{\ar@{.}(10,-14);(14,-10)},{\ar@{.}(12,-14);(16,-10)},{\ar@{.}(14,-14);(18,-10)},
  {\ar@{.}(16,-14);(18,-12)},
 {\ar@{.}(-10,10);(-10,-10)}, {\ar@{.}(-10,-10);(10,-10)}, {\ar@{.}(10,-10);(10,10)}, {\ar@{.}(10,10);(-10,10)},
 {\ar@{.}(-12,10);(-12,-10)}, {\ar@{.}(-14,10);(-14,-10)}, {\ar@{.}(-16,10);(-16,-10)}, {\ar@{.}(-18,10);(-18,-10)},
 {\ar@{.}(12,10);(12,-10)}, {\ar@{.}(14,10);(14,-10)}, {\ar@{.}(16,10);(16,-10)}, {\ar@{.}(18,10);(18,-10)},
 {\ar@{.}(10,12);(-10,12)}, {\ar@{.}(10,14);(-10,14)},
 {\ar@{.}(10,-12);(-10,-12)}, {\ar@{.}(10,-14);(-10,-14)},
 (23,5)*{X'},(-23,5)*{X'}, (-5,17)*{Y'},(-5,-17)*{Y'},
 (21,14)*{Z'},(-21,14)*{Z'},(21,-14)*{Z'},(-21,-14)*{Z'}
 \end{xy}
\]
\caption{Extension $\psi : \R^2 \to \emb{4}{1}$}
\label{fig_outside_zero}
\end{figure}
$X'\to \emb{4}{1}$ is given as a null homotopy $\psi(s,-) :S^1 \to \emb{4}{1}$ ($\abs{s}\ge 1/2$) from
$\psi (\{ 1/2\} \times S^1 )$ to $\iota$.
Since $\psi (\{ 1/2 \} \times [-1/2,1/2])$ and $\psi (\{ -1/2 \} \times [-1/2,1/2])$ are the same cycle of $\emb{4}{1}$,
two homotopies corresponding to $X' \cap \{ s\ge 1/2\}$ and  $X' \cap \{ s\le -1/2\}$ can be taken so that
$\psi (s,t)=\psi (-s,t)$.
This is \eqref{eq_psi_symmetry_1}.
Similarly $Y'\to \emb{4}{1}$ is a homotopy from $\psi (S^1 \times \{ 1/2\} )$ to $\iota$, and $\psi (s,t)=\psi (s,-t)$
on $Y'$.

Eight thick segments $\{ \pm 1/2\} \times \{ \abs{t}\ge 1/2\}$, $\{ \abs{s}\ge 1/2\} \times \{ \pm 1/2\}$ in
Figure \ref{fig_outside_zero} represent the same homotopy from $\psi (1/2, 1/2 )$ to $\iota$.
Thus we can define $\psi$ by
\[
 \psi (r\cos \beta +1/2 ,r\sin \beta +1/2):= \psi (r+1/2 ,1/2) \quad (0 \le \beta \le \pi /4)
\]
on $Z' \cap \{ s,t \ge 1/2\}$, and similarly on other components.
Then it is easy to see $\psi (s,t)=\psi (-s,t)=\psi (s,-t)=\psi (-s,-t)$ on $Z'$.
\end{ex}

Using such an extension $\psi$, we will compute \eqref{eq_evaluation}.

\subsubsection{The first term $\pair{I(H_1 )}{S_{\psi}}$}\label{sssec_first_term}
First let us study on which configurations $x=(x_1 ,\dots ,x_4 ) \in C_{H_1}(S_{\psi})$ the integrand $\omega_{H_1}$
does not vanish, as was done in \S \ref{subsec_additive}.

Let us write $x_p := (s^{(p)},t^{(p)},u^{(p)}) \in (\R^{2k-1})^2 \times \R^1$ ($p=1,\dots ,4$).
Then by \eqref{eq_deform_spinning},
\begin{equation}\label{eq_direction}
\begin{split}
 &S_{\psi}(x_2 ) - S_{\psi}(x_1 )= \\
 &\quad (s^{(2)}-s^{(1)},t^{(2)}-t^{(1)},\psi (s^{(2)},t^{(2)})(u^{(2)})-\psi (s^{(1)},t^{(1)})(u^{(1)})) \\
 &\hskip200pt \in ( \R^{2k-1})^2 \times \R^{2k+2}.
\end{split}
\end{equation}
By \eqref{eq_thin}, the length of the last $2k$ factors of (\ref{eq_direction}) is at most $2\sqrt{2k}\, \delta$.
So $\varphi^{\theta}_{12}(x_1 ,x_2 ) = (S_{\psi}(x_2 ) - S_{\psi}(x_1 ))/\abs{S_{\psi}(x_2 ) - S_{\psi}(x_1 )}$
is near $p^{6k-1}_{\pm}$ only if at least the first $4k-2$ factors of (\ref{eq_direction}) are nearly zero,
that is, $(s^{(1)},t^{(1)})$ is near $(s^{(2)},t^{(2)})$.
Similarly $(s^{(3)},t^{(3)})$ must be near $(s^{(4)},t^{(4)})$.
Moreover $(s^{(2)},t^{(2)})$ must be near $(s^{(3)},t^{(3)})$, since $w=(x_3 -x_2)/\abs{x_3 -x_2}$ is near
$p^{4k-2}_{\pm}$.
Thus we need to consider only $x=(x_1 ,\dots ,x_4 )$ with all $(s^{(p)},t^{(p)})$ ($p=1,\dots ,4$) close to each other.
Notice that they must become closer and closer if we choose $vol_{S^{N-1}}$ with smaller support.

Let $L:=(-1/2-a ,1/2+a)^{4k-2}\times \R^1$ ($a>0$ is a fixed small number) be a small neighborhood of
$[-1/2,1/2]^{4k-2}\times \R^1$, and consider $x \in C_4 (\R^{4k-1}\setminus L)$.

\begin{lem}\label{lem_outside_zero}
If we choose $\psi$ as in \S \ref{sssec_choice_psi} and $vol_{S^{N-1}}$ with sufficiently small support, then
the integration of $\omega_{H_1}$ over $C_4 (\R^{4k-1}\setminus L)$ vanishes.
\end{lem}

\begin{proof}
If $\text{supp}(vol_{S^{N-1}})$ is sufficiently small relative to $\delta >0$, then
$(x_1 ,\dots ,x_4 ) \in C_4 (\R^{4k-1}\setminus L)$ must be such that all $(s^{(p)},t^{(p)})$ ($p=1,\dots ,4$) are
close to each other so that at most one of
$\{ x_1 ,\dots ,x_4 \} \cap (X\setminus L)$ and $\{ x_1 ,\dots ,x_4 \} \cap (Y\setminus L)$ is not empty.

Consider the integration of $\omega_{H_1}$ over $\{ x \in C_4 (\R^{4k-1}\setminus L)\, |\, x \cap X \ne \emptyset \}$.
Notice that in this case $x\in C_4 (X \cup Z)$.
Define an involution $F_1$ of this subspace by
\[
 F_1 (x_1 ,\dots ,x_4 ):=(i_1 x_1 ,\dots ,i_1 x_4 ),
\]
where $i_1 :\R^{4k-1} \to \R^{4k-1}$ is given by $i_1 (s,t,u):=(-s,t,u)$ ($s,t \in \R^{2k-1}$, $u\in \R^1$).
The map $F_1$ preserves the orientation, but $F^*_1 \omega_{H_1}=-\omega_{H_1}$ because
\begin{itemize}
\item
 $\varphi^{\theta}_{12}\circ F_1 =i_1 \circ \varphi^{\theta}_{12}$ by \eqref{eq_psi_symmetry_1} and \eqref{eq_direction}
 (here we abbreviate $i_1 \times \id_{\R^{2k+1}} :\R^{6k}\to \R^{6k}$ to $i_1$), and hence
 $F^*_1 \theta_{12}=(-1)^{2k-1}\theta_{12}=-\theta_{12}$ (since we have assumed
 $i^*_1 vol_{S^{6k-1}}=(-1)^{2k-1}vol_{S^{6k-1}}$; see the remark after Proposition \ref{prop_additive}),
\item
 similarly $F^*_1 \theta_{34}=-\theta_{34}$, and
\item
 $\varphi^{\eta}_{23}\circ F_1 =i_1 \circ \varphi^{\eta}_{23}$ and hence $F^*_1 \eta_{23}=-\eta_{23}$.
\end{itemize}
Hence the integration of $\omega_{H_1}$ over $\{ x \in C_4 (\R^{4k-1}\setminus L)\, |\, x \cap X \ne \emptyset \}$ must
vanish.
Similarly, we can show the vanishing of the integration of $\omega_{H_1}$ over $\{ x\cap Y\ne \emptyset \}$ by considering
an involution $F_2$ given by $F_2 (x_1 ,\dots ,x_4 ):=(i_2 x_1 ,\dots ,i_2 x_4)$, where $i_2 :\R^{4k-1} \to \R^{4k-1}$
is given by $i_2 (s,t,u):=(s,-t,u)$.
The vanishing of the integration of $\omega_{H_1}$ over $\{ x \cap X = x\cap Y = \emptyset \}$ is proved in similar way;
in this case $x \in C_4 (Z\setminus L)$, and hence either $F_1$ or $F_2$ can work.
\end{proof}

By Lemma \ref{lem_outside_zero}, only $x=(x_1 ,\dots ,x_4 )$ with at least one $x_p$ in $L$ (and other $x_q$'s near $L$)
contributes to $\pair{I(H_1 )}{S_{\psi}}$.

Consider the tubular neighborhood $N=\R^{4k-2} \times [-1/2,\, 1/2]$ of $\R^{4k-2}$ in $\R^{4k-1}$.
Define four `fat planes' by
\[
 N_i := \R^{4k-2} \times (\xi_i -\varepsilon ,\, \xi_i + \varepsilon ) \subset N \subset \R^{4k-1}, \quad 1 \le i \le 4
\]
($\xi_i$ and $\varepsilon$ have appeared in the definition of $\psi$).
Recall that we write $q^{2k-1}_{\pm} :=\pi^{-1}(p^{2k-1}_{\pm})$ (we choose $\pi$ so that
$q^{2k-1}_{\pm} \in (-1/2,1/2)^{2k-1}$).

\begin{lem}\label{lem_local_N1}
Let $x_1 , x_2 \in \R^{4k-1}$ be two distinct points with at least one $x_p \in L$.
Then the direction $\varphi^{\theta}_{12}(x_1 ,x_2 )$ is near $p^{6k-1}_{\pm}$ only if
\begin{itemize}
\item $(s^{(1)},t^{(1)})$ is near $(s^{(2)},t^{(2)})$,
\item $x_1 \in N_i$ and $x_2 \in N_{i+2}$ for some $i =1,\dots ,4$ (here the suffixes are understood modulo $4$), and
\item $s^{(2)}$ (resp.\ $t^{(2)}$, $s^{(1)}$, $t^{(1)}$) is near $q_{\pm}$ when $(x_1 ,x_2 )\in N_1 \times N_3$
 (resp.\ $N_2 \times N_4$, $N_3 \times N_1$, $N_4 \times N_2$).
\end{itemize}
\end{lem}

\begin{proof}
We have already seen that $(s^{(1)},t^{(1)})$ must be near $(s^{(2)},t^{(2)})$.
But it is not enough; the direction determined by the last $2k+2$ factors
$\psi (s^{(2)},t^{(2)})(u^{(2)})-\psi (s^{(1)},t^{(1)})(u^{(1)})$ of \eqref{eq_direction} must be near $p^{2k+1}_{\pm}$,
and it is the case only if the two points $\psi (s^{(p)},t^{(p)})(u^{(p)})$ ($p=1,2$) are around a self-intersection
$z_i$ of $f$ (see Figure \ref{fig_blowup}) which is resolved in a direction near $p^{2k+1}_{\pm}$.
Such a situation is realized, for example, if
\begin{itemize}
\item $u^{(1)} \in (\xi_1 -\varepsilon ,\xi_1 +\varepsilon )$ and
 $u^{(2)} \in (\xi_3 -\varepsilon ,\xi_3 +\varepsilon )$, and
\item $\pi (s^{(2)})$ is near $p^{2k+1}_{\pm}$ (recall $\pi : [-1/2,1/2]^{2k-1} \twoheadrightarrow S^{2k-1}$).
\end{itemize}
This is the case of $(x_1 ,x_2 ) \in N_1 \times N_3$ in the Lemma.
\end{proof}

Thus it suffices to consider $x\in C_{H_1}(S_{\psi})$ such that both $(x_1 ,x_2 )$ and $(x_3 ,x_4 )$ satisfy the
condition of Lemma \ref{lem_local_N1}, and all the $(s^{(p)},t^{(p)})\in \R^{4k-2}$ are close to each other.
We divide such configurations into two types.

\noindent{\bf Type I.}
All the four points $S_{\psi}(x_p )$ ($1\le p\le 4$) are near the resolution of a single $z_i$ ($i=1$ or $2$).

\begin{figure}[htb]
\[
 \begin{xy}
 {\ar (0,0);(80,0)},{\ar (0,0);(20,30)},{\ar (0,0);(0,25)}, (0,30)*{\R^1},(12,30)*{\R^{2k-1}},(85,3)*{\R^{2k-1}},
 {\ar@{-}(25,10);(37,28)},{\ar@{-}(65,10);(77,28)},{\ar@{.}(25,10);(65,10)},{\ar@{.}(37,28);(77,28)},
 {\ar@{-}(25,6);(27,9)},{\ar@{.}(28,10.5);(37,24)},{\ar@{-}(65,6);(77,24)}, {\ar@{.}(25,6);(65,6)},
  {\ar@{.}(37,24);(77,24)},
 (33,28)*{N_3},(69,6)*{N_1},
 (33,22)*{\bullet},(29.6,12.9)*{\bullet},(30,22)*{x_4},(28,19)*{x_3},(34,14)*{x_2},(32,11.4)*{x_1},
 {\ar (31,19);(31,14)},
 {\ar@{.}(21,0);(25,6)},{\ar@{.}(61,0);(65,6)},
 (20,2)*{q_+},(60,2)*{q_-}
 \end{xy}
\]
\caption{$x=\varphi^{-1}_{H_1}(v_1 ,v_2 ,w)\in N_{1,3,1,3}$}\label{fig_unlink}
\end{figure}
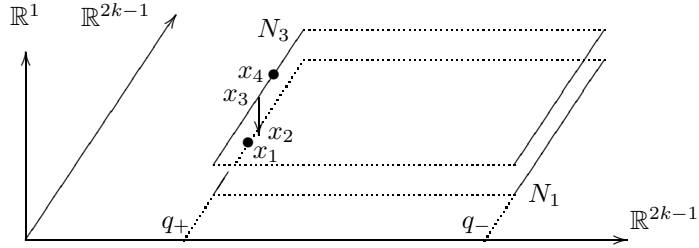

Let us write $N_{p,q,r,s}:=N_p \times N_q \times N_r \times N_s \subset C_4 (\R^{4k-1})$ ($p,q,r,s=1,2,3,4$).
Type I configuration $x=(x_1 ,\dots ,x_4 )$ is in $N_{i, i+2,i,i+2}$ or $N_{i,i+2,i+2,i}$ for some $i$ (modulo $4$).
Figure \ref{fig_unlink} shows an example of $x \in N_{1,3,1,3}$ such that $\varphi_{H_1}(x) =(v_1 ,v_2 ,w)$
($\varphi_{H_1}:=\varphi^{\theta}_{12}\times \varphi^{\theta}_{34}\times \varphi^{\eta}_{23}$) is in the support of
the volume forms.

\noindent{\bf Type II.}
Two points $S_{\psi}(x_1 )$ and $S_{\psi}(x_2 )$ are near the resolution of $z_i$ of $f$, and
$S_{\psi}(x_3 )$ and $S_{\psi}(x_4 )$ are near the resolution of $z_{i+1}$, $i=1$ or $2$ (here $z_3 :=z_1$).

\begin{figure}[htb]
\[
 \begin{xy}
 {\ar (0,0);(80,0)},{\ar (0,0);(20,30)},{\ar (0,0);(0,25)}, (0,30)*{\R^1},(12,30)*{\R^{2k-1}},(85,3)*{\R^{2k-1}},
 {\ar@{-}(23,7);(37,28)},(32,26)*{N_4},
 {\ar@{-}(16,15);(27,15)},{\ar@{-}(29.5,15);(70,15)},(68,18)*{N_3},
 {\ar@{-}(23,3);(30.6,14.4)},{\ar@{-}(31.6,15.9);(37,24)},(39,22)*{N_2},
 {\ar@{-}(16,9);(23,9)},{\ar@{.}(25,9);(27,9)},{\ar@{-}(28,9);(70,9)},(68,6)*{N_1},
 (31,19)*{\bullet},(33,15)*{\bullet},{\ar (33,15);(31,19)},(28,19)*{x_3}, (36,17)*{x_2},
 (29,12)*{\bullet},(33,9)*{\bullet}, (35,7)*{x_1},
 (29,4)*{x_4}="A", {\ar@{-}@/_1mm/"A";(29,12)},
 {\ar@{.}(6,9);(16,9)}, (3,9)*{q_{\pm}},
 {\ar@{.}(21,0);(23,3)}, (19,2)*{q_{\pm}},
 \end{xy}
\]
\caption{$x=\varphi^{-1}_{H_1}(v_1 ,v_2 ,w)\in N_{1,3,4,2}$}\label{fig_link}
\end{figure}
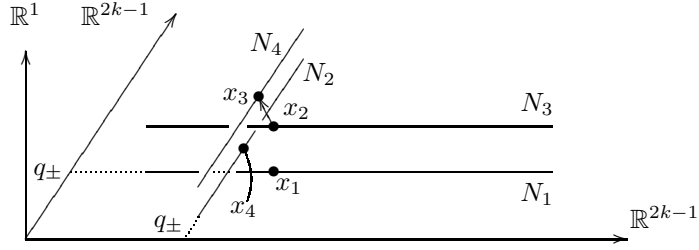

Such a configuration is in $N_{i,i+2,i+1,i+3}$ or $N_{i,i+2,i+3,i+1}$ for some $i$.
Type II inverse image of $(v_1 ,v_2 ,w) \in (S^{6k-1})^2 \times S^{4k-2}$ via $\varphi_{H_1}$ looks like
Figure \ref{fig_link}, because of the following.

\begin{lem}\label{lem_local_N2}
Let $v \in S^{6k-1}$ near $p^{6k-1}_{\pm}$ be given.
If the pair $(x_1 ,x_2 ) \in N_1 \times N_3$ satisfies $\varphi^{\theta}_{12}(x_1 ,x_2 )=v$, then $s^{(2)}$ and
$u^{(2)}$ are uniquely determined by $v$.
If moreover we give any $t^{(2)}$ to fix $x_2$, then accordingly $x_1$ is determined.

Similarly, if $(x_1 ,x_2 ) \in N_2 \times N_4$ satisfies $\varphi^{\theta}_{12}(x_1 ,x_2 )=v$, then
$(x_1 ,x_2 )$ is uniquely determined according to given $s^{(2)}$.
\end{lem}

\begin{proof}
Consider the case $(x_1 ,x_2 ) \in N_1 \times N_3$.
The last $2k+2$ factors of $\varphi^{\theta}_{12}(x_1 ,x_2 )$ are determined by
$\psi (s^{(2)}, t^{(2)})(u^{(2)})-\psi (s^{(1)}, t^{(1)})(u^{(1)})$ as we see in (\ref{eq_direction}).
Comparing them with those of $v$, we can see that $s^{(2)}$, $u^{(1)}$ and $u^{(2)}$ are uniquely determined
so that $s^{(2)}$ is near $q_+$ (see also Lemma \ref{lem_linking_number}).

Thus we can fix $x_2$ if we give any $t^{(2)}$.
Then $x_1$ is also uniquely determined, since $u^{(1)}$ is already determined as above, and $(s^{(1)},t^{(1)})$ is
determined by comparing the first $4k-1$ factors of $v$ with those of $\varphi^{\theta}_{12}(x_1 ,x_2 )$.
\end{proof}

Figure \ref{fig_link} shows an example of $x=(x_1 ,\dots ,x_4) \in N_{1,3,4,2}$.
In this case $x_1,\dots ,x_4$ are uniquely determined by $v_1 =\varphi^{\theta}_{12}(x)$ and
$v_2 =\varphi^{\theta}_{34}(x)$ up to $t^{(2)}$ and $s^{(4)}$ by Lemma \ref{lem_local_N2}, and $t^{(2)}$ and $s^{(4)}$
are determined by $w=(x_3 -x_2 )/\abs{x_3 -x_2}$.

Firstly we compute the contribution of Type II configurations.
Choose neighborhoods $U^{N-1}_{\pm}$ of $p^{N-1}_{\pm}$ ($N=6k$ or $4k-1$) so that the support of $vol_{S^{N-1}}$ is
contained in $U^{N-1}_+ \sqcup U^{N-1}_-$.
Also choose neighborhoods $V^{2k-1}_{\pm}$ of $q_{\pm}$ so that all the four points $(x_1 ,\dots ,x_4 )$ of Type II
configurations are mapped into the same $V^{2k-1}_l \times V^{2k-1}_m$ ($l,m=\pm$) by $pr_{4k-2}:\R^{4k-1}\to \R^{4k-2}$.
Define
\[
 N^{l,m}_i := V^{2k-1}_{l} \times V^{2k-1}_{m} \times (\xi_i -\varepsilon ,\, \xi_i +\varepsilon )
 \subset N_i ,\quad l,m =\pm .
\]
Define $N^{l,m}_{p,q,r,s}:=N^{l,m}_p \times N^{l,m}_q \times N^{l,m}_r \times N^{l,m}_s$ ($l,m=\pm$,
$\{ p,q,r,s \} =\{ 1,2,3,4\}$).
Type II configurations are in $N^{l,m}_{i,i+2,i+1,i+3}$ or $N^{l,m}_{i,i+2,i+3,i+1}$ for some $i=1,\dots ,4$ and
$l,m=\pm$.
There are $4\times 4\times 2=32$ such components.
Via the direction map $\varphi_{H_1}$, each component is mapped to some
$U^{6k-1}_l \times U^{6k-1}_{l'} \times U^{4k-2}_{l''}$ homeomorphically.
Hence each component contributes to $\pair{I(H_1 )}{S_{\psi}}$ by $\pm (1/2)^3 =\pm 1/8$.

The signs are given as follows.

\begin{lem}\label{lem_signs}
\begin{enumerate}
\item
 All the signs for $N^{l,m}_{p,q,r,s}$, $l,m=\pm$, are same for fixed $p,q,r,s$.
\item
 If the orientations are chosen so that the sign for $N^{l,m}_{1,3,2,4}$ is $+1$, then the signs for $N^{l,m}_{2,4,1,3}$
 and $N^{l,m}_{3,1,4,2}$ are $-1$, and those for other five components $N^{l,m}_{2,4,3,1}$, $N^{l,m}_{3,1,4,2}$,
 $N^{l,m}_{4,2,1,3}$, $N^{l,m}_{1,3,4,2}$ and $N^{l,m}_{4,2,3,1}$ are $+1$.
\end{enumerate}
\end{lem}

\begin{proof}
(1) is because the difference between the integrations of $\omega_{H_1}$ over $N^{l,m}_{p,q,r,s}$ and
$N^{-l,m}_{p,q,r,s}$ comes from the antipodal map of $S^{2k-1}$, which preserves the orientation.

For (2), first we show that the signs for $N^{l,m}_{1,3,2,4}$ and $N^{l,m}_{3,1,4,2}$ are different.
This follows from the diagram which is commutative up to isotopy
\[
 \xymatrix{
  N^{l,m}_{1,3,2,4} \ar[d]_-G \ar[r]^-{\varphi_{H_1}} & U^{6k-1}_l \times U^{6k-1}_m \times U^{4k-2}_-
   \ar[d]^-{\iota_{S^{6k-1}} \times \iota_{S^{6k-1}} \times g} \\
  N^{l,m}_{2,4,1,3} \ar[r]^-{\varphi_{H_1}} & U^{6k-1}_l \times U^{6k-1}_m \times U^{4k-2}_-
 }
\]
where
\begin{multline*}
 G(x_1 ,\dots ,x_4 ) := (t^{(1)},s^{(1)},u^{(1)}-\xi_1 +\xi_2 \, ; \, t^{(2)},s^{(2)},u^{(2)}-\xi_3 +\xi_4 \, ; \\
 t^{(3)},s^{(3)},u^{(3)}-\xi_3 +\xi_1 \, ; \, t^{(4)},s^{(4)},u^{(4)}-\xi_4 +\xi_3 ),
\end{multline*}
and $g$ is induced from $(s,t,u)\mapsto (s,t,-u)\in (\R^{2k-1})^2 \times \R^1$.
$G$ preserves the orientation but $\iota_{S^{6k-1}} \times g \times \iota_{S^{6k-1}}$ does not, so the signs for
$N^{l,m}_{1,3,2,4}$ and $N^{l,m}_{3,1,4,2}$ are different.
The signs for $N^{l,m}_{2,4,1,3}$ and $N^{l,m}_{3,1,4,2}$ are same, since their difference comes from an automorphism
of $H_1$, which always preserves the orientation.

Next we show that all the signs for other five components are same as that for $N^{l,m}_{1,3,2,4}$.

\noindent
\underline{(i) $N^{l,m}_{1,3,2,4}$ and $N^{l,m}_{4,2,3,1}$};
the signs for them are same, since the difference arises from an automorphism of $H_1$.

\noindent
\underline{(ii) $N^{l,m}_{1,3,2,4}$ and $N^{l,m}_{2,4,3,1}$};
their signs are same because of the following commutative diagram (up to isotopy);
\[
 \xymatrix{
  N^{l,m}_{1,3,2,4} \ar[d]_-{G'} \ar[r]^-{\varphi_{H_1}} & U^{6k-1}_l \times U^{6k-1}_m \times U^{4k-2}_-
   \ar[d]^-{i_1 \times \id \times g} \\
  N^{l,m}_{2,4,3,1} \ar[r]^-{\varphi_{H_1}} & U^{6k-1}_l \times U^{6k-1}_m \times U^{4k-2}_+
 }
\]
where
\begin{multline*}
 G'(x_1 ,\dots ,x_4 ) := (s^{(2)},t^{(1)},u^{(2)} \, ; \, s^{(1)},t^{(2)},u^{(1)} ; \\
 s^{(3)}-s^{(2)}+s^{(1)},t^{(3)},u^{(3)} \, ; \, s^{(4)}-s^{(2)}+s^{(1)},t^{(4)},u^{(4)}),
\end{multline*}
and $i_1$ was given in Lemma \ref{lem_outside_zero}.
The vertical arrows preserve the orientations, so the signs are same.
Similar diagrams shows that the signs for $N^{l,m}_{1,3,2,4}$ and $N^{l,m}_{3,1,2,4}$ are same.

\noindent
\underline{(iii) $N^{l,m}_{4,2,1,3}$ and $N^{l,m}_{1,3,4,2}$};
their signs are same as those of $N^{l,m}_{3,1,2,4}$ and $N^{l,m}_{1,3,2,4}$ respectively, since their differences
come from automorphisms of $H_1$.
\end{proof}

By Lemma \ref{lem_signs}, 24 components contribute to $\pair{I(H_1 )}{S_{\psi}}$ by $\pm 1$ and other 8 components
by $\mp 1$.
Thus their contribution is $\pm (1/8)\times (24-8) =\pm 2$.

Next consider Type I configurations.

\begin{lem}\label{lem_dim_reason}
The contribution of Type I configurations is zero.
\end{lem}

\begin{proof}
Consider the case when all the four points $S_{\psi}(x_p )$ are near $z_1$.
Then $x\in C_4 ((N_1 \sqcup N_3 )\cap L)$.
But on this space, the direction map $\varphi_{H_1}$ is invariant under the translation
\[
 (x_1 ,\dots ,x_4 ) \mapsto (\tau_v x_1 ,\dots ,\tau_v x_4 ),\quad v\in \R^{2k-1}
\]
where $\tau_v (s,t,u):=(s,t+v,u)\in \R^{2k-1}\times \R^{2k-1}\times \R^1$.
This is because of \eqref{eq_direction} and $\psi (s,t)(u)=\psi (s,t+v)(u)$ if $(s,t,u)\in (N_1 \sqcup N_3 )\cap L$.
Hence the image of $\varphi_{H_1}$ must be of positive codimension and the integrand
$\omega_{H_1}=\varphi^*_{H_1} (vol_{S^{6k-1}}\times vol_{S^{6k-1}}\times vol_{S^{4k-2}})$ must vanish on
$C_4 ((N_1 \sqcup N_3 )\cap L)$.

In the case when four points $S_{\psi}(x_p )$ are near $z_2$, a similar translation $\tau'_v (s,t,u):=(s+v,t,u)$
would prove the vanishing of $\omega_{H_1}$ on $C_4 ((N_2 \sqcup N_4 )\cap L)$.
\end{proof}

Thus the first term of \eqref{eq_evaluation} is equal to $\pm 2 \times (1/2)=\pm 1$.

\subsubsection{Remaining terms}\label{sssec_remaining}
To complete the proof of Theorem \ref{thm_Haefliger}, we will prove that the second and the third terms of
\eqref{eq_evaluation} do not contribute to $\pair{\calH}{S_{\psi_{\delta}}}$.

\noindent
{\bf Contribution of $H_2$.}
Now we compute the second term of \eqref{eq_evaluation}.
Recall $C_{H_2}(S_{\psi}) \subset C^o_3 (\R^{4k-1})\times \R^{6k}$.

\begin{lem}\label{lem_ev_H2}
$\pair{I(H_2 )}{S_{\psi_{\delta}}}=O(\delta )$.
\end{lem}

\begin{proof}
Divide $\pair{I(H_2 )}{S_{\psi}}$ into two integrations;
the integration over
\[
 R_{\ge 1}:=\{ (x;y) \in C_{H_2}(S_{\psi}) \, | \, y \not\in \R^{4k}\times (-1,1)^{2k}\}
\]
($x=(x_1 ,x_2 ,x_3 )$) and that over $R_{<1}:=C_{H_2}(S_{\psi}) \setminus R_{\ge 1}$.

The integration over $R_{\ge 1}$ is well defined and continuous at $\delta =0$, since $y\in \R^{6k}$ is
far from the image of $S_{\psi_{\delta}}$ (see \eqref{eq_thin}).
But when $\delta =0$, this integral is zero, since all the three points $S_{\psi}(x_p )$ ($p=1,2,3$) are in
$\R^{4k}\times \{ 0 \}^{2k}$ and hence the image of the direction map
$\varphi_{H_2}:=\varphi^{\theta}_{14}\times \varphi^{\theta}_{24}\times \varphi^{\theta}_{34}$ is of positive codimension
$\ge 2k-1$ in $(S^{6k-1})^3$ (see Lemma \ref{lem_localize2}).
Hence the integration over $R_{\ge 1}$ is $O(\delta )$.

Next consider the integration over $R_{<1}$.
We have only to consider $(x;y)$ with the first $6k-1$ factors of $S_{\psi}(x_1 ),S_{\psi}(x_2 ),S_{\psi}(x_3 )$ and $y$
close to each other; otherwise the image of $\varphi_{H_2}$ cannot be in $(\text{supp}(vol))^3$ (see condition (2) just
before Figure \ref{fig_contribution}).
In particular all $(s^{(p)},t^{(p)})\in \R^{4k-2}$ ($p=1,2,3$) must be close to each other.

Consider the case that $x=(x_1 ,x_2 ,x_3 )\in C_3 (\R^{4k-1}\setminus L)$.
Similarly as in Lemma \ref{lem_outside_zero}, if we choose $vol_{S^{6k-1}}$ with sufficiently small support, then we
may assume that at most one of $x\cap (X\setminus L)$ and $x\cap (Y\setminus L)$ is non-empty.
We can prove the vanishing of the integration of $\omega_{H_2}$ over $\{ (x;y) \,| \, x \cap X \ne \emptyset \}$ in a
similar way as in Lemma \ref{lem_outside_zero} by considering an involution
$F_1 : (x_1 ,x_2 ,x_3 ;y)\mapsto (i_1 x_1 ,i_1 x_2 ,i_1 x_3 ;i_1 y)$ ($i_1$ was defined in Lemma
\ref{lem_outside_zero}), which preserves the orientation but satisfies $F^*_1 \omega_{H_2}=-\omega_{H_2}$ because
$\varphi^{\theta}_{p4} \circ F_1 = i_1 \circ \varphi^{\theta}_{p4}$ ($p=1,2,3$) on $X$ by \eqref{eq_psi_symmetry_1} and
$i^*_1 vol_{S^{6k-1}}=-vol_{S^{6k-1}}$.
The vanishing of the integrations over $\{ (x;y) \,| \, x \cap Y \ne \emptyset \}$ and
$\{ x\cap X = x\cap Y=\emptyset \}$ can also be proved by similar involutions $F_i$, $i=1,2$.

So we may assume that one of $x_p$ is in $L$ (and other two points are near $L$).
Since the first $6k-1$ factors of $S_{\psi}(x_p )$ ($p=1,2,3$) are close to each other, an analogous argument to the
proof of Lemma \ref{lem_localize2} shows that only the integration over the subspace of $(x;y)$ with $x$ in
$C_3 (N_1 \sqcup N_3 )$ or $C_3 (N_2 \sqcup N_4 )$ or $C_3 (N'_1 \sqcup N'_2 )$, where
$N'_i :=\R^{4k-2}\times (\zeta_i -\varepsilon ,\zeta_i +\varepsilon )$ ($i=1,2$) correspond to the crossing $c$ of $f$
(see \S \ref{sssec_choice_psi} and Figure \ref{fig_f}); otherwise two or more $S_{\psi}(x_p )$'s are in
$\R^{4k}\times \{ 0\}^{2k}$ and hence the image of the map $\varphi_{H_2}$ is of positive codimension $\ge 2k-1$ in
$(S^{6k-1})^3$.
But similarly as in Lemma \ref{lem_dim_reason}, on these spaces we can define translations $\tau$, $\tau'$ under which
$\varphi_{H_2}$ is invariant, and hence the integrand $\omega_{H_2}$ must vanish by dimensional reason.
\end{proof}

\noindent
{\bf Contribution of the correction term $c$.}
Lastly we compute the third term of \eqref{eq_evaluation}.
This is an integration of $D(S_{\psi})^* \mu$ over $\R^{4k-1}$ (cf.\ Lemma \ref{lem_localize3}).
See Lemma \ref{lem_hatomega_closed} and Definition \ref{def_c} for the definition of $c$.

\begin{lem}\label{lem_ev_c}
$\pair{c}{S_{\psi}}=0$.
\end{lem}

\begin{proof}
First we show that integrations of $D(S_{\psi})^* \mu$ over $X$, $Y$ and $Z$ (see \S \ref{sssec_choice_psi}) vanish.
For $X$, this is because
\begin{itemize}
\item
 $i^*_1 D(S_{\psi})^* \mu =D(S_{\psi})^* \mu$ on $X$, since $D(S_{\psi})\circ i_1 =i_1 \circ D(S_{\psi})$ on $X$ by
 \eqref{eq_psi_symmetry_1} (where $i_1 : \I_{4k-1} (\R^{6k})\to \I_{4k-1} (\R^{6k})$ is given by $f\mapsto i_1 \circ f$;
 see Remark \ref{rem_mu_symmetric}) and we can choose $\mu$ so that $i^*_1 \mu =\mu$ (see Remark \ref{rem_mu_symmetric}),
 and
\item
 $i_1$ is an orientation reversing diffeomorphism of $X$.
\end{itemize}
Similar arguments hold for $Y$ and $Z$.

So we may restrict the integration to $[-1/2,1/2]^{4k-2}\times \R^1$.
If $(s,t,u)\not\in N_3 \sqcup N_4$, then $\psi (s,t)(u)$ does not depend on $(s,t)\in \R^{4k-2}$, so $D(S_{\psi})$ is
invariant under the translations $\tau$ and $\tau'$ defined in the proof of Lemma \ref{lem_dim_reason}.
If $(s,t,u)\in N_3$ (resp.\ $N_4$), then $\psi (s,t)(u)$ does not depend on $t\in \R^{2k-1}$ (resp.\ $s$), so
$D(S_{\psi})$ is invariant under the translations $\tau$ (resp.\ $\tau'$).
Thus  the image of $[-1/2,1/2]^{4k-2}\times \R^1$ via $D(S_{\psi})$ must be of dimension $<4k-1$ and hence
a ($4k-1$)-form $D(S_{\psi})^* \mu$ must vanish on $[-1/2,1/2]^{4k-2}\times \R^1$.
\end{proof}

Thus we have completed the proof of Theorem \ref{thm_Haefliger};
only Type II configurations for $H_1$ contribute to $\pair{\calH}{S_{\psi}}$ by $\pm 1$, and hence
$\pair{\calH}{S_{\psi}}=\pm 1$.

\subsection{Non-triviality of $\calH$ in general dimensions}\label{subsec_H_general}

Here we complete the proof of Theorem \ref{thm_main1}.
Suppose that $n>j\ge 2$, $n-j\ge 3$ is odd and $m:=2n-3j-3>0$.
Put $n-j=2k+1$ ($k\ge 1$) and consider $S_{\psi} \in \emb{6k}{4k-1}$ as above.
Notice that $n=6k-m$ and $j=4k-m-1$, and in particular $4k-m-1>0$.

Since $S_{\psi}$ is of the form (\ref{eq_deform_spinning}), we can find
$l_m \in \Omega^m \emb{6k-m}{4k-m-1}=\Omega^m \emb{n}{j}$ such that $S_{\psi}=\gr^m (l_m )$.
Explicitly we can define $l_m :\R^m \to \emb{n}{j}$ by
\begin{equation}\label{eq_graphing}
\begin{split}
 &l_m (t_1 ,\dots ,t_m)(x_1 ,\dots ,x_j ) := \\
 &\quad ((x_1 ,\dots ,x_{j-1}), \psi(t_1 ,\dots ,t_m , x_1 ,\dots ,x_{j-1})(x_j)) \in \R^{j-1} \times \R^{2k+2}=\R^n
\end{split}
\end{equation}
and regard it as in $\Omega^m \emb{n}{j}$.
We think of $[l_m ]$ as the generator of $H_m (\emb{n}{j})$ via the Hurewicz isomorphism ($\emb{n}{j}$ is
$(2n-3j-4)$-connected; see \cite{Budney08}).

Consider $\calH =[I(H)+c]\in H^m_{DR}(\emb{n}{j})$.
The following completes the proof of Theorem \ref{thm_main1}.

\begin{thm}\label{thm_Haefliger_general}
The Kronecker pairing $\pair{\calH}{l_m}$ is equal to $\pm 1$.
\end{thm}

\begin{proof}
Define the spaces $\hat{C}_{H_i}$, $i=1,2$, by $\hat{C}_{H_1} := \R^m \times C^o_4 (\R^j )$ and
\[
 \hat{C}_{H_2} :=
 \{ (t,(x_1 ,x_2 ,x_3 ),y)\in \R^m \times C^o_3 (\R^j ) \times \R^n \, | \, l_m (t)(x_p )\ne y,\ p=1,2,3 \} .
\]
These spaces are also defined as the following pullback square;
\[
 \xymatrix{
  \hat{C}_{H_i} \ar[r]^-{\hat{l}_m} \ar[d] & C_{H_i} \ar[d]^-{\pi_{H_i}} \\
  \R^m \ar[r]^-{l_m}                       & \emb{n}{j}
 }
\]
Then $\pair{\calH}{l_m}$ is equal to
\begin{equation}\label{eq_ev_general_dim}
 \frac{1}{4}\int_{\hat{C}_{H_1}} \hat{l}^*_m \omega_{H_1} + \frac{1}{12} \int_{\hat{C}_{H_2}} \hat{l}^*_m \omega_{H_2}
 +\int_{\R^m \times \R^j}D(l_m )^* \mu ,
\end{equation}
where $D(l_m ):\R^m \times \R^j \to \I_j (\R^n )$ is defined by $(t,x)\mapsto d(l_m (t))_x$.

Recall $C_{H_1}(S_{\psi})=C^o_4 (\R^{4k-1})$ and $C_{H_2}(S_{\psi})\subset C^o_3 (\R^{4k-1})\times \R^{6k}$.
We regard $\hat{C}_{H_i} \subset C_{H_i}(S_{\psi})$ ($i=1,2$);
\begin{align*}
 \hat{C}_{H_1} &\cong \{ (x_1 ,\dots ,x_4 ) \in C_{H_1}(S_{\psi}) \, |\, pr_m (x_1 )=\dots =pr_m (x_4 )\} , \\
 \hat{C}_{H_2} &\cong \{ (x_1 ,x_2 ,x_3 ;y) \in C_{H_2}(S_{\psi}) \, |\, pr_m (x_1 )=pr_m (x_2 )=pr_m (x_3 )=pr_m (y) \}
\end{align*}
($pr_m :\R^N \to \R^m$ ($N=4k-1$ or $6k$) is the first $m$ projection), via diffeomorphisms given by respectively
\[
 (t,x)   \longmapsto ((t,x_1 ),\dots ,(t,x_4 )), \quad
 (t,x,y) \longmapsto ((t,x_1),(t,x_2 ),(t,x_3 ),(t,y)).
\]
The direction maps $C_{H_i} \to S^{N-1}$, $N=n$ or $j$, composed by $\hat{l}_m$ are regarded as
\begin{equation}\label{eq_direction_proj}
 \begin{split}
 &\varphi^{\theta}_{12} : \hat{C}_{H_1}\longrightarrow S^{n-1}=S^{6k-1}\cap \{ x\in \R^{6k} \, ; \, pr_m (x) =0 \}, \\
 &\varphi^{\theta}_{12}((t,x_1 ),(t,x_2 ))=
 \frac{(x_2 -x_1 ,S_{\psi}(t, x_2 )-S_{\psi}(t, x_1 ))}{\abs{(x_2 -x_1 ,S_{\psi}(t, x_2 )-S_{\psi}(t, x_1 ))}},
 \end{split}
\end{equation}
and so on.
Then the integrations relating to $L\subset \R^m \times \R^j =\R^{4k-1}$ (a neighborhood of
$[-1/2 ,1/2]^m \times [-1/2,1/2]^{j-1}\times \R^1$) can be computed in similar ways as in the previous subsection;

\noindent{\bf The first term.}
Type I contribution vanishes by the translations $\tau$ or $\tau'$ as in Lemma \ref{lem_dim_reason}.
Type II configurations contribute by $\pm 2$; each component
\begin{align*}
 &(N^{l,m}_{i,i+2,i+1,i+3}) \cap \{ pr_m (x_1 )=\dots =pr_m (x_4 )\} \ \text{ or} \\
 &(N^{l,m}_{i,i+2,i+3,i+1}) \cap \{ pr_m (x_1 )=\dots =pr_m (x_4 )\}
\end{align*}
for some $i=1,\dots ,4$ and $l,m=\pm$ is mapped via the direction map to some
\[
 U^{6k-1}_l \times U^{6k-1}_{l'} \times U^{4k-2}_{l''}\cap \{ \text{first } m \text{ projections}=0\} ,
\]
which is $U^{n-1}_l \times U^{n-1}_{l'} \times U^{j-1}_{l''}$ (see \S \ref{sssec_first_term}).
The sign arguments are slightly different, but the result is same as Lemma \ref{lem_signs};
in the first diagram in the proof of Lemma \ref{lem_signs}, the map $G$ restricted to $\hat{C}_{H_1}$ preserve the
orientation but the left vertical map does not.
Both the vertical maps in the second diagram restricted to $\hat{C}_{H_1}$ have the same orientation sign $(-1)^a$,
where $a=\min \{ m,2k-1\}$.

\noindent{\bf The second term.}
Similarly as in Lemma \ref{lem_ev_H2}, the integration over
$R_{\ge 1}:=\{ (t,x,y)\in \hat{C}_{H_2}\, | \, y\not\in \R^{j+1}\times (-1,1)^{2k}\}$ (notice $j+1+2k=n$) is
$O(\delta )$.
In $R_{<1}:=\hat{C}_{H_2}\setminus R_{\ge 1}$, only $(t,x,y)$ with $pr_{j-1}(x_p )$ ($p=1,2,3$) close to each other and
$x \in C_3 ((N_i \sqcup N_{i+2})\cap L)$ ($i=1,2$) or $C_3((N'_1 \sqcup N'_2 )\cap L)$ may contribute to the integral
(where $N_i := \R^m \times \R^{j-1}\times (\xi_i -\varepsilon ,\xi_i +\varepsilon )$,
$N'_i := \R^m \times \R^{j-1}\times (\zeta_i -\varepsilon ,\zeta_i +\varepsilon )$); otherwise the image of
$\varphi_{H_2}$ is of positive codimension.
But on these spaces the direction map is invariant under the similar translations $\tau$ or $\tau'$ to those in
Lemma \ref{lem_dim_reason} and the integrand vanishes by dimensional reason.

\noindent{\bf The third term.}
Similarly as in Lemma \ref{lem_ev_c}, since the derivation map $D(l_m )$ is invariant on $L$ under the translations
given by using $\tau$ or $\tau'$, the integration over $L$ vanishes by dimensional reason.

The signs appearing in the proof of vanishing of integrations over $\R^{4k-1}\setminus L$ are slightly different from
those in the previous subsection.
Recall the involutions $F_1$ and $F_2$ on ${C}_{H_i}(S_{\psi})$.
They clearly preserve $\hat{C}_{H_i}$, and are given by
\begin{alignat*}{2}
 F_l ((t,x_1 ),\dots ,(t,x_4 ))&=(i_l (t,x_1 ),\dots ,i_l (t,x_4 )) &\quad &\text{on } \hat{C}_{H_1}, \\
 F_l ((t,x_1 ),(t,x_2 ),(t,x_3 );(t,y))&=(i_l (t,x_1 ),i_l (t,x_2 ),i_l (t,x_3 );i_l (t,y))&\quad &\text{on }\hat{C}_{H_2}
\end{alignat*}
($l=1,2$), here we regard $\hat{C}_{H_i}$ as a subspace of $C_{H_i}(S_{\psi})$ as above.
The maps $i_1 ,i_2:\R^N\to \R^N$ ($N=4k-1$ or $6k$) are as given in Lemma \ref{lem_outside_zero}.
For $N=4k-1$, they are written explicitly as
\begin{align*}
 i_1 (t,x) &=
  \begin{cases}
   (-t_1 ,\dots ,-t_m ;-x_1 ,\dots ,-x_{2k-1-m},x_{2k-m},\dots ,x_j ) & m\le 2k-1, \\
   (-t_1 ,\dots ,-t_{2k-1}, t_{2k},\dots ,t_m ;x_1 ,\dots ,x_j )      & m>2k-1,
  \end{cases} \\
 i_2 (t,x) &=
  \begin{cases}
   (t_1 ,\dots ,t_m ;x_1 ,\dots ,x_{2k-1-m},-x_{2k-m},\dots ,-x_{j-1}, x_j ) & m\le 2k-1, \\
   (t_1 ,\dots ,t_{2k-1}, -t_{2k},\dots ,-t_m ;-x_1 ,\dots -x_{j-1}, x_j )   & m>2k-1.
  \end{cases}
\end{align*}
First consider the first and the second terms of \eqref{eq_ev_general_dim}.
Let $X,Y,Z\subset \R^m \times \R^j =\R^{4k-1}$ be subsets defined similarly as in \S \ref{sssec_choice_psi}, and set
\begin{align*}
 \hat{C}_{H_1}(X) &:= \{ (t,(x_1 ,\dots ,x_4 ))\in \hat{C}_{H_1} \, | \, (t,x_p )\in X,\ \forall p \} , \\
 \hat{C}_{H_2}(X) &:= \{ (t,(x_1 ,x_2 ,x_3 ),y)\in \hat{C}_{H_2} \, | \, (t,x_p )\in X,\ \forall p \}
\end{align*}
and so on.
Then the actions of $F_1$ and $F_2$ on the forms $\omega_{H_i}$ and orientations of the spaces are described as in
Table \ref{table:signs}, which is a consequence of the equations
\begin{align*}
 \varphi \circ F_1 &=
  \begin{cases}
   i_{1,2k-m-1} \circ \varphi & m\le 2k-1 \\
   \varphi                    & m>2k-1
  \end{cases}
 \quad \text{on } \hat{C}_{H_i}(X), \\
 \varphi \circ F_2 &=
  \begin{cases}
   i_{2k-m,j-1} \circ \varphi & m\le 2k-1 \\
   i_{1,j-1}    \circ \varphi & m>2k-1
  \end{cases}
 \quad \text{on } \hat{C}_{H_i}(Y),
\end{align*}
where $\varphi$ is one of the direction maps, and $i_{p,q}:\R^n \to \R^n$ ($p<q$) is given by
\[
 i_{p,q} (a_1 ,\dots ,a_n ) :=(a_1 ,\dots ,a_{p-1},-a_p ,\dots ,-a_q , a_{q+1},\dots ,a_n )
\]
(in particular $i_1$ and $i_2$ we have used can be written as $i_1 =i_{1,2k-1}$, $i_2 =i_{2k,4k-2}$).
\begin{table}[ht]
{\begin{tabular}{|c||c|c|}\hline
                                           & $m\le 2k-1$                 & $m>2k-1$ \\
\hline\hline
orientation sign of $F_1$                  & $(-1)^m$                    & $-1$ \\
\hline
$F^*_1 \omega_{H_i}$ on $\hat{C}_{H_i}(X)$ & $(-1)^{2k-1-m}\omega_{H_i}$ & $+\omega_{H_i}$ \\
\hline
orientation sign of $F_2$                  & $+1$                        & $(-1)^{m-2k+1}$ \\
\hline
$F^*_2 \omega_{H_i}$ on $\hat{C}_{H_i}(Y)$ & $-\omega_{H_i}$             & $(-1)^{j-1}\omega_{H_i}$ \\ \hline
\end{tabular}}
\caption{Signs of $F_l$}
\label{table:signs}
\end{table}
Thus the integrations of $\omega_{H_i}$ over $\hat{C}_{H_i}(X)$ and $\hat{C}_{H_i}(Y)$ vanish (when $m>2k-1$, we use
$j=4k-1-m$ and hence $m-2k+j$ is odd).

For the third term of \eqref{eq_ev_general_dim}, we have to study the signs arising from the involutions $i_1$ and $i_2$.
The orientation signs of $i_1 ,i_2 :\R^m \times \R^j \to \R^m \times \R^j$ are always $-1$.
But $i_1,i_2$ always preserve the integrand, because
\begin{align*}
 D(l_m )\circ i_1 &=
 \begin{cases}
  i_{1,2k-m-1} \circ D(l_m ) & \text{if }m\le 2k-1, \\
  D(l_m )                    & \text{if }m>2k-1,
 \end{cases}\quad \text{on }X\cup Z, \\
 D(l_m )\circ i_2 &=
 \begin{cases}
  i_{2k-m,j-1} \circ D(l_m )  & \text{if }m\le 2k-1, \\
  i_{1,j-1} \circ D(l_m ) & \text{if }m>2k-1,
 \end{cases}\quad \text{on }Y\cup Z,
\end{align*}
and we can choose $\mu$ so that $i'_1$, $i'_2$ and $i''_2$ preserve $\mu$ (see Remark \ref{rem_mu_symmetric}).
Hence the integrations of $\mu$ over $X$, $Y$ and $Z$ vanish.
\end{proof}

\section{Vanishing results}\label{sec_vanish}

In \S\S \ref{subsection_boundary_survey}, ..., \ref{subsection_hidden_infty} we prove Theorem \ref{cochain} assuming
$n-j \ge 2$ is even, by studying the boundary strata of compactified configuration spaces.
Some results here hold even if $n-j\ge 3$ is odd and can be used to prove Theorem \ref{thm_closed}
(see \S \ref{subsec_proof_closed}).
In \S \ref{subsec_indep_vol} we complete the proof of Proposition \ref{prop_indep_vol}.

\subsection{Boundary strata}\label{subsection_boundary_survey}

Let $\Gamma$ be a graph with $s$ i-vertices and $t$ e-vertices.
We denote by $C_{s,t}$ the fiber of $\pi_{\Gamma}$.
The compactified configuration space $C_{s,t}$ is a manifold with corners.
The boundary $\partial C_{s,t}$ consists of configurations where some points in the configuration are allowed to
`collide together.'
Moreover $\partial C_{s,t}$ is stratified via `complexities' of collisions.

But here we do not need the complete description of all strata.
For our purpose only the most `generic' part, the codimension one strata, are needed.
Such strata correspond to `coinstantaneous collisions' of points, and are parametrized by subsets of the set
$V(\Gamma )$ of vertices of $\Gamma$, as we will explain in \S \ref{subsection_boundary_strata}.

\subsection{Codimension one strata}\label{subsection_boundary_strata}

To any subset $A \subset V(\Gamma )$ with $\sharp A \ge 2$, a codimension one stratum $C_A \subset \partial C_{s,t}$
is assigned.
Namely $C_A$ consists of configurations where the points labeled by $A$ `simultaneously collide together.'
More precisely, any point in $C_A$ can be written as a limit point
\begin{equation}\label{first_limit}
 \lim_{\tau \to 0} (x_1 (\tau ),\dots ,x_s (\tau ),y_{s+1} (\tau ),\dots ,y_{s+t} (\tau ))
\end{equation}
such that $(x_1 (\tau ),\dots ,y_{s+t} (\tau )) \in C^o_{s,t}$ for $\tau >0$, and $x_p (\tau )$, $y_q (\tau )$ can be
written as
\[
 x_p (\tau ) =
 \begin{cases}
  x_p \quad (\text{constant}) & \text{if } p \not\in A, \\
  z + \tau v_p & \text{if } p \in A,
 \end{cases}\quad
 y_q (\tau ) =
 \begin{cases}
  y_q \quad (\text{constant}) & \text{if } q \not\in A, \\
  z + \tau w_q & \text{if } q \in A,
 \end{cases}
\]
for some $v_p \in \R^j \setminus \{ 0\}$, $w_q \in \R^n \setminus \{ 0\}$ and $z \in \R^n$.

There are other types of codimension one strata, denoted by $C^{\infty}_A$, parametrized by all the
non-empty subsets $A \subset V(\Gamma )$.
The stratum $C^{\infty}_A \subset \partial C_{s,t}$ consists of configurations where the points labeled by $A$
`escape to infinity.'
More precisely $C^{\infty}_A$ is the set of limit points
\begin{equation}\label{infty_limit}
 \lim_{\tau \to \infty} (x_1 (\tau ),\dots ,x_s (\tau ),y_{s+1} (\tau ),\dots ,y_{s+t} (\tau ))
\end{equation}
where $(x_1 (\tau ),\dots ,y_{s+t} (\tau )) \in C^o_{s,t}$ ($0 <\tau <\infty$) is of the form
\[
 x_p (\tau ) =
 \begin{cases}
  x_p \quad (\text{constant}) & \text{if } p \not\in A, \\
  x_p + \tau v_p              & \text{if } p \in A,
 \end{cases}\quad
 y_q (\tau ) =
 \begin{cases}
  y_q \quad (\text{constant}) & \text{if } q \not\in A, \\
  y_q + \tau w_q              & \text{if } q \in A,
 \end{cases}
\]
for some $v_p \in \R^j \setminus \{ 0\}$, $w_q \in \R^n \setminus \{ 0\}$, $x_p \in \R^j$ and $y_q \in \R^n$.

All the codimension one strata is of the form $C_A$ or $C^{\infty}_A$, hence we have
\[
 \partial C_{s,t}= \left( \bigcup_{A \subset V(\Gamma ), \ \sharp A \ge 2}\overline{C_A} \right) \cup
 \left( \bigcup_{A' \subset V(\Gamma ), \ \sharp A' \ge 1}\overline{C^{\infty}_{A'}} \right) .
\]
We will call $C_A$ ($\sharp A \ge 2$) a {\it stratum of non-infinity type}, and $C^{\infty}_A$
($A \ne \emptyset$) a {\it stratum at infinity}.

Define the subsets $\Sigma_A$ and $\Sigma^{\infty}_A$ of $C_{\Gamma}$ by
\[
 \Sigma_A :=\bigcup_{f \in \emb{n}{j}}C_A (f) , \quad \Sigma^{\infty}_A :=\bigcup_{f \in \emb{n}{j}}C^{\infty}_A (f)
\]
where $C_A (f)$ and $C^{\infty}_A (f)$ are codimension one strata of $\pi_{\Gamma}^{-1}(f)$ described as above.
Then $\Sigma_A$ and $\Sigma^{\infty}_A$ fiber over $\emb{n}{j}$.

Notice that it is enough to describe $\Int \Sigma^{(\infty )}_A$ for the proof of Theorem \ref{cochain}.
In \S \ref{subsection_strata_1} and \S \ref{subsection_strata_infinity} we will describe these strata explicitly,
following \cite{BottTaubes94, CCL02, Rossi_thesis, Watanabe07}.

\subsection{Explicit description of non-infinity type strata}\label{subsection_strata_1}

Let $A \subset V(\Gamma )$ be a subset with $\sharp A \ge 2$
(recall that $V(\Gamma )$ denotes the set of vertices of a graph $\Gamma$).
Here we study the strata $\Sigma_A$ of non-infinity type.
Denote by $E(\Gamma )$ the set of edges of $\Gamma$.

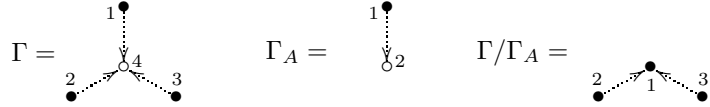
\begin{figure}[htb]
\[
 \begin{xy}
 (0,0)*{\circ}="O",(0,8)*{\bullet}="A",(-6.928,-4)*{\bullet}="B",(6.928,-4)*{\bullet}="C",
 {\ar@{.>}"A";"O"_<{1}^>{4}},{\ar@{.>}"B";"O"^<{2}},{\ar@{.>}"C";"O"_<{3}},(-12,2)*{\Gamma =},
 (35,0)*{\circ}="o",(35,8)*{\bullet}="a", {\ar@{.>}"a";"o"_<{1}^>{2}},(23,2)*{\Gamma_A =},
 (70,0)*{\bullet}="p",(63.072,-4)*{\bullet}="q",(76.928,-4)*{\bullet}="r",
 {\ar@{.>}"q";"p"^<{2}_>{1}}, {\ar@{.>}"r";"p"_<{3}},(53,2)*{\Gamma /\Gamma_A =},
 \end{xy}
\]
\caption{Examples of $\Gamma_A$ and $\Gamma /\Gamma_A$ for $A=\{ 1,4 \}$ ($n$ odd)}\label{fig_subgraph_noninfty}
\end{figure}

\begin{defn}\label{subgraph}
The {\em subgraph} $\Gamma_A$ of $\Gamma$ associated with $A$ is a (possibly non-admissible) graph with
$V(\Gamma_A )=A$ and $E(\Gamma_A ) = \{ pq \in E(\Gamma )\, |\, p,q \in A,\ p \ne q \}$ (hence small loops are ignored).
If $A=\{ p_1 <\dots <p_k \}$, then the vertex of $\Gamma_A$ which was labeled by $p_a$ in $\Gamma$ is re-labeled
by $a$.
The labels of edges are suitably decreased.

The {\it quotient graph} $\Gamma / \Gamma_A$ is a graph obtained by `collapsing $\Gamma_A$ to a point $v_A$.'
More precisely,
\begin{align*}
 V(\Gamma / \Gamma_A ) &:= (V(\Gamma ) \setminus A) \sqcup \{ v_A \} ,\\
 E(\Gamma / \Gamma_A ) &:= \{ pq \in E(\Gamma ) \, |\, p,q \not \in A \} \sqcup
 \{ pv_A \,|\, p \not\in A,\ pq \in E(\Gamma ) \text{ for }\exists q \in A \} .
\end{align*}
The vertex $v_A$ is internal if there is an i-vertex in $A$, and is external otherwise.
We label $v_A$ by $\min \{ p\, ;\, p\in A\}$, and the labels of other vertices and edges are suitably decreased
(see Figure \ref{fig_subgraph_noninfty} for an example).
\end{defn}

There is a projection
\[
 p_A : \Int \Sigma_A \longrightarrow C^o_{\Gamma / \Gamma_A}
\]
which maps the limit point (\ref{first_limit}) to $(\dots ,x_p , \dots ,z,\dots ,y_q ,\dots )$, $p,q \not\in A$.
The point $z$ corresponds to the vertex $v_A$, hence $z \in \R^n$ if $A$ contains no i-vertex, and $z \in \R^j$ if
there is an i-vertex in $A$.

The fiber of $p_A$ is thought of as the space of `infinitesimal configurations' at the colliding point.
We will define a fibration $\rho_A :\hat{B}_A \to B_A$ which describes such infinitesimal configurations.

\begin{defn}\label{def_I_j}
Define $\I_j (\R^n )$ as the set of all $j$-frames in $\R^n$.
In other words
\[
 \I_j (\R^n ) := \{ \text{linear injective maps } \R^j \hookrightarrow \R^n \} .
\]
We give $\I_j (\R^n )$ a natural structure of an open submanifold of $(\R^n \setminus \{ 0\})^j$.
\end{defn}

Let $a$ and $b$ be the numbers of i- and e-vertices in $A$.
Define a manifold $B_A$ by
\[
 B_A :=
 \begin{cases}
  \{ * \}     & a=0, \\
  \I_j (\R^n ) & a>0.
 \end{cases}
\]
When $a=0$, we define $\hat{B}_A$ to be $C^o_{\sharp A}(\R^n )$ modulo scaling and translation;
\[
 \hat{B}_A := C^o_b (\R^n ) / \R^n \rtimes \R_{>0} ,
\]
where the action of $\R^n \rtimes \R_{>0}$ is defined by
\[
 (y_p )_{p \in A} \longmapsto (\alpha (y_p -\beta ))_{p \in A},\quad \forall \alpha >0,\ \forall \beta \in \R^n .
\]
The map $\rho_A : \hat{B}_A \to B_A =\{ *\}$ is defined as the canonical one.

When $a>0$, $\hat{B}_A$ is $\I_j (\R^n ) \times C^o_{a,b}$ modulo scaling and translation in the directions of $j$-planes;
\[
  \hat{B}_A := (\I_j (\R^n ) \times C^o_{a,b}) / \R^j \rtimes \R_{>0}, \\
\]
where the action of $\R^j \rtimes \R_{>0}$ is defined by
\[
 (\iota ;(x_p ,y_q )_{p,q \in A}) \longmapsto (\iota ;(\alpha (x_p -\beta ), \alpha (y_q -\iota (\beta )))_{p,q \in A}),
 \quad \forall \alpha >0,\ \forall \beta \in \R^j .
\]
The map $\rho_A : \hat{B}_A \to B_A =\I_j (\R^n )$ is defined as the natural projection.

Finally define $D_A : C^o_{\Gamma / \Gamma_A} \to B_A$ to be the canonical map if $a=0$, and
\[
 D_A (f;\dots ,x_p ,\dots ,z,\dots ,y_q ,\dots )_{p,q \not\in A}) := (df_z : T_z \R^j \longrightarrow T_{f(z)}\R^n )
 \in \I_j (\R^n )
\]
if $a>0$ (each tangent spaces are naturally identified with $\R^j$ and $\R^n$).

\begin{prop}[\cite{BottTaubes94, CCL02, Rossi_thesis, Watanabe07}]\label{Interior_strata}
The fibration $p_A : \Int \Sigma_A \to C^o_{\Gamma / \Gamma_A}$ is the pull-back of $\rho_A$ via $D_A$;
\[
 \xymatrix{
  \Int \Sigma_A \ar[r]^-{\hat{D}_A} \ar[d]_-{p_A} & \hat{B}_A \ar[d]^-{\rho_A} \\
  C^o_{\Gamma / \Gamma_A} \ar[r]^-{D_A} & B_A
 }
\]
In particular $\Int \Sigma_A \approx C^o_{\Gamma /\Gamma_A} \times \hat{B}_A$ if $A$ has no i-vertex.
\end{prop}

Notice that the differential form $\omega_{\Gamma / \Gamma_A} \in \Omega^*_{DR}(C^o_{\Gamma / \Gamma_A} )$ can be
defined similarly to \S \ref{subsection_form}, by using the direction maps $\varphi$ corresponding to the
edges of $\Gamma / \Gamma_A$.
Similarly, the maps $\varphi_e$ for any edges of $\Gamma_A$ are well defined on $\hat{B}_A$;
if $e=\overrightarrow{pq}$ is an $\eta$-edge,
\[
 \hat{\varphi}^{\eta}_e (\iota ; (x_r ,y_s )_{r,s \in A}) := \frac{x_q -x_p}{\abs{x_q -x_p}} \in S^{j-1},
\]
and if $e$ is a $\theta$-edge,
\[
 \hat{\varphi}^{\theta}_e (\iota ; (x_r ,y_s )_{r,s \in A}) := \frac{z_q -z_p}{\abs{z_q -z_p}} \in S^{n-1},
\]
where $z_p =x_p$ or $\iota (y_p )$ according to whether $p$ is internal or external.
Hence
\[
 \hat{\omega}_{\Gamma_A}:=\bigwedge_{e \in E(\Gamma_A )}\varphi^*_e vol \in \Omega^*_{DR} (\hat{B}_A)
\]
can be defined.
Then we have
\[
 \omega_{\Gamma}|_{\Int \Sigma_A} =\pm (p^*_A \omega_{\Gamma / \Gamma_A}) \wedge (\hat{D}^*_A \hat{\omega}_{\Gamma_A})
\]
and hence
\begin{equation}\label{p_A}
 (p_A )_* \omega_{\Gamma}|_{\Int \Sigma_A}=\pm \omega_{\Gamma / \Gamma_A} \wedge D^*_A (\rho_A )_* \hat{\omega}_{\Gamma_A}
\end{equation}
by the compatibility of fiber-integrations with pullbacks (for signs see \S\S \ref{subsec_orientation},
\ref{subsec_principal_cancel}).

Denote $\pi^{\partial_A}_{\Gamma} :=\pi_{\Gamma}|_{\Int \Sigma_A} : \Int \Sigma_A \to \emb{n}{j}$.
Notice that $\pi^{\partial_A}_{\Gamma}=\pi_{\Gamma / \Gamma_A} \circ p_A$.
We will often use the following criterion to show the vanishing of an integration along $\Int \Sigma_A$.

\begin{lem}[\cite{CCL02}]\label{vanish_1}
Let $a$ and $b$ be the numbers of i- and e-vertices in $A$ respectively.
Then the fiber integration $(\pi^{\partial_A}_{\Gamma})_* \omega_{\Gamma}|_{\Int \Sigma_A}$ vanishes unless
\begin{alignat*}{2}
 &\deg \hat{\omega}_{\Gamma_A} = nb-(n+1) &\quad &\text{if } a=0, \\
 &0 \le \deg \hat{\omega}_{\Gamma_A} - (ja+nb-(j+1)) \le nj &\quad &\text{if } a>0.
\end{alignat*}
\end{lem}

\begin{proof}
The form $(\rho_A )_* \hat{\omega}_{\Gamma_A}$ vanishes unless
$0 \le \deg (\rho_A )_* \hat{\omega}_{\Gamma_A} \le \dim B_A$.
We have
\[
 \deg (\rho_A )_* \hat{\omega}_{\Gamma_A} = \deg \hat{\omega}_{\Gamma_A}
 - \dim (\text{fiber of }\rho_A : \hat{B}_A \to B_A )
\]
and, by definition of the fibration $\rho_A : \hat{B}_A \to B_A$,
\[
 (\dim B_A ,\, \dim (\text{fiber of }\rho_A )) =
 \begin{cases}
  (0,\, nb -(n+1))      & \text{if } a=0, \\
  (nj,\, ja +nb -(j+1)) & \text{if } a>0.
 \end{cases}
\]
Hence $(\rho_A )_* \hat{\omega}_{\Gamma_A}$ vanishes unless the conditions of the Lemma are satisfied.
Then the formulas $(\pi^{\partial_A}_{\Gamma})_* = \pm (\pi_{\Gamma / \Gamma_A})_* \circ (p_A )_*$ and (\ref{p_A})
complete the proof.
\end{proof}

\subsection{Explicit description of strata at infinity}\label{subsection_strata_infinity}

Let $\Gamma$ be a graph and $A$ a non-empty subset of $V(\Gamma )$.
Below we describe $\Int \Sigma^{\infty}_A$ following \cite{BottTaubes94, CCL02, Rossi_thesis, Watanabe07}.

\begin{defn}
Define the {\it complementary graph} $\Gamma^c_A$ by letting
\[
 V(\Gamma^c_A ) := V(\Gamma ) \setminus A, \quad
 E(\Gamma^c_A ) := \{ \overrightarrow{pq} \in E(\Gamma ) \, | \, p,q \not\in A \}
\]
($\Gamma^c_A :=\emptyset$ if $A=V(\Gamma )$).
A {\it graph with infinity} is a graph with a specified vertex $v^{\infty}$ (called {\it vertex at infinity}),
which is not regarded as being internal nor external.

Define $\Gamma^{\infty}_A$, a graph with infinity, by `shrinking $\Gamma^c_A$ to a point.'
Namely $\Gamma^{\infty}_A$ is defined similarly as a quotient graph $\Gamma / \Gamma^c_A$ (see Definition
\ref{subgraph}), but its vertex $v_A$ is replaced by $v^{\infty}$, a vertex at infinity (see Figure
\ref{fig_subgraph_infty} for an example).
By definition $\Gamma^{\infty}_A :=\Gamma$ if $A=V(\Gamma )$.
\end{defn}
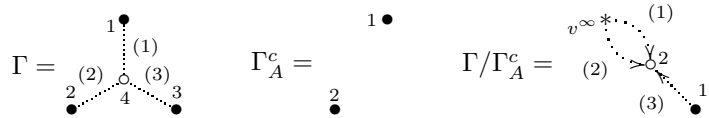
\begin{figure}[htb]
\[
 \begin{xy}
 (0,0)*{\circ}="O",(0,8)*{\bullet}="A",(-6.928,-4)*{\bullet}="B",(6.928,-4)*{\bullet}="C",
 {\ar@{.}"A";"O"_<{1}^(.45){(1)}},{\ar@{.}"B";"O"^<{2}^(.6){(2)}_>{4}},{\ar@{.}"C";"O"_<{3}_(.6){(3)}},(-12,2)*{\Gamma =},
 (35,8)*{\bullet},(33,8)*{\sb 1},(28.072,-4)*{\bullet},(28.072,-2)*{\sb 2},(21,2)*{\Gamma^c_A =},
 (70,2)*{\circ}="o",(64,8)*{*}="a",(76,-4)*{\bullet}="b",
 {\ar@{.>}@/^/"a";"o"^(.5){(1)}^>{2}}, {\ar@{.>}@/_/"a";"o"_<{v^{\infty}}_(.5){(2)}},{\ar@{.>}"b";"o"_<{1}^(.5){(3)}},
 (51,2)*{\Gamma /\Gamma^c_A =},
 \end{xy}
\]
\caption{Examples of $\Gamma^c_A$ and $\Gamma^{\infty}_A$ for $A= \{ 3,4\}$ ($n$ even)}\label{fig_subgraph_infty}
\end{figure}

$\Int \Sigma^{\infty}_A$ fibers over $C^o_{\Gamma^c_A}$;
\[
 p^{\infty}_A : \Int \Sigma^{\infty}_A \longrightarrow C^o_{\Gamma^c_A}
\]
which maps the limit point (\ref{infty_limit}) in \S \ref{subsection_boundary_strata} to $(x_p ,y_q )_{p,q \not\in A}$.
In other words $p^{\infty}_A$ forgets the points escaping to infinity.

As in \S \ref{subsection_strata_1}, we will define a space $\hat{B}^{\infty}_A$ which describes infinitesimal
configurations around infinity.
The space $\hat{B}^{\infty}_A$ is a subquotient of $C^o_{a,b}$ ($a,b$ are the numbers of i- and e-vertices in $A$)
modulo scaling;
\[
 \hat{B}^{\infty}_A :=
 \{ (x_p ,y_q )_{p,q\in A} \in C^o_{a+b}(\R^n \setminus 0) \, |\, x_p \in \R^j \times \{ 0\}^{n-j}\}/\R_{>0} .
\]
The origin $0 \in \R^n$ corresponds to $v^{\infty}$, which we use to fix the coordinates.
So in this case translation is not needed (compare it with the definition of $\hat{B}_A$).

\begin{prop}[\cite{BottTaubes94, CCL02, Rossi_thesis, Watanabe07}]\label{face_at_infinity}
$\Int \Sigma^{\infty}_A$ is homeomorphic to $C^o_{\Gamma^c_A} \times \hat{B}^{\infty}_A$.
\end{prop}

The form $\omega_{\Gamma^{\infty}_A} \in \Omega^*_{DR}(\hat{B}^{\infty}_A)$ can be defined as in \S \ref{subsection_form}
since the direction maps $\hat{\varphi}$ are invariant under scaling.
The vertex $v^{\infty}$ corresponds to $0 \in \R^n$.
More precisely, for each $e \in E(\Gamma^{\infty}_A)$, define the maps $\varphi_e : \hat{B}^{\infty}_A \to S^{N-1}$
($N=j$ or $n$) by
\[
 \varphi_e (x_p ,y_q )_{p,q \in A} :=
 \begin{cases}
  (z_q -z_p ) / \abs{z_q -z_p} & e = \overrightarrow{pq},\ p,q \in A, \\
  -z_p / \abs{z_p}             & e = \overrightarrow{pv^{\infty}} ,\ p \in A,
 \end{cases}
\]
where $z_p$ denotes $x_p$ or $y_p$ according to whether $p$ is an i-vertex or not.
Then
\[
 \omega_{\Gamma^{\infty}_A} := \bigwedge_{e \in E(\Gamma^{\infty}_A )} \omega_e .
\]
Under the identification in Proposition \ref{face_at_infinity},
\[
 \omega_{\Gamma}|_{\Int \Sigma^{\infty}_A}=\pm pr^*_1 \omega_{\Gamma^c_A}\wedge pr^*_2 \omega_{\Gamma^{\infty}_A}.
\]
Define $\pi^{\partial^{\infty}_A}_{\Gamma} : \Int \Sigma^{\infty}_A \to \emb{n}{j}$ as the restriction of $\pi_{\Gamma}$
onto the stratum $\Int \Sigma^{\infty}_A$ of fibers.
Then we have the following (which should be compared with Lemma \ref{vanish_1}).

\begin{lem}[\cite{CCL02}]\label{vanish_infty}
Let $a$ and $b$ be the numbers of i- and e-vertices of $\Gamma^{\infty}_A$ respectively (other than $v^{\infty}$).
Then the integration $(\pi^{\partial^{\infty}_A}_{\Gamma})_* \omega_{\Gamma}|_{\Int \Sigma^{\infty}_A}$ vanishes unless
$\deg \omega_{\Gamma^{\infty}_A} = \dim \hat{B}^{\infty}_A$, or equivalently, unless
\[
 \deg \omega_{\Gamma^{\infty}_A} = ja+nb-1.
\]
\end{lem}

\subsection{Orientations of boundary strata}\label{subsec_orientation}

Let $\Gamma$ be a graph and $A$ a non-empty subset of $V(\Gamma )$.

\begin{defn}
The boundary face $\Sigma_A$ (or the subgraph $\Gamma_A$) is said to be {\it principal} if $A$ consists of exactly
two vertices.
Similarly the boundary face $\Sigma^{\infty}_A$ is said to be {\it principal} if $A$ consists of exactly
one vertex.

All the other boundary strata are said to be {\it hidden}.
\end{defn}

Here we study the induced orientations of the non-infinity type principal strata $C_A$ from that of
$\text{fib}(\pi_{\Gamma}) \approx C_{s,t}$ (see \S \ref{subsection_boundary_strata}).
We are not interested in the orientations of strata at infinity and hidden strata, since the integrations along these
strata will be proved to vanish (see below).

Let $s$ and $t$ be the numbers of i- and e-vertices of $\Gamma$ respectively.
The fiber $C_{s,t}$ is equipped with the natural orientation as the subspace of $(\R^j )^s \times (\R^n )^t$.

Let $A \subset V(\Gamma )$ be a subset with $\sharp A=2$ (thus $\Sigma_A$ is principal), and $a$ and $b$ (with $a+b=2$)
the numbers of i- and e-vertices in $A$.

\noindent
\underline{\bf Case 1}: $a=2$, $b=0$.

Let $A=\{ x_p , x_q \}$ ($p<q \le s$).
In this case $\hat{B}_A = \I_j (\R^n ) \times S^{j-1}$ and hence
$\Int \Sigma_A \approx C_{\Gamma / \Gamma_A} \times S^{j-1}$ (see Proposition \ref{Interior_strata}). 
A neighborhood of $\Int \Sigma_A$ in $C_{\Gamma}$ is identified with $[0,1) \times C_{\Gamma / \Gamma_A} \times S^{j-1}$
by the homeomorphism onto the image
\begin{equation}\label{eq_coord1}
\begin{split}
 &(\varepsilon ;f;x_1 , \dots ,x_{q-1}, x_{q+1}, \dots ,x_s ;y_{s+1},\dots ,y_{s+t};v) \\
 &\hskip100pt \longmapsto (f; x_1 ,\dots ,x_{q-1}, x_p +\varepsilon v ,x_{q+1}, \dots ;y_{s+1},\dots )
\end{split}
\end{equation}
where $f\in \emb{n}{j}$ and $v \in S^{j-1}$.
It is not hard to see that the local coordinate of $C_{\Gamma}$ given by (\ref{eq_coord1}) has the orientation sign
$(-1)^{j(q-1)+(j-1)nt}$.
Putting $\varepsilon =0$ in the local coordinate (\ref{eq_coord1}), we obtain the natural orientation of
$C_{\Gamma / \Gamma_A} \times S^{j-1}$.
Thus the induced orientation of $C_{\Gamma /\Gamma_A} \times S^{j-1}$ as a boundary face of $C_{\Gamma}$ has the sign
$(-1)^{j(q-1)+(j-1)nt+1}$.

\noindent
\underline{\bf Case 2}: $a=b=1$ or $a=0$, $b=2$.

Let $A=\{ x_p , y_q \}$ ($p \le s < q$ or $s<p<q$).
In this case $\hat{B}_A =\I_j (\R^n ) \times S^{n-1}$ and hence
$\Int \Sigma_A \approx C_{\Gamma / \Gamma_A} \times S^{n-1}$ (see Proposition \ref{Interior_strata}). 
A neighborhood of $\Int \Sigma_A$ in $C_{\Gamma}$ is identified with $[0,1) \times C_{\Gamma / \Gamma_A} \times S^{n-1}$
by the homeomorphism onto the image
\begin{equation}\label{eq_coord2}
\begin{split}
 &(\varepsilon ;f;x_1 , \dots ,x_s ;y_{s+1},\dots y_{q-1}, y_{q+1}, \dots ,y_{s+t};w) \\
 &\hskip100pt \longmapsto (f; x_1 \dots ;y_{s+1},\dots ,y_{q-1}, z_p +\varepsilon w ,y_{q+1},\dots )
\end{split}
\end{equation}
where $w \in S^{n-1}$ and $z_p =f(x_p )$ or $y_p$ according to whether $p$ is internal or not.
The local coordinate of $C_{\Gamma}$ given by (\ref{eq_coord2}) has the orientation sign
$(-1)^{n(s+q-1)+js}$.
Putting $\varepsilon =0$ in the local coordinate (\ref{eq_coord2}), we obtain the natural orientation of
$C_{\Gamma / \Gamma_A} \times S^{n-1}$.
Thus the induced orientation of $C_{\Gamma /\Gamma_A}\times S^{n-1}$ as a boundary face of $C_{\Gamma}$ has the sign
$(-1)^{n(s+q-1)+js+1}$.

\subsection{Principal faces}\label{subsec_principal_cancel}
In this subsection we compute the fiber integration along the principal faces.

\begin{thm}[\cite{CCL02, Rossi_thesis, Watanabe07}]\label{principal}
The integration of $\omega_{\Gamma}$ along the principal face $\Sigma_A$ of non-infinity type (thus $\sharp A =2$)
vanishes unless the two vertices are joined by an edge in $\Gamma$.
\end{thm}

\begin{proof}
Let $a$ and $b$ be the numbers of i- and e-vertices in $A$ respectively ($a+b=2$ since $A$ is principal).
If two vertices in $A$ are not connected by an edge, then we have
$\hat{\omega}_{\Gamma_A}=1 \in \Omega^0_{DR}(\hat{B}_A)$.
Thus, if $A$ has no i-vertex ($a=0$, $b=2$), then the first equality of Lemma \ref{vanish_1} does not hold;
\[
 \deg \hat{\omega}_{\Gamma_A} = 0 \ne 0\cdot j + 2n - (n+1) =n-1.
\]
If $A$ has an i-vertex ($(a,b)=(1,1)$ or $(2,0)$), then the second inequality of Lemma \ref{vanish_1} does not
hold since
\[
 \deg \hat{\omega}_{\Gamma_A} - (ja+nb-(j+1)) =
 \begin{cases}
 -(j-1) & \text{if }(a,b)=(2,0) \\
 -(n-1) & \text{if }(a,b)=(1,1).
 \end{cases}
\]
\end{proof}

\begin{thm}[\cite{CCL02, Rossi_thesis, Watanabe07}]\label{infty_principal}
Let $\Gamma$ be an admissible graph.
Then the integration of $\omega_{\Gamma}$ along the principal face $\Sigma^{\infty}_A$ (thus $\sharp A =1$) always
vanishes.
\end{thm}

\begin{proof}
Let $p$ be the only vertex in $A$.
The graph $\Gamma^{\infty}_A =\Gamma / \Gamma^c_A$ has two vertices; one is $p$ and the other is $v^{\infty}$.
$\sharp E(\Gamma^{\infty}_A )$ is equal to the valency of $p$ in $\Gamma$.

By Lemma \ref{vanish_infty}, $(\pi^{\partial^{\infty}_A}_{\Gamma})_* \omega_{\Gamma}|_{\Int \Sigma^{\infty}_A}$
vanishes unless
\[
 \deg \omega_{\Gamma^{\infty}_A}=
 \begin{cases}
  j-1 & \text{if } p \text{ is internal}, \\
  n-1 & \text{if } p \text{ is external},
 \end{cases}
\]
hence $(\pi^{\partial^{\infty}_A}_{\Gamma})_* \omega_{\Gamma}|_{\Int \Sigma^{\infty}_A}$ does not vanish only if,
in $\Gamma$,
\begin{itemize}
\item $p$ is a uni-valent i-vertex with exactly one adjacent $\eta$-edge, or
\item $p$ is a uni-valent e-vertex with exactly one adjacent $\theta$-edge.
\end{itemize}
But neither case occurs since $\Gamma$ is admissible (see Definition \ref{def_admissible}).
\end{proof}

\begin{thm}[\cite{CCL02, Rossi_thesis, Watanabe07}]\label{3.1_first_half}
The sum of integrations of $\omega_{\Gamma}$ along all the principal faces $\Sigma_A$ of non-infinity type is equal to
$I(\delta \Gamma )$.
\end{thm}

\begin{proof}
By the above Theorems \ref{principal} and \ref{infty_principal}, we only need to consider the principal faces
$\Sigma_A$ such that the two vertices of $\Gamma_A$ are joined by an edge $e$.
We will show
\[
 (\pi^{\partial_A}_{\Gamma})_* \omega_{\Gamma}|_{\Int \Sigma_A} = (-1)^{\tau (e)}(\pi_{\Gamma /e})_* \omega_{\Gamma /e}
\]
for any principal strata $\Sigma_A$, where $e\in E(\Gamma )$ is the only edge of the subgraph $\Gamma_A$ and $\tau (e)$
is the sign given in Proposition \ref{def_signs}.
Then we will obtain
\[
 \sum_{\genfrac{}{}{0pt}{}{A \subset V(\Gamma )}{\sharp A=2}}
 (\pi^{\partial_A}_{\Gamma})_* \omega_{\Gamma}|_{\Int \Sigma_A}
 = \sum_{e \in E(\Gamma ) \setminus \{ \text{loops}\}} (-1)^{\tau (e)}(\pi_{\Gamma /e})_* \omega_{\Gamma /e}
 =I(\delta \Gamma ).
\]
First we consider the case when $n$ and $j$ are odd.
We divide the proof into four cases.

\noindent
\underline{\bf Case (a)}.
$A=\{ p,q \}$ consists of two e-vertices (so we can assume $s < p < q$ and the edge $e = \overrightarrow{pq}$
is a $\theta$-edge).

In this case $\hat{B}_A = S^{n-1}$ and $B_A =\{ * \}$, so $\Int \Sigma_A \approx C^o_{\Gamma /e} \times S^{n-1}$.
The induced orientation of $C^o_{\Gamma /e} \times S^{n-1}$ from $C_{\Gamma}$ is $(-1)^{n(s+q-1)+js+1}$ by Case 2 in
\S \ref{subsec_orientation}, and is equal to $(-1)^q$ since $n$ and $j$ are odd.
This sign is $(-1)^{\tau (e)}$ (see Proposition \ref{def_signs}).

Under the identification $\Int \Sigma_A \approx C^o_{\Gamma /e} \times S^{n-1}$, the map
$\varphi_e : C^o_{\Gamma} \to S^{n-1}$ restricts to the projection $pr_2 : C^o_{\Gamma /e} \times S^{n-1} \to S^{n-1}$.
Hence via the diffeomorphism $\Int \Sigma_A \approx S^{n-1} \times C^o_{\Gamma /e}$, the form
$\omega_{\Gamma_A}|_{\Int \Sigma_A}$ corresponds to
\[
  (-1)^{\tau (e)} p^*_A \omega_{\Gamma /e} \wedge pr^*_2 vol_{S^{n-1}}
 \in \Omega^*_{DR}(S^{n-1} \times C^o_{\Gamma /e})
\]
and hence we have
\begin{align*}
 (\pi^{\partial_A}_{\Gamma})_* \omega_{\Gamma}|_{\Int \Sigma_A}
 &= (-1)^{\tau (e)} (\pi_{\Gamma / e})_* \circ (p_A )_* (p^*_A \omega_{\Gamma / e} \wedge pr^*_2 vol_{S^{n-1}}) \\
 &= (-1)^{\tau (e)} (\pi_{\Gamma / e})_* \left( \int_{S^{n-1}} vol_{S^{n-1}} \right) \omega_{\Gamma / e} \\
 &= (-1)^{\tau (e)} (\pi_{\Gamma / e})_* \omega_{\Gamma / e}.
\end{align*}

\noindent
\underline{\bf Case (b)}.
$A= \{ p,q\}$ contains both an e- and an i-vertex (thus we can assume $p \le s <q$, and
$e=\overrightarrow{pq}$ is a $\theta$-edge).

In this case $\Int \Sigma_A \approx C^o_{\Gamma / e} \times S^{n-1}$.
Similarly as in Case (a) above, the induced orientation of $C^o_{\Gamma / e} \times S^{n-1}$ from $C_{\Gamma}$
is $(-1)^q =(-1)^{\tau (e)}$ and $\varphi_e : C^o_{\Gamma} \to S^{n-1}$ restricts to the projection
$pr_2 : C^o_{\Gamma /e} \times S^{n-1} \to S^{n-1}$.
Thus, as in the Case (a),
\[
 (\pi_{\Gamma})_* \omega_{\Gamma}|_{\Int \Sigma_A}  = (-1)^{\tau (e)} (\pi_{\Gamma /e})_* \omega_{\Gamma /e}.
\]

\noindent
\underline{\bf Case (c)}.
Both two points $p,q$ of $A$ are internal and $e=\overrightarrow{pq}$ is an $\eta$-edge.

In this case $\Int \Sigma_A \approx C^o_{\Gamma /e} \times S^{j-1}$.
Proof is the same as the above cases, since by Case 1 in \S \ref{subsec_orientation}, the induced orientation of
$C^o_{\Gamma / e} \times S^{j-1}$ from $C_{\Gamma}$ is $(-1)^q =(-1)^{\tau (e)}$.

\noindent
\underline{\bf Case (d)}.
$A=\{ p,q\}$ consists of two i-vertices (thus we can assume $p<q\le s$) and the edge $e=\overrightarrow{pq}$ is a
$\theta$-edge.

In this case $\Int \Sigma_A \approx C^o_{\tilde{\Gamma}} \times S^{j-1}$, where $\tilde{\Gamma}$ is $\Gamma /e$
with its small loop $e$ removed.
The right hand side is nothing but the space $C^o_{\Gamma /e}$ ($\Gamma /e$ is a graph with small loop; see \S
\ref{subsection_integral}), up to the orientation.
Again by Case 1 in \S \ref{subsec_orientation}, the induced orientation of $\Int \Sigma_A \subset C_{\Gamma}$ has the
sign $(-1)^{\tau (e)}$.
Hence
\[
 (\pi_{\Gamma})_* \omega_{\Gamma}|_{\Int \Sigma_A}  = (-1)^{\tau (e)} (\pi_{\Gamma /e})_* \omega_{\Gamma /e}
\]
as desired.

The proof of the case when $n$ and $j$ are even is similar.
In Case (a), the induced orientation of $\Sigma_A \subset C_{\Gamma}$ has the sign $-1$.
To integrate the form $\theta_e$ first, we must put $\theta_e$ at the top of $\omega_{\Gamma}$.
Such a re-ordering yields the sign $(-1)^{i-1}$ if $e$ is the $i$-th edge.
Hence
\[
 (\pi^{\partial_A}_{\Gamma})_* \omega_{\Gamma}|_{\Int \Sigma_A} = (-1)^i (\pi_{\Gamma / e})_* \omega_{\Gamma / e}
 = (-1)^{\tau (e)}(\pi_{\Gamma / e})_* \omega_{\Gamma / e}.
\]
The remaining three cases are proved in similar ways.
In Case (d), we have to put the $S^{j-1}$-factor at the end of the $(S^{j-1})^u$-part, so the sign $(-1)^u$ appears.

When $n$ is even, we have to consider one more case;

\noindent
\underline{\bf Case (e)}.
$A=\{ p,q\}$ consists of two internal vertices (thus we can assume $p<q\le s$) which are joined by an $\eta$-edge
$\overrightarrow{pq}_{\eta}$ and a $\theta$-edge $\overrightarrow{pq}_{\theta}$.

In this case, the induced orientation of $\Sigma_A \subset C_{\Gamma}$ has the sign $-1$ as in Case (a).
But we need to put $S^{j-1}$ at the end of $S^{j-1}$-factor, and not to move the forms $\theta$ and $\eta$.
Hence
\[
 (\pi^{\partial_A}_{\Gamma})_* \omega_{\Gamma}|_{\Int \Sigma_A} = (-1)^{u+1}(\pi_{\Gamma / e})_* \omega_{\Gamma / e}
 = (-1)^{\tau (e)}(\pi_{\Gamma / e})_* \omega_{\Gamma / e}.
\]
\end{proof}

By the above Theorem \ref{3.1_first_half}, the proof of Theorem \ref{cochain} is reduced to showing that
hidden faces do not contribute to the fiber integration.
The following, whose proof will be given in \S \ref{subsection_hidden_first} and \S \ref{subsection_hidden_infty},
will complete the proof of Theorem \ref{cochain}.

\begin{thm}\label{3.1_second_half}
Let $\Gamma$ be an admissible graph.
Then all the integrations of $\omega_{\Gamma}$ along hidden boundary faces of the fiber of
$\pi_{\Gamma} : C_{\Gamma} \to \emb{n}{j}$ vanish if
(1) $n-j \ge 2$ is even and $\Gamma$ is a tree, or
(2) both $n>j\ge 3$ are odd and $\Gamma$ has at most one loop component.
\end{thm}

\subsection{Hidden faces; the non-infinity type}\label{subsection_hidden_first}

Here we show that all the hidden faces $\Sigma_A$ of non-infinity type, $\sharp A \ge 3$, do not contribute
to the fiber integration.

\begin{lem}[\cite{Rossi_thesis, Watanabe07}]\label{disconnected}
Suppose $A \subset V(\Gamma )$ is such that the subgraph $\Gamma_A$ is not connected.
Then $(\pi^{\partial_A}_{\Gamma})_* \omega_{\Gamma}|_{\Int \Sigma_A}$ vanishes.
\end{lem}

\begin{proof}
If $\Sigma_A$ is principal (so $\sharp A =2$), then the claim of this Lemma is exactly that of Theorem \ref{principal}.
So we assume $\sharp A \ge 3$.

Suppose $\Gamma_A = \Gamma_{A_1} \sqcup \Gamma_{A_2}$ for non-empty subsets $A_1 ,A_2 \subset A$ (thus we can assume
$\sharp A_1 \ge 2$).
Let $a_i$ and $b_i$ be the numbers of i- and e-vertices of $A_i$, $i=1,2$.

We define the space $\tilde{B}_A$, which contains $\hat{B}_A$ as an open subset, by
\[
 \tilde{B}_A := (C^0_{a_1 ,b_1} \times C^o_{a_2 ,b_2}) / \sim
\]
here $\sim$ is defined by using the translation and the scaling.
In other words, `a point in $A_1$ may collide with a point in $A_2$.'
Since there are no edges joining a point in $A_1$ to that of $A_2$, the direction maps $\varphi$ corresponding to
the edges of $\Gamma_A$ are well-defined on $\tilde{B}_A$, and so is the associated differential form (we denote it by
$\tilde{\omega}_{\Gamma_A}$).
The restriction of $\tilde{\omega}_{\Gamma_A}$ onto $\hat{B}_A$ is $\hat{\omega}_{\Gamma_A}$.

Consider a free action of $\R^N$ on $\tilde{B}_A$ defined by the translations of points in $A_1$ (points in $A_2$ are
fixed).
Here $N=j$ or $n$ according to whether $A_1$ contains an i-vertex or not.
Let $p:\tilde{B}_A \to \tilde{B}_A / \R^N$ be the quotient map.
Since the direction maps $\tilde{\varphi}:\tilde{B}_A \to S^{j-1}$ or $S^{n-1}$ factor through $p$, there exists a
form $\omega'_{\Gamma_A} \in \Omega^*_{DR}(\tilde{B}_A /\R^N )$ such that
$p^* \omega'_{\Gamma_A} = \tilde{\omega}_{\Gamma_A}$.
We have a map of fibrations
\[
 \xymatrix{
 \hat{B}_A \ar@{^{(}->}[r] \ar[rd]_-{\rho_A} & \tilde{B}_A \ar[d]^-{\tilde{\rho}_A} \ar@{->>}[r]^-p
  & \tilde{B}_A / \R^N \ar[ld]^-{\rho'_A} \\
  & B_A &
 }
\]
(for definition of $\rho_A$ see Proposition \ref{Interior_strata})
and it holds that
\[
 (\rho_A )_* \omega_{\Gamma_A} =(\tilde{\rho}_A )_* \tilde{\omega}_{\Gamma_A} = (\rho'_A )_* \omega'_{\Gamma_A}.
\]
This implies $(\rho_A )_* \omega_{\Gamma_A}=0$, since the fiber of $\rho'_A$ is of strictly less dimension than those of
$\rho_A$ and $\tilde{\rho}_A$.
This together with the formula (\ref{p_A}) completes the proof.
\end{proof}

Thanks to Lemma \ref{disconnected}, below we can assume that $\Gamma_A$ is connected.

In Theorem \ref{cochain} we assumed $g \le 1$, that is, our graph has at most one loop component (other than small
loops), and so does its connected subgraph $\Gamma_A$.

\begin{prop}\label{vanish_tree}
Suppose $n-j$ is even.
If $\sharp A \ge 3$ and $\Gamma_A$ is a tree, then
$(\pi^{\partial_A}_{\Gamma})_* \omega_{\Gamma}|_{\Int \Sigma_A}$ vanishes.
\end{prop}

\begin{proof}
Since $\Gamma_A$ is a tree, there are at least two uni-valent vertices in $\Gamma_A$.
All the possibilities of uni-valent vertices are listed in Figure \ref{types_ends}.
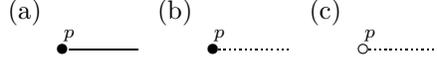
\begin{figure}[htb]
\[
 \begin{xy}
 (-5,5)*{\text{(a)}},
 (0,0)*{\bullet}="A",
 {\ar@{-}"A";(10,0)^<{p}},
 (15,5)*{\text{(b)}},
 (20,0)*{\bullet}="B",
 {\ar@{.}"B";(30,0)^<{p}},
 (35,5)*{\text{(c)}},
 (40,0)*{\circ}="C",
 {\ar@{.}"C";(50,0)^<{p}}
 \end{xy}
\]
\caption{Uni-valent vertices}\label{types_ends}
\end{figure}
We will prove the vanishing of integration along $\Sigma_A$ in the successive Lemmas;
types (a) and (c) in Lemma \ref{hidden_1} and type (b) in Lemma \ref{hidden_2}.
The assumption that $n-j$ is even will be used in the proof of Lemma \ref{hidden_2}.
\end{proof}

\begin{lem}[\cite{CCL02}]\label{hidden_1}
If $\sharp A \ge 3$ and $A$ has a uni-valent vertex $p$ of types (a) or (c) in Figure \ref{types_ends}, then the
integration $(\pi^{\partial_A}_{\Gamma})_* \omega_{\Gamma}|_{\Int \Sigma_A}$ vanishes.
\end{lem}

\begin{proof}
Let $q \in A$ be the vertex joined to $p$ in $\Gamma_A$, which must be internal in the case (a), while in the case (c)
it may be both internal or external.
There is a fiberwise free action $\R_{>0}$ on $\hat{B}_A$ defined on each fiber by
\[
 (\dots ,z_p , \dots )\longmapsto (\dots ,az_p +(1-a)z_q ,\dots ),\quad a\in \R_{>0}.
\]
Then $\hat{\omega}_{\Gamma_A} \in \Omega^*_{DR}(\hat{B}_A)$ is basic with respect to the quotient
$\hat{B}_A \to \hat{B}_A /\R_{>0}$.
Since $\sharp A \ge 3$, the fiber of $\hat{B}_A / \R_{>0} \to B_A$ is of strictly less dimension than that of
$\rho_A$.
Hence the similar argument as in Lemma \ref{disconnected} completes the proof.
\end{proof}

\begin{lem}\label{hidden_2}
Suppose $n-j$ is even.
If $\sharp A \ge 3$ and $\Gamma_A$ is a tree all of whose uni-valent vertices are of type (b) in Figure \ref{types_ends},
then $(\pi_{\Gamma})_* \omega_{\Gamma}|_{\Int \Sigma_A}$ vanishes.
\end{lem}

\begin{proof}
The vertex $q$ of $A$ which is joined to $p$ may be both internal or external.
Since valency of $q$ is greater than one, there exist vertices $r_1 ,\dots ,r_a$ ($r_i \ne p$, $a \ge 1$) which are
joined to $q$.

Suppose one of them, say $r_1$, is uni-valent.
By our assumption $r_1$ is also of type (b).
Consider a fiberwise involution $\chi_1 : \hat{B}_A \to \hat{B}_A$ defined by
\[
 \chi_1 (\iota ; \dots ,x_p , \dots ,x_{r_1},\dots ) := (\iota ; \dots ,x_{r_1}, \dots ,x_p ,\dots )
\]
(other coordinates are not changed).
This involution changes the orientation of the fiber by $(-1)^j$, while
$\chi_1^* \hat{\omega}_{\Gamma_A}=(-1)^{n-1}\hat{\omega}_{\Gamma_A}$ since $\chi^*_1 \theta_{pq}=\theta_{r_1 q}$,
$\chi^*_1 \theta_{r_1 q}=\theta_{pq}$ and $\theta$'s are of degree $n-1$.
Thus
\[
 (\rho_A )_* \hat{\omega}_{\Gamma_A} =(-1)^{j+n-1}(\rho_A )_* \hat{\omega}_{\Gamma_A}
 =-(\rho_A )_* \hat{\omega}_{\Gamma_A}
\]
since $n-j$ is even, and it must vanish.

Thus we may assume all the $r_i$'s are at least bi-valent.
That is, we can assume that all the uni-valent vertices $p$ (of type (b)) and adjacent $q$ are such that no other
uni-valent vertex is joined to $q$.
Then we can find at least two pairs $(p,q)$ of vertices such that $p$ is uni-valent vertex joined to exactly one
bi-valent vertex $q$; since otherwise $\Gamma_A$ cannot be a tree.
Such pairs, say $(p,q)$ and $(p' ,q' )$, are of types (b-1) or (b-2) or (b-3) in Figure \ref{fig_b1}, where an asterisk
can be both internal or external.
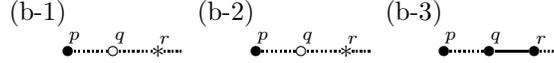
\begin{figure}[htb]
\[
 \begin{xy}
 (-4,5)*{\text{(b-1)}},
 (0,0)*{\bullet}="A",(6,0)*{\circ}="B",(12,0)*{*}="C",
 {\ar@{.}"A";"B"^<{p}},{\ar@{.}"B";"C"^<{q}},{\ar@{.}"C";(15,0)^<{r}},
 (21,5)*{\text{(b-2)}},
 (25,0)*{\bullet}="D",(31,0)*{\circ}="E",(37,0)*{*}="F",
 {\ar@{.}"D";"E"^<{p}},{\ar@{.}"E";"F"^<{q}},{\ar@{.}"F";(40,0)^<{r}},
 (46,5)*{\text{(b-3)}},
 (50,0)*{\bullet}="G",(56,0)*{\bullet}="H",(62,0)*{\bullet}="I",
 {\ar@{.}"G";"H"^<{p}},{\ar@{-}"H";"I"^<{q}},{\ar@{.}"I";(65,0)^<{r}}
 \end{xy}
\]
\caption{Types (b-1), (b-2), (b-3)}\label{fig_b1}
\end{figure}

If $(p,q)$ is of type (b-1), define a fiberwise involution $\chi_2 :\hat{B}_A \to \hat{B}_A$ by
\[
 \chi_2 (\iota ; \dots ,y_q , \dots ) := (\iota ; \dots ,\iota(x_p )+ z_r -y_q ,\dots ),
\]
where $z_r$ is $\iota(x_r )$ or $y_r$ according to whether $r$ is internal or not.
This involution has orientation sign $(-1)^n$, while $\chi^*_2 \hat{\omega}_{\Gamma_A}=(-1)^{n-1}\hat{\omega}_{\Gamma_A}$
similarly to $\chi_1$.
Thus
\[
 (\rho_A )_* \hat{\omega}_{\Gamma_A} =(-1)^{n+n-1}(\rho_A )_* \hat{\omega}_{\Gamma_A}
 =-(\rho_A )_* \hat{\omega}_{\Gamma_A}
\]
and it must vanish.

So finally we can assume that both pairs $(p,q)$ and $(p' ,q' )$ are of types (b-2) or (b-3).
The proof of this case appeared in \cite{Rossi_thesis, Watanabe07}; consider a fiberwise involution
$\chi_3 :\hat{B}_A \to \hat{B}_A$ defined by
\[
 \chi_3 (\iota ; \dots ,x_p , \dots ,x_{p'},\dots )
 := (\iota ; \dots , x_q -x_{q'} + x_{p'}, \dots ,x_q - x_{q'} + x_p , \dots ).
\]
The orientation sign of $\chi_3$ is $(-1)^j$, and it satisfies
$\chi^*_3 \hat{\omega}_{\Gamma_A}=(-1)^{n-1}\hat{\omega}_{\Gamma_A}$ since $\chi^*_3 \theta_{pq}=\theta_{p' q'}$,
$\chi^*_3 \theta_{p'q'}=\theta_{pq}$.
Since $n-j$ is even,
\[
 (\rho_A )_* \hat{\omega}_{\Gamma_A} =(-1)^{n+j-1}(\rho_A )_* \hat{\omega}_{\Gamma_A}
 =-(\rho_A )_* \hat{\omega}_{\Gamma_A}
\]
as in the case of $\chi_1$, and it must vanish.
\end{proof}

When both $n>j\ge 3$ are odd, we have to consider the case that $\Gamma_A$ has one loop component to prove Theorem
\ref{cochain}.

\begin{prop}[\cite{Rossi_thesis, Watanabe07}]
Suppose $n>j\ge 3$ are odd.
If $\sharp A \ge 3$ and $\Gamma_A$ is connected with one loop component, then
$(\pi_{\Gamma})_* \omega_{\Gamma}|_{\Int \Sigma_A}$ vanishes.
\end{prop}

\begin{proof}
Let $s$ and $t$ be the numbers of i- and e-vertices in $A$ respectively.
Define an involution $F$ of $\hat{B}_A$ by
\begin{multline*}
 F(\iota ; x_1 ,\dots ,x_s ; y_{s+1} ,\dots ,y_{s+t}) := \\
 (\iota ; x_1 ,2x_1 -x_2 ,\dots ,2x_1 -x_s ;2\iota (x_1 )-y_{s+1},\dots ,2\iota (x_1 )-y_{s+t}),
\end{multline*}
whose orientation sign is $(-1)^{j(s-1)+nt}=(-1)^{s+t-1}$.

Let $\alpha$ and $\beta$ be the numbers of $\eta$- and $\theta$-edges of $\Gamma_A$ respectively.
Since $\hat{\varphi}_e \circ F =\iota_{S^{N-1}}\circ \hat{\varphi}_e$ for any edge $e$ ($N=j$ or $n$ according to whether
$e$ is an $\eta$-edge or a $\theta$-edge), we have
\[
 F^* \omega_{\Gamma_A} =(-1)^{j\alpha + n\beta} \omega_{\Gamma_A} =(-1)^{\alpha + \beta}\omega_{\Gamma_A}.
\]
But by our assumption, $\Gamma_A$ has exactly one loop component, so $\alpha+ \beta$ is equal to the number of
vertices of $\Gamma_A$, that is, $s+t$.
Thus
\[
 (\rho_A )_* \omega_{\Gamma_A} = (-1)^{s+t-1}(-1)^{s+t}(\rho_A )_* \omega_{\Gamma_A} =-(\rho_A )_* \omega_{\Gamma_A}
\]
and hence $(\rho_A )_* \omega_{\Gamma_A} = 0$.
The formula (\ref{p_A}) in \S \ref{subsection_strata_1} completes the proof.
\end{proof}

\subsection{Hidden faces; strata at infinity}\label{subsection_hidden_infty}

In this subsection we will prove that the hidden strata $\Sigma^{\infty}_A$ at infinity (thus $\sharp A \ge 2$)
do not contribute to integrals.

\begin{lem}[\cite{CCL02, Rossi_thesis, Watanabe07}]\label{lem_tree_infty}
If $\sharp A \ge 2$, then the integration along $\Sigma^{\infty}_A$ vanishes.
\end{lem}

\begin{proof}
First we will show that $\omega_{\Gamma^{\infty}_A}$ cannot satisfy the equation in Lemma \ref{vanish_infty} when
$A \subset V(\Gamma )$ is a proper subset.

If $\Gamma$ is an admissible graph, then each e-vertex of $\Gamma^{\infty}_A$ ($=\Gamma / \Gamma^c_A$) is of at least
tri-valent, and each i-vertex is an endpoint of some $\theta$-edge.
This implies
\begin{equation}\label{eq_deg_positive}
 2\sharp \{ \theta\text{-edges of }\Gamma^{\infty}_A \} \ge a+3b,
\end{equation}
where $a$ and $b$ are numbers of i- and e-vertices of $\Gamma_A$.

Since $\Gamma^{\infty}_A \setminus \{ v^{\infty}\}$ is a one-dimensional open object with at most one loop component and
with at least one open edge,
\begin{equation}\label{number_edge}
 \sharp \{ \eta \text{-edges of } \Gamma^{\infty}_A \} + \sharp \{ \theta \text{-edges of } \Gamma^{\infty}_A \} \ge a+b.
\end{equation}
By using estimations \eqref{eq_deg_positive}, \eqref{number_edge} and $n-j-2\ge 0$, we see that
$\omega_{\Gamma^{\infty}_A}$ does not satisfy the criterion of Lemma \ref{vanish_infty};
\begin{align*}
 \deg \omega_{\Gamma^{\infty}_A}
 & \ge (n-1)\sharp \{ \theta \text{-edges of } \Gamma^{\infty}_A \}
  +(j-1)(a+b-\sharp \{ \theta \text{-edges of } \Gamma^{\infty}_A \} ) \\
 & = (n-j)\sharp \{ \theta \text{-edges of } \Gamma^{\infty}_A \} + (j-1)(a+b) \\
 & \ge \frac{1}{2}(n-j)(a+3b) + (j-1)(a+b) \\
 & = (ja+nb-1) +\frac{1}{2}(n-j-2)(a+b)+1 \\
 & > ja+nb-1 \\
 & = \dim \hat{B}^{\infty}_A .
\end{align*}
Next consider the case $A=V(\Gamma )$ (thus $\Gamma^{\infty}_A =\Gamma$).
Define a fiberwise free action of $\R_{>0}$ on $\hat{B}^{\infty}_A$ by
\[
 (x_1 ,\dots ,y_{a+1},\dots ) \longmapsto
 (x_1 , \alpha x_2 +(1-\alpha )x_1 , \dots ,\alpha y_{a+1} + (1-\alpha )f_0 (x_1 ), \dots ),
\]
where $f_0 : \R^j \hookrightarrow \R^n$ is the standard inclusion given by $x \mapsto (x,0,\dots ,0)$.
This action is non-trivial since $\sharp A \ge 2$.
The differential form $\omega_{\Gamma^{\infty}_{V(\Gamma)}}$ is basic with respect to this action, hence similar argument
to Lemma \ref{disconnected} completes the proof.
\end{proof}

\begin{rem}
In the end of the proof we used that the long knots are standard near infinity, so $\omega_{\Gamma^{\infty}_{V(\Gamma )}}$
contains no information about the base space $\emb{n}{j}$.
\end{rem}

Thus we have shown that all of hidden and infinity contributions vanish, and completed the proof of
Theorem \ref{3.1_second_half}.

\subsection{Proof of Theorem \ref{thm_closed}}\label{subsec_proof_closed}
In the above proofs, we have used
\begin{itemize}
\item symmetry of the fiber, and
\item dimension counting.
\end{itemize}
Below we check that some of the arguments are valid even if $n-j \ge 3$ is odd (in particular $n=6k$, $j=4k-1$), and
can be used to prove Theorem \ref{thm_closed}.

\begin{lem}\label{lem_dH=0}
$d\calH =I(\delta H)$ modulo the contributions of hidden faces.
\end{lem}

\begin{proof}
Exactly similar to Theorem \ref{3.1_first_half}.
We used Theorems \ref{principal} and \ref{infty_principal} to prove Theorem \ref{3.1_first_half}, which are proved by
only dimension counting.
\end{proof}

\begin{lem}\label{lem_dH_hidden}
The hidden faces except for the face $\Sigma_A$ corresponding to $A=V(H_2 )$ do not contribute to the integral.
\end{lem}

\begin{proof}
The strata at infinity do not contribute; the proof of Lemma \ref{lem_tree_infty} does not use the parities of $n,j$.
By Lemma \ref{disconnected}, which uses only dimension counting, we have only to consider the faces $\Sigma_A$ for $A$
such that the subgraph $\Gamma_A$ is connected.

For $A \subset V(H_1 )$ with $\sharp A =3$, use Lemma \ref{hidden_1}.

For $A=V(H_1 )$, consider a fiberwise involution $F:\hat{B}_A \to \hat{B}_A$ defined by
\[
 F(\iota ; (x_1 ,\dots ,x_4 )) := (\iota ; (2x_2 -x_1 ,x_2,x_3 ,x_4)).
\]
The orientation sign of $F$ is $(-1)^j$, while $F^* \hat{\omega}_{H_1}=(-1)^n \hat{\omega}_{H_1}$ because
$F^* \theta_{12} = (-1)^n \theta_{12}$ and $F^*$ preserves $\eta_{23}$, $\theta_{34}$.
Since $n-j$ is odd,
\[
 (\rho_A )_* \hat{\omega}_{H_1} =(-1)^{n+j}(\rho_A )_* \hat{\omega}_{H_1} = -(\rho_A )_* \hat{\omega}_{H_1}
\]
and hence the integration along $\Sigma_{V(H_1)}$ must vanish.

For $A \subset V(H_2 )$ with $\sharp A=3$, then $A$ contains two i-vertices joined to the e-vertex labeled by $4$.
Then using the involution $\chi_2$ appeared in the proof of Lemma \ref{hidden_2}, vanishing for type (b-1),
we can complete the proof.
\end{proof}

Next consider the contribution of `anomalous face' $\Sigma_A$, $A=V(H_2 )$ (hence $\Gamma_A =H_2$).
Recall from \S \ref{subsection_strata_1} that the face $\Sigma_A$ is described by the pullback square in Proposition
\ref{Interior_strata} with $B_A =\I_j (\R^n )$ and $C_{\Gamma / \Gamma_A}=\emb{n}{j}\times \R^j$.
The contribution of $\Sigma_A$ is given by $p_* D^*_A (\rho_A )_* \hat{\omega}_{H_2}$, where
$p:\emb{n}{j} \times \R^j \to \emb{n}{j}$ is the first projection.

At present we cannot determine whether this contribution vanishes or not, so
\[
 dI(H) = \frac{1}{6}p_* D^*_A (\rho_A )_* \hat{\omega}_{\Gamma_A}.
\]
But we can define a correction term $c$ which kills this contribution as follows.

\begin{lem}\label{lem_hatomega_closed}
The form $(\rho_A )_* \hat{\omega}_{H_2} \in \Omega^{2n-2j-2}_{DR}(\I_j (\R^n ))$ is closed.
\end{lem}

\begin{proof}
By the generalized Stokes theorem,
\[
 d(\rho_A )_* \hat{\omega}_{H_2} = \pm (\rho^{\partial}_A )_* \hat{\omega}_{H_2},
\]
where $\rho^{\partial}_A$ is the restriction of $\rho_A$ onto the boundary of the fiber.
The principal faces correspond to the graph (up to labeling)
\[
 \begin{xy}
 (0,0)*{\bullet}="A", (10,0)*{\bullet}="B", (20,0)*{\bullet}="C",
 {\ar@{.}"A";"B"^<{1}},{\ar@{.}"B";"C"^>{3}},
 (10,2)*{\sb 2}
 \end{xy}
\]
obtained by contracting the edge $i4$ ($i=1,2,3$) of $H_2$.
Consider the involution
\[
 F : (x_1 ,x_2 ,x_3 ) \longmapsto (x_1 ,x_2 ,2x_2 -x_3)
\]
of the principal face.
The orientation sign of $F$ is $(-1)^j$, while $F^* \hat{\omega}_{H_2}=(-1)^n \hat{\omega}_{H_2}$
since $F^* \theta_{12} = \theta_{12}$ and $F^*_1 \theta_{23} = \theta_{32}=(-1)^n \theta_{23}$.
Since $j+n$ is odd, the integration along the principal faces vanish.

The hidden (but non-anomalous) contributions are proved to vanish in completely similar way as in
Lemma \ref{lem_dH_hidden}.
Since there is no anomalous face of the fiber of $\hat{B}_A \to B_A$, we have proved that
$d(\rho_A )_* \hat{\omega}_{H_2} = \pm (\rho^{\partial}_A )_* \hat{\omega}_{H_2}=0$.
\end{proof}

Thus we have a cohomology class $[(\rho_A )_* \hat{\omega}_{H_2}] \in H^{2n-2j-2}_{DR}(\I_j (\R^n ))$.
But in fact $\I_j (\R^n )$ is homotopy equivalent to Stiefel manifold of $j$-frames in $\R^n$ \cite{Rossi_thesis}, and
it is known that, when $n-j$ is odd, its cohomology ring with coefficients in $\R$ is given by
\[
 H^* (\I_j (\R^n );\R ) \cong
 \begin{cases}
  H^* (S^{2n-5}\times S^{2n-9} \times \dots \times S^{2(n-j)+1} \times S^{n-1};\R ) & n \text{ is even}, \\
  H^* (S^{2n-7}\times S^{2n-11} \times \dots \times S^{2(n-j)+1};\R ) & n \text{ is odd}.
 \end{cases}
\]
Hence we can find a form $\mu \in \Omega^{2n-2j-3}_{DR}(\I_j (\R^n ))$ such that $d\mu =(\rho_A )_* \hat{\omega}_{H_2}/6$
(the factor $S^{n-1}$ in the right hand side does not cause any trouble; if $2n-2j-2=n-1$ then $n=2j+1$ and $n$ becomes
odd).

\begin{defn}\label{def_c}
Define $c:=-p_* D^*_A \mu \in \Omega^{2n-3j-3}_{DR}(\emb{n}{j})$ ($p:\emb{n}{j} \times \R^j \to \emb{n}{j}$ is the first
projection).
\end{defn}

Then Theorem \ref{thm_closed} is easily proved;
\begin{alignat*}{2}
 &d(I(H)+c) = p_* D^*_A (\rho_A )_* \hat{\omega}_{H_2}/6 -dp_* D^*_A \mu & \quad & \\
 &\ = p_* D^*_A (\rho_A )_* \hat{\omega}_{H_2}/6 -p_* dD^*_A \mu \pm p^{\partial}_* D^*_A \mu & \quad &
  \text{(Stokes theorem)} \\
 &\ = p_* D^*_A (\rho_A )_* \hat{\omega}_{H_2}/6 -p_* D^*_A d\mu & \quad & (D_A \text{ is constant near }
  \partial C_1 (\R^j ))\\
 &\ =0 & \quad & \text{(by definition of }\mu ).
\end{alignat*}

\begin{rem}\label{rem_mu_symmetric}
We considered involutions $i_{p,q} : \R^n \to \R^n$ in the proof of Theorem \ref{thm_Haefliger_general}.
They induce involutions $i_{p,q}$ on $\I_j (\R^n )$ given by $\iota \mapsto i_{p,q} \circ \iota$.
This lifts to $\hat{i}_{p,q} : \hat{B}_A \to \hat{B}_A$ ($A=V(H_2 )$) defined by
$ \hat{i}_1 (\iota ;(x_1 ,x_2 ,x_3 );y) := (i_{p,q} \circ \iota ;(x_1 ,x_2 ,x_3 );i_{p,q} (y))$.
This has the orientation sign $(-1)^{p-q+1}$ on the fiber, and $\hat{i}^*_1 \hat{\omega}_A =(-1)^{p-q+1}\hat{\omega}_A$
since $i^*_{p,q} \theta_{*4}=(-1)^{p-q+1}\theta_{*4}$ (see the remark after Proposition \ref{prop_additive}).
Hence $i^*_{p,q} (\rho_A )_* \hat{\omega}_A =(\rho_A )_* \hat{\omega}_A$.
This implies that, replacing $\mu$ with $(\mu +i^*_{p,q} \mu )/2$, we may assume that $i^*_{p,q} \mu =\mu$.
Since $i_{p,q}$'s for different $p,q$ commute with each other, we can arrange $\mu$ so that it is preserved by all of these
involutions by repeating the same procedure as above.
\end{rem}

\subsection{Independency on volume forms}\label{subsec_indep_vol}

There is another vanishing result, which is needed in the proof of Propositions \ref{prop_indep_vol},
\ref{prop_indep_vol2} (independency on the choices of volume forms of the map $I$ on cohomology).
Recall that we assume that $g$, $n$ and $j$ are such that the integration map $I$ is a cochain map and that $n-j>2$.

\begin{lem}\label{tilde_I_closed}
The differential form $\tilde{I} (\Gamma ) \in \Omega^*_{DR}(\emb{n}{j} \times [0,1])$ from the proof of Proposition
\ref{prop_indep_vol} is closed if $n-j>2$.
\end{lem}

\begin{proof}
As already done for $I(\Gamma)$ in \S \ref{subsec_principal_cancel}, \S \ref{subsection_hidden_first} and \S
\ref{subsection_hidden_infty}, we must show
\begin{itemize}
\item
 $d_{\emb{n}{j}\times [0,1]}\tilde{I}(\Gamma ) = \tilde{I}(\delta \Gamma )$ modulo the contributions of hidden faces
 of the fibers of $\pi_{\Gamma_i}$, and hence it vanishes since $\Gamma$ is a cocycle, and
\item
 the contributions of hidden faces also vanish.
\end{itemize}
In this section we have proved the vanishing results by using symmetry and dimension counting.
We have to repeat these proofs for $\tilde{I}$.
The symmetry arguments can be applied to the cases here, since the factor $[0,1]$ does not cause any trouble.

We can check that the dimension-counting arguments also work.
The key ingredients in the proofs are the equation (\ref{p_A}) and Lemmas \ref{vanish_1}, \ref{vanish_infty}.
For $\tilde{I}$, we need to replace the criteria in Lemmas \ref{vanish_1}, \ref{vanish_infty} with
\begin{alignat*}{2}
 &\deg \omega_{\Gamma_A} = nb-n \ \text{ or } \ nb-(n+1) &\quad &\text{if } a=0, \\
 &0 \le \deg \omega_{\Gamma_A} - (ja+nb-(j+1)) \le nj+1  &\quad &\text{if } a>0.
\end{alignat*}
and
\[
 \deg \omega_{\Gamma^{\infty}_A} = ja+nb \ \text{ or } \ ja+nb-1
\]
respectively, since these criteria are consequences of the following pullback square
\[
\xymatrix{
  \Int \Sigma_A \times [0,1] \ar[r]^-{\hat{D}_A \times \id} \ar[d]_-{p_A \times \id}
   & \hat{B}_A \times [0,1] \ar[d]^-{\rho_A \times \id} \\
  C^o_{\Gamma / \Gamma_A} \times [0,1] \ar[r]^-{D_A \times \id} & B_A \times [0,1]
 }
\]
and of $\Int \Sigma^{\infty}_A \times [0,1] \approx C_{\Gamma^c_A} \times \hat{B}^{\infty}_A \times [0,1]$.
These replacements do not affect all the arguments in \S \ref{subsection_hidden_first} and
\S \ref{subsection_hidden_infty}, except for Lemma \ref{lem_tree_infty};
a problem may occur when $n-j=2$ and $A=\{ p,q\}$ consists of two i-vertices, since in such a case
$\deg \omega_{\Gamma^{\infty}_A}=(j-1)+(n-1)=2j=\dim (\hat{B}^{\infty}_A \times [0,1])$ can happen.
\end{proof}

Similar arguments show the following.

\begin{lem}\label{lem_indep_H_vol}
Suppose $n-j \ge 3$ is odd.
Then the cohomology class $\calH \in H^{2n-3j-3}_{DR}(\emb{n}{j})$ is independent of the choices of (anti-)symmetric
volume forms.
\end{lem}

\begin{proof}
First let $vol_{S^{j-1}}$ be fixed, and let $v_0$, $v_1$ be two symmetric volume forms of $S^{n-1}$ with total integral
one.
Define the maps $I_0$, $I_1$, $\tilde{I}$ and the form $\tilde{v}$ as in the proof of Proposition \ref{prop_indep_vol}.
Then we have
\begin{itemize}
\item $d_{\emb{n}{j}\times [0,1]}\tilde{I}(H) = \tilde{I}(\delta H)$ modulo the contributions of hidden
faces of the fibers of $\pi_{H_i}$, and hence vanishes since $H$ is a cocycle, and
\item the contributions of hidden faces except for $\Sigma_{V(H_2 )}$ also vanish.
\end{itemize}
Let $A:=V(H_2 )$.
Then $(\rho_A \times \id )_* \hat{\omega}_A \in \Omega^{2n-2j-2}_{DR}(\I_j (\R^n )\times [0,1])$ is a closed form;
the proof is same as Lemma \ref{lem_hatomega_closed}.
But since $H^{2n-2j-2}_{DR}(\I_j (\R^n )\times [0,1])=0$ when $n-j$ is odd, we have a form $\tilde{\mu}$ which satisfies
$(\rho_A \times \id )_* \hat{\omega}_A /6=d\tilde{\mu}$.
Using $\tilde{c}:=-\tilde{p}_* (D_A \times \id )^* \tilde{\mu}$ (where
$\tilde{p}:\R^j \times \emb{n}{j}\times [0,1]\to \emb{n}{j}\times [0,1]$ is the projection), we have a closed form
\[
 \tilde{\calH}:=\tilde{I}(H)+\tilde{c} \in \Omega^{2n-3j-3}_{DR}(\emb{n}{j} \times [0,1]).
\]
By the generalized Stokes theorem,
\[
 (I_1 (H)+\tilde{c}|_{\emb{n}{j}\times \{ 1\}}) - (I_0 (H)+\tilde{c}|_{\emb{n}{j}\times \{ 0\}}) =\pm d p_* \tilde{\calH}.
\]
But $\tilde{c}|_{\emb{n}{j}\times \{ \varepsilon \}}$ ($\varepsilon =0,1$) comes from
$\tilde{\mu}|_{\I_j (\R^n )\times \{ \varepsilon \}}$ which satisfies
$d\tilde{\mu}|_{\I_j (\R^n )\times \{ \varepsilon \}}=(\rho_A \times \id_{\{ \varepsilon\}})_* \hat{\omega}_A /6$.
Hence $\tilde{c}|_{\emb{n}{j}\times \{ \varepsilon \}}$ works as a correction term for $I_{\varepsilon}(H)$.
Thus we see that $\calH$ is independent of the choices of symmetric $vol_{S^{n-1}}$.

Similar arguments work when we fix $vol_{S^{n-1}}$ and use two different $vol_{S^{j-1}}$'s
(in this case we can choose $\tilde{\mu}$ of the form $q^* \mu$ where $q:\I_j (\R^n )\times [0,1] \to \I_j (\R^n )$ is
the projection, because the graph $H_2$ contains no $\eta$-edges and hence we do not use $vol_{S^{j-1}}$ to define the
correction term).
\end{proof}

\providecommand{\bysame}{\leavevmode\hbox to3em{\hrulefill}\thinspace}
\providecommand{\MR}{\relax\ifhmode\unskip\space\fi MR }
\providecommand{\MRhref}[2]{%
  \href{http://www.ams.org/mathscinet-getitem?mr=#1}{#2}
}
\providecommand{\href}[2]{#2}

\end{document}